\input epsf

\documentclass[11pt]{article}
\usepackage{float}
\usepackage{epstopdf}
\usepackage{graphicx}
\usepackage{latexsym,amsmath,amsfonts,amscd,amssymb,subfigure}
\usepackage{graphics}
\usepackage{graphicx}
\usepackage[dvips]{color}
\usepackage{epsfig}
\usepackage{changebar}
\usepackage{bm}

\numberwithin{equation}{section}

\begin{document}

\baselineskip=1.5pc
\vspace{.5in}

\begin{center}

{\Large \bf An absolutely convergent  fixed-point fast sweeping WENO method on triangular meshes for steady state of hyperbolic conservation laws}

\end{center}

\vspace{.15in}

\centerline{
Liang Li\footnote{School of Mathematics and Statistics, Huang Huai University, Zhumadian, Henan 463000, P.R. China. E-mail: liliangnuaa@163.com.},
Jun Zhu\footnote{
State Key Laboratory of Mechanics and Control of Mechanical Structures and Key Laboratory of Mathematical Modelling and High Performance Computing of Air Vehicles (NUAA), MIIT, Nanjing University of Aeronautics and Astronautics, Nanjing, Jiangsu 210016, P.R.
China. E-mail: zhujun@nuaa.edu.cn. Research was supported by NSFC grant 11872210 and Grant No. MCMS-I-0120G01.},
Yong-Tao Zhang\footnote{Department of Applied and Computational
Mathematics and Statistics, University of Notre Dame, Notre Dame,
IN 46556, USA. E-mail: yzhang10@nd.edu
}$\textsuperscript{,}
\renewcommand*{\thefootnote}{\fnsymbol{footnote}}
\setcounter{footnote}{0}\footnote{Corresponding author.}$
}
\baselineskip=1.7pc

\begin{abstract}
High order accuracy fast sweeping methods have been developed in the literature to efficiently solve steady state solutions of hyperbolic partial differential
equations (PDEs) on rectangular meshes, while they were not available yet on unstructured meshes.
In this paper, we extend high order accuracy fast sweeping methods to unstructured triangular meshes by applying fixed-point iterative sweeping techniques to a fifth-order finite volume unstructured WENO scheme, for solving steady state solutions of hyperbolic conservation laws.
Similar as other fast sweeping methods, fixed-point fast sweeping methods use the Gauss-Seidel iterations and alternating sweeping strategy to cover
characteristics of hyperbolic PDEs in each sweeping order to achieve fast convergence
rate to steady state solutions. An advantage of fixed-point fast sweeping methods
which distinguishes them from other fast sweeping methods
is that they are explicit and do not require inverse operation of nonlinear local systems. This also provides certain  convenience in designing high order fast sweeping methods on unstructured meshes. As in the first order fast sweeping methods on triangular meshes, we introduce multiple reference points to determine alternating sweeping directions on unstructured meshes. All the
cells on the mesh are ordered according to their centroids' distances to those reference points, and the resulted orderings provide sweeping directions for iterations.
To make the residue of the fast sweeping iterations converge to machine zero / round off errors, we follow the approach in our early work of developing
the absolutely convergent fixed-point fast sweeping WENO methods on rectangular meshes, and adopt high order WENO scheme
with unequal-sized sub-stencils, specifically here a fifth-order finite volume unstructured WENO scheme for spatial discretization. Extensive numerical experiments, including problems with complex domain geometries, are performed to show the accuracy, computational efficiency, and absolute convergence of the presented fifth-order fast sweeping scheme on triangular meshes. Furthermore, the proposed method is compared with the forward Euler time marching method and the popular third order total variation diminishing Rung-Kutta (TVD-RK3) time-marching method for steady state computations. Numerical examples show that the developed fixed-point fast sweeping WENO method is the most efficient scheme among them, and especially it can save up to $70\%$ CPU time costs than TVD-RK3 in order to converge to steady state solutions.
\end{abstract}

\bigskip
{\bf Key Words:} Fixed-point fast sweeping methods, triangular mesh, high order WENO schemes, unequal-sized sub-stencils, steady state, hyperbolic conservation laws.

\normalsize
\section{Introduction}

Steady state problems of hyperbolic PDEs, such as hyperbolic conservation laws and Hamilton-Jacobi equations, are common mathematical models arising in many applications, e.g. compressible fluid mechanics, optimal control, geometric optics, image processing and computer vision, etc. One of the most important properties of these hyperbolic type boundary value problems is that their solution information propagates along characteristics starting from the boundary. Fast sweeping methods, a class of  iterative methods originally developed in the literature to solve static Hamilton-Jacobi equations (see e.g. \cite{ZHAO, QZZ, QZZ2, FLZ}), take advantage of such property to efficiently solve steady state problems of hyperbolic PDEs.
The methods use alternating sweeping strategy to cover a family of characteristics in a certain direction simultaneously in each sweeping order. Combined with the Gauss-Seidel iterations, these methods can achieve a fast convergence speed. High order accuracy fast sweeping methods have been developed and studied extensively on rectangular meshes. In \cite{ZZQ}, high order weighted essentially non-oscillatory (WENO) fast sweeping schemes for solving static Hamilton-Jacobi equations on rectangular meshes were developed,
where an explicit strategy in the iterative schemes was designed to avoid directly solving very complicated local nonlinear equations derived from high order WENO discretizations.
The methods were combined with accurate boundary treatment techniques such as the inverse Lax-Wendroff methods in \cite{XZZS}. In an implicit way, fast sweeping methods were also applied in discontinuous Galerkin (DG) methods to efficiently solve static Eikonal equations \cite{LSZZ,WZ2,ZCLZS}. This kind of implicit high order fast sweeping methods are very efficient and often have linear computational complexity as the first order fast sweeping methods \cite{ZHAO, QZZ, QZZ2}); however the algorithms are much more complicated than explicit
methods such as the WENO fast sweeping methods \cite{ZZQ}, which makes these implicit high order fast sweeping methods difficult to be applied in solving more complicated hyperbolic PDEs.

Although high order accuracy fast sweeping methods were well developed on rectangular meshes, they have not been explored much on unstructured meshes to solve steady state problems of general hyperbolic PDEs. A key factor
in developing a high order accuracy fast sweeping method on unstructured meshes, which is efficient and relatively easy to be implemented, is to design an approach to effectively incorporate an unstructured nonlinearly stable scheme into the fast sweeping iteration framework to accelerate the convergence to steady state solutions. In this paper, we develop high order accuracy fast sweeping methods on unstructured triangular meshes by applying fixed-point iterative sweeping techniques to a fifth-order finite volume unstructured WENO scheme, for solving steady state solutions of hyperbolic conservation laws. The fixed-point fast sweeping methods were first designed in \cite{ZZC} for solving static Hamilton-Jacobi equations on rectangular meshes.
They can be considered as a generalization of the explicit high order fast sweeping methods in \cite{ZZQ}.
As other fast sweeping methods, the fixed-point fast sweeping methods use the Gauss-Seidel iterations and alternating sweeping strategy to cover
characteristics of hyperbolic PDEs in each sweeping order to achieve fast convergence
rate to steady state solutions.
However, this kind of fast sweeping methods have nice properties such as that they are fully explicit and do not involve any inverse operation of nonlinear local systems, and they can be easily adopted to solve complex hyperbolic systems with any monotone numerical fluxes and high order nonlinearly stable schemes. For example, the fixed-point fast sweeping methods were applied to sparse-grid WENO schemes for efficiently solving multidimensional Eikonal equations in \cite{ZMYZ}; in \cite{WZ, LIZZ, LZSZ},
high order accuracy fixed-point fast sweeping WENO methods for
efficiently solving steady state problems of compressible flows modeled by nonlinear hyperbolic conservation laws were developed. As shown in \cite{LIZZ}, it is interesting
to find that the special case of using the parameter $\gamma=1$ and the Lax-Friedrichs flux in the third order fixed-point fast sweeping schemes  is equivalent to the third order Lax-Friedrichs fast sweeping methods in \cite{ZZQ, W.Chen}.
Here,
the ``explicit'' property of fixed-point fast sweeping methods also provides convenience for us to develop high order fast sweeping methods on unstructured meshes.

As in the fast sweeping methods on rectangular meshes, there are also two important components in designing fast sweeping  methods on unstructured meshes: (1) a nonlinearly stable scheme to discretize the hyperbolic PDEs, which is often called a ``local solver''; (2) systematic orderings of all grid points / cells, which can cover all directions of the characteristics.
The local solver used in this paper is in the class of WENO schemes.
WENO schemes are popular high order accuracy numerical methods for spatial discretizations of nonlinear hyperbolic conservation laws, e.g., see \cite{BCCD,JS,LOC,SZS,S2}. They have the
advantage of attaining uniform high order accuracy in smooth regions of the solution while maintaining sharp and essentially
non-oscillatory transitions of discontinuities. On unstructured meshes, WENO schemes were constructed in e.g.
\cite{DK, DK2, F, HS,D.Levy, TTD, Y.-T.Zhang,ZS}. These unstructured WENO schemes follow the idea in original structured WENO schemes \cite{JS, S}, i.e., the sub-stencils have comparable size and are smaller than the big stencil of the scheme. They can be classified to two types according to their differences in WENO reconstructions on unstructured
meshes, as that defined in \cite{ZS, Y.Liu}. Specifically, the major difference between two types of these schemes
is the different method to construct sub-stencils and find linear weights.
The first type (type I) reconstruction \cite{DK, DK2, F, TTD}  has an order
of accuracy not higher than that of the reconstruction on
each sub-stencil. This is similar as essentially non-oscillatory (ENO) schemes \cite{Harten, SO}.
The nonlinear weights in type I WENO reconstructions do not contribute towards the increase of the order of accuracy, but they are designed purely for the purpose of
nonlinear stability, or to avoid spurious oscillations.
Because type I WENO schemes just need to choose the linear
weights as arbitrary positive numbers for better
linear stability (e.g. the centered sub-stencil is assigned
a larger linear weight than the others), they are easier to construct than the type II WENO
schemes discussed in the following. The second type (type II) consists of WENO schemes whose order
of accuracy is higher than that of the reconstruction on
each sub-stencil (see e.g. \cite{HS,D.Levy, Y.-T.Zhang,ZS}). A crucial step in building a type II WENO scheme on unstructured meshes is to construct lower order polynomials
whose weighted average will give the same result as the high order reconstruction at each Gaussian
quadrature point for the flux integral on the element boundary.
This step is actually the most difficult step in designing a robust type II high order WENO schemes on unstructured meshes,
since we can not guarantee the quality of the unstructured meshes when the domain geometry is very complicated.
Especially, when the spatial domain has higher dimensions (e.g., three dimensional problems)
and complex geometry, the quality of the unstructured meshes is hard to control.
Distorted local mesh geometries can be easily generated.
The local linear system for finding linear weights could have very large condition number or is even singular
at the places where mesh quality is bad. This is the reason why type II WENO schemes are more difficult to
construct than type I WENO schemes on unstructured meshes. However they have a much more compact stencil than type I WENO schemes of the same accuracy, which is an advantage in
applications. In \cite{Y.Liu}, the approaches of type II and type I WENO schemes are combined to achieve a robust unstructured
finite volume WENO reconstruction, and the appearance of negative and very large linear weights are avoided
no matter how bad the quality of the unstructured meshes is. Recently, high order WENO schemes
with unequal-sized sub-stencils were developed on unstructured meshes \cite{ZQ3,ZS3}. As the type I WENO schemes,
the linear weights of these WENO schemes with unequal-sized sub-stencils
can be chosen as arbitrary positive numbers as long as their sum equals 1, hence they are easier to construct than
the type II WENO schemes on unstructured meshes. Also they have the same big stencil as the type II WENO schemes,
hence they are as compact as the type II WENO schemes. Another advantage of WENO schemes
with unequal-sized sub-stencils is on solving steady state problems. Studies on high order WENO schemes on unequal-sized sub-stencils reveal that they improve the convergence of high order WENO
schemes with equal-sized sub-stencils to steady state solutions so that the residue of time-marching iterations settles down to machine zero / round off errors \cite{JUNZ4,ZS3}. The property of an iterative scheme that the residue of its iterations can settle down to machine zero or the round off error level in a finite number of iterations is called ``absolute convergence'' in \cite{LIZZ}. In this paper, to make the residue of the fast sweeping iterations converges to machine zero / round off errors, we follow the approach in our early work of developing
the absolutely convergent fixed-point fast sweeping WENO methods on rectangular meshes \cite{LIZZ}, and adopt high order WENO scheme
with unequal-sized sub-stencils, specifically here a fifth-order finite volume unstructured WENO discretization \cite{ZS3} as the local solver of the fast sweeping method.

The efficient convergence of fast sweeping methods is due to the fact that all
directions of characteristics can be divided into a finite number of groups, and any characteristic
can be decomposed into a finite number of pieces that belong to one of the
above groups, then systematic orderings can be designed to follow the causality of each
group of directions simultaneously.
Orderings of all grid points or cells on a rectangular mesh are natural. For example in a two-dimensional (2-D) problem, all directions of the characteristics can be divided into four groups, namely, up-right, up-left, down-left, and down-right. All grid points or cells are ordered by using their indexes to provide four orderings to cover all these four groups of characteristics. However, such natural orderings do not exist on unstructured meshes or triangular meshes considered here. In this paper, we apply the approach in the first order fast sweeping methods on triangular meshes \cite{QZZ} to the high order fast sweeping WENO method developed here. Specifically, multiple reference points are introduced to determine alternating sweeping directions on unstructured triangular meshes. All the
cells on the mesh are ordered according to their centroids' distances to those reference points, and the resulted orderings provide sweeping directions for the fixed-point iterations.

To compute steady state of hyperbolic conservation laws, the forward Euler time marching method is preferred rather than a Runge-Kutta method, because it is very simple with a single-step and one stage, and
time direction accuracy has no effects on the numerical accuracy of steady state solutions.
However, a high order linear spatial scheme (e.g., a fifth order linear scheme) with the forward Euler time marching is linearly unstable. Nonlinearly stable spatial discretizations such as WENO schemes can help to stabilize the computation by the forward Euler method, but often a very small CFL number and large amount of iterations are required to converge
to steady state, as shown in \cite{WZ}.
Hence, for both linear and nonlinear stability of high order schemes in solving steady state problems of hyperbolic conservation laws, a common approach is to use a total variation diminishing Runge-Kutta method, e.g. the popular third order scheme (TVD-RK3) \cite{GST,SO}, to time-march numerical solution to steady state.
In \cite{LIZZ,WZ}, it was shown that on rectangular meshes, the fast sweeping technique can largely improve the stability of high order WENO schemes with the forward Euler time marching, and the fast sweeping WENO schemes converge to steady state solution much faster than the popular TVD-RK3 time-marching approach. In this paper, we will draw the similar conclusion for the developed absolutely convergent fixed-point fast sweeping WENO method on unstructured triangular meshes. Numerical experiments show that the proposed fast sweeping scheme is more efficient than both the forward Euler time marching method and the TVD-RK3 scheme, and it converges under comparable CFL numbers as the TVD-RK3 scheme for all examples, while the direct forward Euler time marching method needs a smaller CFL number to converge in general.

The organization of the paper is as follows. The detailed description of the new absolutely convergent fixed-point fast sweeping WENO method on unstructured triangular meshes is provided in Section 2. In Section 3, we perform numerical experiments to study the proposed method, and carry out comparisons of different methods.
Extensive numerical examples, including the shock reflection, supersonic flow past a circular cylinder, and the problems of supersonic and subsonic flows past an airfoil which have complex domain geometries, etc., are solved to show the accuracy, computational efficiency, and absolute convergence of the presented fifth-order fast sweeping WENO scheme on unstructured triangular meshes.
Concluding remarks are given in Section 4.

\section{Description of the numerical methods}
\label{sec1}
We consider steady state problems of hyperbolic conservation laws with appropriate boundary conditions. The two dimensional (2D) case has the following general form
\begin{equation}
\textbf{f}(\textbf{u})_x+\textbf{g}(\textbf{u})_y=\textbf{R},
\label{1}
\end{equation}
where ${\bf u}$ is the vector of the unknown conservative variables, $\textbf{f}(\textbf{u})$ and $\textbf{g}(\textbf{u})$ are the vectors of flux functions,
and $\textbf{R}(x,y)$ is the source term.
For example, the steady Euler system of equations in compressible fluid dynamics has
that $\textbf{u} = (\rho, \rho u, \rho v, E)^T$, $\textbf{f}(\textbf{u}) = (\rho u, \rho u^2 + p, \rho u v, u (E + p))^T$, and $\textbf{g}(\textbf{u}) = (\rho v, \rho u v, \rho v^2 +p, v (E + p))^T$.
Here $\rho$ is the density of fluid, $(u,v)^T $ is the velocity vector, $p$ is the pressure, and $E=\frac{p}{\gamma'-1}+\frac{1}{2}\rho(u^2+v^2)$ is total energy where the constant $\gamma'=1.4$ for the case of air.
A spatial discretization of (\ref{1}) by a high order WENO scheme leads to a large nonlinear system of algebraic
equations with the size determined by the number of spatial triangular cells.

\subsection{The local solver}

\begin{figure}%[H]
\centering
\subfigure[]{
\centering
\includegraphics[width=3.5in]{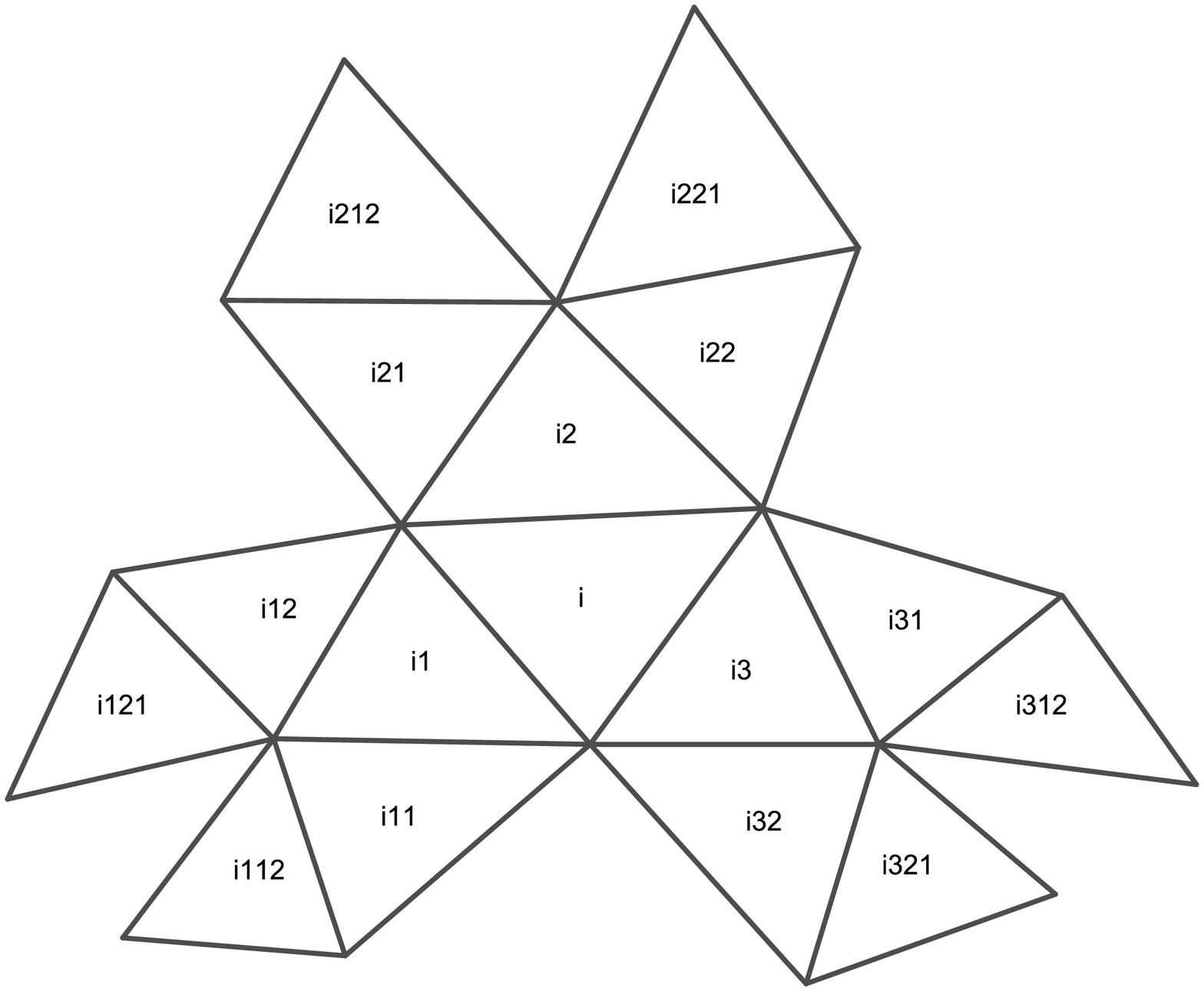}
}

\subfigure[]{
\centering
\includegraphics[width=1.4in]{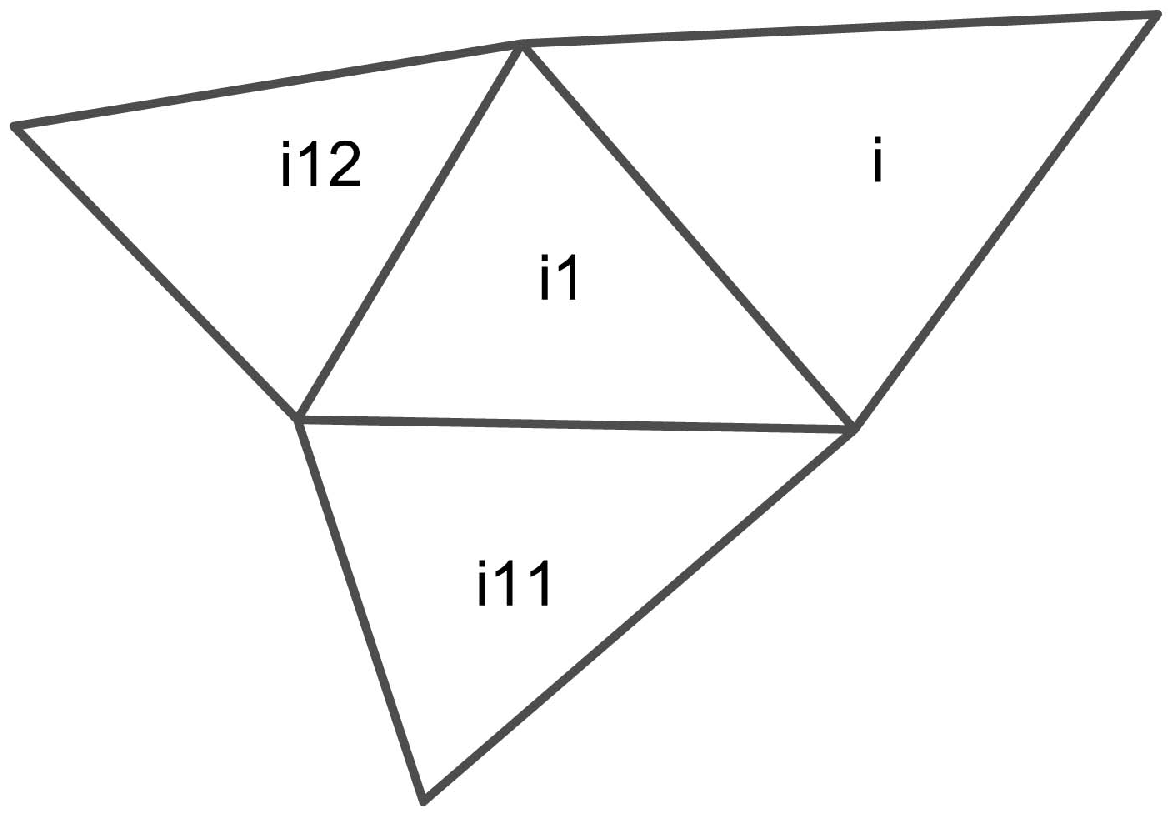}
}
\subfigure[]{
\centering
\includegraphics[width=1.4in]{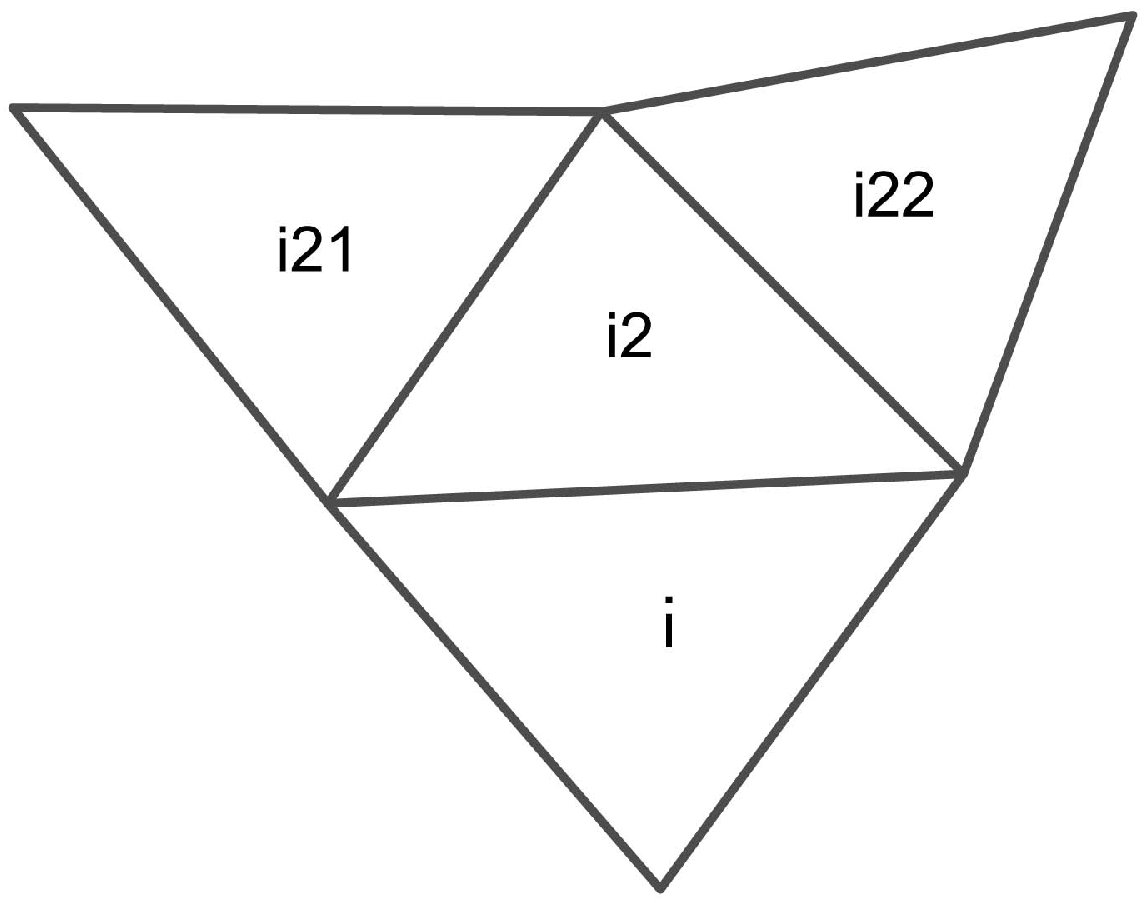}
}
\subfigure[]{
\centering
\includegraphics[width=1.4in]{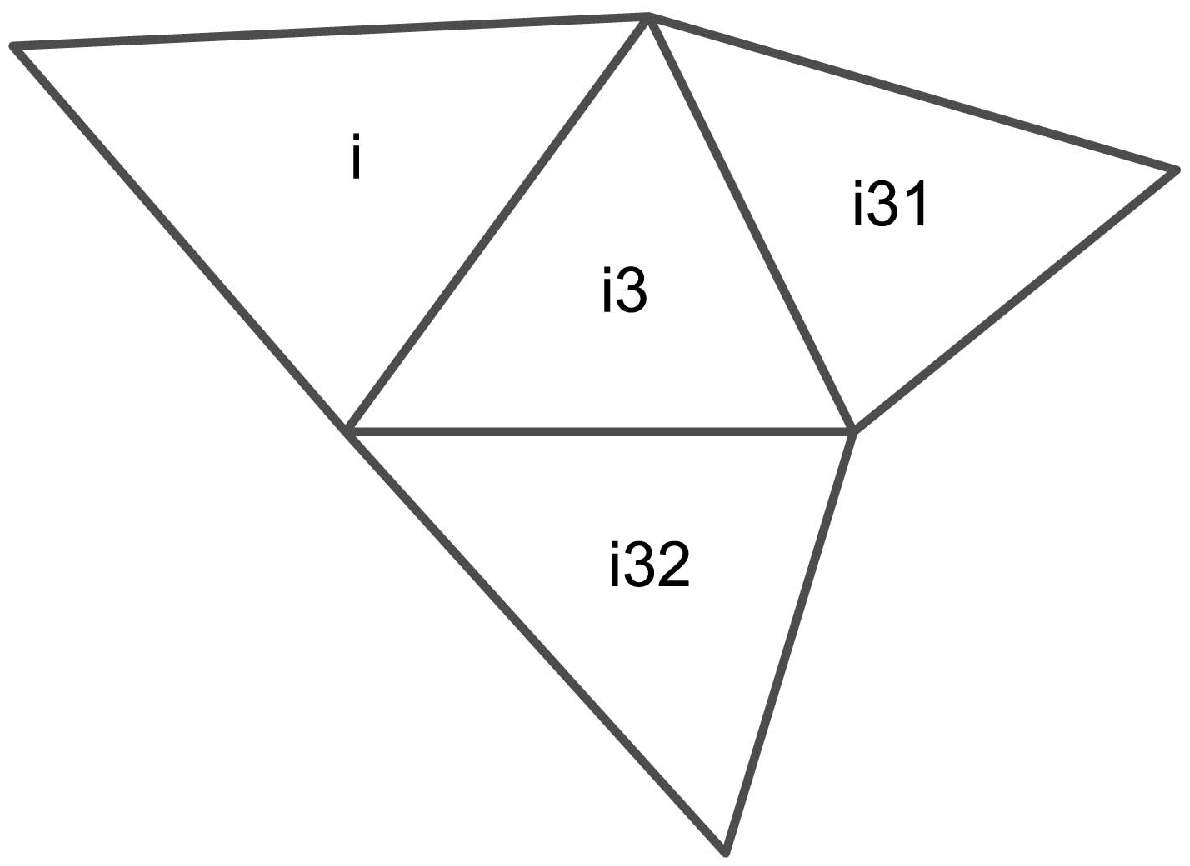}
}
\subfigure[]{
\centering
\includegraphics[width=1.4in]{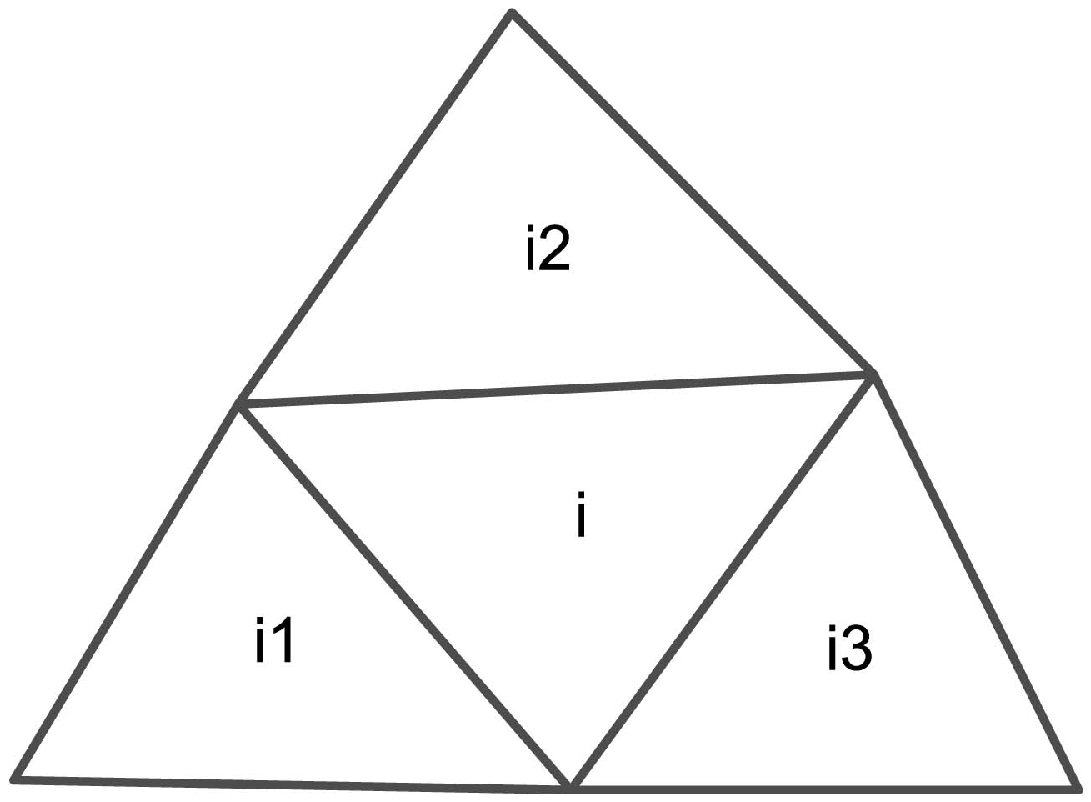}
}
\caption{\label{0}
Unequal-sized sub-stencils of a target cell $\triangle_{i}$ for the fifth order finite volume unstructured WENO local solver. (a) sub-stencil $T_{1}$; (b) sub-stencil $T_{2}$; (c) sub-stencil $T_{3}$, (d) sub-stencil $T_{4}$; (e) sub-stencil $T_{5}$.}
\end{figure}

In this section, we describe the fifth-order finite volume WENO scheme with unequal-sized sub-stencils on triangular meshes for discretization of (\ref{1}), which is the local solver of the developed fast sweeping method.
The computational domain is partitioned by an unstructured triangular mesh with computational cells
$\{\triangle_{i}, i=1, \cdots, M\}$, and $M$ is the total number of cells in the mesh.
We integrate (\ref{1}) over a target cell $\triangle_{i}$ to obtain
\begin{equation}
\frac{1}{|\triangle_{i}|}\int_{\partial\triangle_{i}}\textbf{F}\cdot\overrightarrow{\textbf{n}}ds=\overline{\textbf{R}}_{i},
\label{2}
\end{equation}
in which  $\overline{\textbf{R}}_{i}=\frac{1}{|\triangle_{i}|}\int_{\triangle_{i}}\textbf{R}(x,y)dxdy$, $\textbf{F}=(\textbf{f},\textbf{g})$, $\partial\triangle_{i}$ is the boundary of the target cell $\triangle_{i}, |\triangle_{i}|$ is the area of $\triangle_{i}$, and $\overrightarrow{\textbf{n}}$ is the outward unit normal to the boundary of $\triangle_{i}$. The line integrals in (\ref{2}) are discretized by a three-point Gaussian quadrature formula $\cite{HS}$ on each edge of the cell $\triangle_{i}$
\begin{equation}
\int_{\partial\triangle_{i}}\textbf{F}\cdot\overrightarrow{n}ds\approx\sum_{k=1}^{3}|\partial\triangle_{i_{k}}|\sum_{m=1}^{3}\sigma_{m}\textbf{F}(\textbf{u}(x_{G_{k_{m}}},y_{G_{k_{m}}}))\cdot\overrightarrow{\textbf{n}}_{k},
\label{3}
\end{equation}
where $\{|\partial\triangle_{i_{k}}|\}_{k=1}^3$ are the lengths of the edges, $\{\sigma_{m}\}_{m=1}^3$ are the quadrature weights, $\{(x_{G_{k_{m}}},y_{G_{k_{m}}})\}_{k, m=1}^3$ are the Gaussian quadrature points, $\{\overrightarrow{\textbf{n}}_{k}\}_{k=1}^3$ are the outward unit normals at the quadrature points. For linear stability, $\textbf{F}(\textbf{u}(x_{G_{k_{m}}},y_{G_{k_{m}}}))\cdot\overrightarrow{\textbf{n}}_{k}, m=1, \cdots, 3; k=1, \cdots, 3$ are approximated by monotone numerical fluxes such as the Lax-Friedrichs flux
\begin{equation}
\begin{split}
\textbf{F}(\textbf{u}(x_{G_{k_{m}}},y_{G_{k_{m}}}))\cdot\overrightarrow{\textbf{n}}_{k}&\approx\frac{1}{2}\big[(\textbf{F}(\textbf{u}^{+}(x_{G_{k_{m}}},y_{G_{k_{m}}}))+\textbf{F}(\textbf{u}^{-}(x_{G_{k_{m}}},y_{G_{k_{m}}})))\cdot\overrightarrow{\textbf{n}}_{k}\\
&-\alpha(\textbf{u}^{+}(x_{G_{k_{m}}},y_{G_{k_{m}}})-\textbf{u}^{-}(x_{G_{k_{m}}},y_{G_{k_{m}}}))\big].
\label{4}
\end{split}
\end{equation}
Here $\alpha$ is taken as an upper bound for the eigenvalues of the flux Jacobian in the $\overrightarrow{\textbf{n}}_{k}$ direction, and $\textbf{u}^{-}$ and $\textbf{u}^{+}$ are the reconstructed values of the numerical solution $\textbf{u}$ inside and outside of the target cell $\triangle_{i}$ at different Gaussian quadrature points, based on the cell average values of $\textbf{u}$ on cells of the computational stencil.
We use the fifth-order WENO method \cite{ZS3} to reconstruct the function values of $\textbf{u}(x,y)$ at different Gaussian quadrature points $(x_{G_{k_{m}}},y_{G_{k_{m}}})$ on the cell boundaries. Reconstruction on the target cell $\triangle_{i}$ provides values for $\textbf{u}^{-}$ in (\ref{4}), with detailed algorithm given in the following.
Values for $\textbf{u}^{+}$ in (\ref{4}) are provided by the reconstruction on one of the neighboring cells of $\triangle_{i}$.

\bigskip
{\bf Reconstruction algorithm:}

{\bf Step 1.} Select a big central spatial stencil $T_{1}=\big\{\triangle_{i},\triangle_{i1},\triangle_{i2},\triangle_{i3},\triangle_{i11},\triangle_{i12},\triangle_{i21},\triangle_{i22},\triangle_{i31},\\ \triangle_{i32}, \triangle_{i112},\triangle_{i121},\triangle_{i212},\triangle_{i221},\triangle_{i312},\triangle_{i321}\}$  (see Figure \ref{0}), which are formed by several layers of neighbor cells of the target cell $\triangle_{i}$. Based on the cell average values of ${\bf u}$ on cells of $T_{1}$, we construct a quartic polynomial ${\bf p}_{1}(x,y)$ on $T_{1}$ to obtain a fifth-order approximation of the conservative variable ${\bf u}$. Since in general the number of cells on $T_{1}$ are more than the degree of freedom of the quartic polynomial, the least-square method is used to find ${\bf p}_{1}(x,y)$ as in \cite{HS}, namely, it is required that ${\bf p}_{1}(x,y)$ has the same cell average of ${\bf u}$ on the target cell $\triangle_{i}$ and matches the cell averages of ${\bf u}$ on the other cells of the stencil $T_{1}$ in the least-square way:
\begin{equation}
\begin{split}
&\frac{1}{|\triangle_{i}|}\int_{\triangle_{i}}\textbf{p}_{1}(x,y)dxdy=\overline{\textbf{u}}_{i},\quad \textbf{p}_{1}(x,y)=\arg\min\sum_{ll\in{W_1}}\left(\frac{1}{|\triangle_{ll}|}\int_{\triangle_{ll}}\textbf{p}(x,y)dxdy-\overline{\textbf{u}}_{ll}\right)^{2},\\
&W_1=\{i1,i2,i3,i11,i12,i21,i22,i31,i32,i112,i121,i212,i221,i312,i321\}.
\label{5}
\end{split}
\end{equation}
Note that if the number of cells located inside the stencil $T_1$ is less than that is required for reconstructing the polynomial because some of the cells coincide with each other, neighboring cells in the next layer of the mesh will be added to the stencil to provide enough number of cells for the reconstruction.

{\bf Step 2.} Form three sectorial sub-stencils $T_{2}=\{\triangle_{i},\triangle_{i1},\triangle_{i11},\triangle_{i12}\},T_{3}=\{\triangle_{i},\triangle_{i2},\triangle_{i21},\triangle_{i22}\},\\ T_{4}=\{\triangle_{i},\triangle_{i3},\triangle_{i31},\triangle_{i32}\}$ and one central sub-stencil $T_{5}=\{\triangle_{i},\triangle_{i1},\triangle_{i2},\triangle_{i3}\}$ (see Figure \ref{0}).
Note that in forming these sectorial sub-stencils $T_2,T_3,T_4$, we use straight lines to connect the barycenter of the target cell $\triangle_i$ with its three vertices and split the plane into three sectors. Every sectorial sub-stencil consists of the target cell $\triangle_i$ and its neighboring cells whose barycenters lie in the same sector. Such kind of distributions of sub-stencils are important to obtain the non-oscillatory performance for non-smooth solutions, as shown in \cite{Y.Liu}.
Then four linear polynomials ${\bf p}_{m}(x,y), m=2, \cdots, 5$ on such sub-stencils are constructed to obtain  second-order approximations of the conservative variable ${\bf u}$. Again, the least-square method is used to find ${\bf p}_{m}(x,y)$ by requiring that ${\bf p}_{m}(x,y)$ has the same cell average of ${\bf u}$ on the target cell $\triangle_{i}$ and matches the cell averages of ${\bf u}$ on the other cells of the sub-stencil $T_{m}$ in the least-square way:
\begin{equation}
\begin{split}
&\frac{1}{|\triangle_{i}|}\int_{\triangle_{i}}\textbf{p}_{m}(x,y)dxdy=\overline{\textbf{u}}_{i},\quad \textbf{p}_{m}(x,y)=\arg\min\sum_{ll\in{W_{m}}}\left(\frac{1}{|\triangle_{ll}|}\int_{\triangle_{ll}}\textbf{p}(x,y)dxdy-\overline{\textbf{u}}_{ll}\right)^{2},\\
&m=2,\cdots,5; W_{2}=\{i1,i11,i12\}, W_{3}=\{i2,i21,i22\}, W_{4}=\{i3,i31,i32\}, W_{5}=\{i1,i2,i3\}.
\label{6}
\end{split}
\end{equation}

{\bf Step 3.} Define the linear weights. We emphasize that different from the previous unstructured WENO schemes \cite{HS, ZS}, the WENO reconstruction with unequal-sized sub-stencils here uses the big stencil $T_1$ as one of the sub-stencils and the other sub-stencils $T_2,T_3,T_4,T_5$ are much smaller, while in the WENO schemes \cite{HS, ZS}, every sub-stencil is only part of the big stencil and their union is the big stencil.
To use the reconstruction $\textbf{p}_{1}(x,y)$ on the sub-stencil $T_1$ for the final high order accuracy, we apply the similar ideas for the central WENO schemes \cite{LPR,LPR2} and rewrite $\textbf{p}_{1}(x,y)$ as
\begin{equation}
\textbf{p}_{1}(x,y)=\gamma_{1}\bigg(\frac{1}{\gamma_{1}}\textbf{p}_{1}(x,y)-\sum_{m=2}^{5}\frac{\gamma_{m}}{\gamma_{1}}\textbf{p}_{m}(x,y)\bigg)+\sum_{m=2}^{5}{\gamma_{m}}\textbf{p}_{m}(x,y)
\label{7}
\end{equation}
with $\sum_{m=1}^{5}\gamma_{m}=1$ and $\gamma_{1}\neq{0}$. In these expressions, $\gamma_{m}$ for $m=1, \cdots,5$ are the linear weights. Following the practice in \cite{DK,Y.Liu,ZQ,ZZSQ}, we take the linear weights as $\gamma_{1}=0.96$ and $\gamma_{2}=\gamma_{3}=\gamma_{4}=\gamma_{5}=0.01$.

{\bf Step 4.} Compute the smoothness indicators $\beta_{k}$, which measure how smooth the functions $\textbf{p}_{k}$ for $k=1,\cdots,5$ are in the target cell $\triangle_{i}$. We use the same recipe for the smoothness indicators as that in \cite{HS,JS}:
\begin{equation}
\beta_{k}=\sum_{|\alpha|=1}^{r}\int_{\triangle_{i}}|\triangle_{i}|^{|\alpha|-1}(\textbf{D}^{\alpha}\textbf{p}_{k}(x,y))^{2}dxdy, \quad k=1,\cdots,5,
\label{8}
\end{equation}
where $\textbf{D}$ is the differential operator and $\alpha$ is the multi-index for derivative orders.  $r=4$ for $l=1$, and $r=1$ for $l=2,\cdots,5$.

{\bf Step 5.} Compute the nonlinear weights based on the linear weights and the smoothness indicators. We adopt the WENO-Z type nonlinear weights as specified in \cite{BCCD,CCD,DB}, the nonlinear weights are defined as
\begin{equation}
\omega_{k}=\frac{\overline{\omega}_{k}}{\sum_{k=1}^{5}\overline{\omega}_{k}},\quad \overline{\omega}_{k}=\gamma_{k}(1+\frac{\tau}{\varepsilon+\beta_{k}}),\quad k=1,\cdots,5.
\label{9}
\end{equation}
According to \cite{ZQ,ZQ3}, we define $\tau=(\frac{|\beta_{1}-\beta_{2}|+|\beta_{1}-\beta_{3}|+|\beta_{1}-\beta_{4}|+|\beta_{1}-\beta_{5}|}{4})^{2}$, and take $\varepsilon=10^{-6}$.

{\bf Step 6.} The final reconstruction polynomial for the approximation of $u(x,y)$ at any points of the target cell $\triangle_{i}$ is given as
\begin{equation}
\textbf{u}(x,y)\approx\omega_{1}\bigg(\frac{1}{\gamma_{1}}\textbf{p}_{1}(x,y)-\sum_{m=2}^{5}\frac{\gamma_{m}}{\gamma_{1}}\textbf{p}_{m}(x,y)\bigg)+\sum_{m=2}^{5}{\omega_{m}}\textbf{p}_{m}(x,y).
\label{10}
\end{equation}

\subsection{Absolutely convergent fixed-point fast sweeping WENO scheme}
After the spatial discretization of the equation (\ref{1}) by the finite volume scheme described in the last section, we have
\begin{equation}
\begin{split}
0&=-\frac{1}{|\triangle_{i}|}\sum_{k=1}^{3}|\partial\triangle_{i_{k}}|\sum_{m=1}^{3}\sigma_{m}\frac{1}{2}\big[(\textbf{F}(\textbf{u}^{+}(x_{G_{k_{m}}},y_{G_{k_{m}}}))+\textbf{F}(\textbf{u}^{-}(x_{G_{k_{m}}},y_{G_{k_{m}}})))\cdot\overrightarrow{\textbf{n}}_{k}\\
&-\alpha(\textbf{u}^{+}(x_{G_{k_{m}}},y_{G_{k_{m}}})-\textbf{u}^{-}(x_{G_{k_{m}}},y_{G_{k_{m}}}))\big]+\overline{\textbf{R}}_{i}, \quad i=1, \cdots, M,
\label{11}
\end{split}
\end{equation}
where $M$ is the number of triangular cells. The aforementioned fifth-order WENO reconstructions at Gaussian quadrature points $(x_{G_{k_{m}}},y_{G_{k_{m}}})$ on the cell boundaries lead to a nonlinear algebraic system
\begin{equation}
L(\overline{\textbf{u}}_{i},\overline{\textbf{u}}_{i1},\cdots,\overline{\textbf{u}}_{i321})=0,\quad i=1, \cdots, M.
\label{11.5}
\end{equation}
Note that here $\overline{\textbf{u}}_{i}$ is the cell average of the numerical solution of the unknown function ${\bf u}$ on the target cell $\triangle_{i}$. $\overline{\textbf{u}}_{i1},\cdots,\overline{\textbf{u}}_{i321}$ are the cell averages on the WENO scheme's stencil of the target cell $\triangle_{i}$. $L$ denotes the spatial discretization operator applied to the PDE using this fifth-order finite volume WENO local solver, which is a nonlinear function of the cell averages of the numerical solution on the computational stencil of the WENO scheme.
The fixed-point fast sweeping schemes for solving the nonlinear system (\ref{11.5}) are based on iterative schemes of time marching type for solving steady state problems. Time marching methods for solving steady state problems are essentially Jacobi type fixed-point iterative schemes. The forward Euler (FE) time marching method with time step size $\Delta t_n$ to solve the nonlinear system (\ref{11.5}) is the following Jacobi type fixed-point iterative scheme
\begin{equation}
\overline{\textbf{u}}_{i}^{n+1}=\overline{\textbf{u}}_{i}^{n}+\triangle{t_{n}}L(\overline{\textbf{u}}_{i}^{n},\overline{\textbf{u}}_{i1}^{n},\cdots,\overline{\textbf{u}}_{i321}^{n}),\quad i=1, \cdots, M,
\label{12}
\end{equation}
where $\overline{\textbf{u}}_{i}^{n+1}$ and $\overline{\textbf{u}}_{i}^{n}$ etc. are the numerical solution values
at iteration steps $n+1$ and $n$.
The popular TVD-RK3 time marching method to solve the nonlinear system (\ref{11.5}) is the following Jacobi type fixed-point iterative scheme:
\begin{equation}
\begin{split}
&\overline{\textbf{u}}_{i}^{(1)}=\overline{\textbf{u}}_{i}^{n}+\triangle{t_{n}}L(\overline{\textbf{u}}_{i}^{n},\overline{\textbf{u}}_{i1}^{n},\cdots,\overline{\textbf{u}}_{i321}^{n}),\quad i=1,\cdots,M,\\
&\overline{\textbf{u}}_{i}^{(2)}=\frac{3}{4}\overline{\textbf{u}}_{i}^{n}+\frac{1}{4}\overline{\textbf{u}}_{i}^{(1)}+\frac{1}{4}\triangle{t_{n}}L(\overline{\textbf{u}}_{i}^{(1)},\overline{\textbf{u}}_{i1}^{(1)},\cdots,\overline{u}_{i321}^{(1)}),\quad i=1,\cdots,M,\\
&\overline{\textbf{u}}_{i}^{n+1}=\frac{1}{3}\overline{\textbf{u}}_{i}^{n}+\frac{2}{3}\overline{\textbf{u}}_{i}^{(2)}+\frac{2}{3}\triangle{t_{n}}L(\overline{\textbf{u}}_{i}^{(2)},\overline{\textbf{u}}_{i1}^{(2)},\cdots,\overline{\textbf{u}}_{i321}^{(2)}),\quad i=1,\cdots,M.
\label{13}
\end{split}
\end{equation}

To obtain faster convergence to steady state solutions of high order unstructured WENO schemes for solving hyperbolic PDEs than the above Jacobi type fixed-point iterations, we apply the fast sweeping techniques so that
the important characteristics property of hyperbolic PDEs can be utilized in the iterations, as in the WENO
fast sweeping methods on structured meshes \cite{WZ, LIZZ}. The fast sweeping techniques can be applied to either the FE scheme (\ref{12}) or the TVD-RK3 scheme (\ref{13}). However, as that found in \cite{WZ,ZMYZ}, the fast sweeping scheme resulted by applying the fast sweeping techniques to
the TVD-RK3 scheme is less efficient than the obtained method by applying the fast sweeping techniques to the FE scheme. Hence here we apply the fast sweeping techniques to the FE scheme (\ref{12}) and obtain the high order FE type fixed-point fast sweeping scheme on unstructured meshes. As shown in the numerical experiments of the following section, the FE fixed-point fast sweeping scheme permits larger CFL numbers than the original FE scheme. In fact it has comparable CFL numbers as that of the TVD-RK3 scheme, and it is much more efficient than the TVD-RK3 scheme to reach steady state of the solutions.
The form of the FE fixed-point fast sweeping scheme is:
\begin{equation}
\overline{\textbf{u}}_{i}^{n+1}=\overline{\textbf{u}}_{i}^{n}+\triangle{t_{n}}L(\overline{\textbf{u}}_{i}^{n},\overline{\textbf{u}}_{i1}^{\ast},\cdots,\overline{\textbf{u}}_{i321}^{\ast}),\quad i=k_{1},\cdots,k_{M}.
\label{14}
\end{equation}
Note that the fast sweeping methods have two essential components, i.e., the Gauss-Seidel philosophy and alternating direction sweeping iterations. The Gauss-Seidel philosophy requires that the newest numerical values of $\overline{\textbf{u}}$ are used in the finite volume WENO reconstruction stencils as long as they are available. Alternating direction sweepings cover characteristics in different directions.
Here the iteration direction in the FE fixed-point fast sweeping scheme (\ref{14}), which is marked as ``$i=k_{1},\cdots,k_{M}$'', means that the iterations in the scheme (\ref{14}) do not just proceed in only one direction ``$i = 1,\cdots,M$'' as that in the Jacobi type
schemes (\ref{12}) and (\ref{13}), but in the alternating directions repeatedly.

\begin{figure}%[H]
\centering
\subfigure[]
{\centering\includegraphics[width=3.4in]{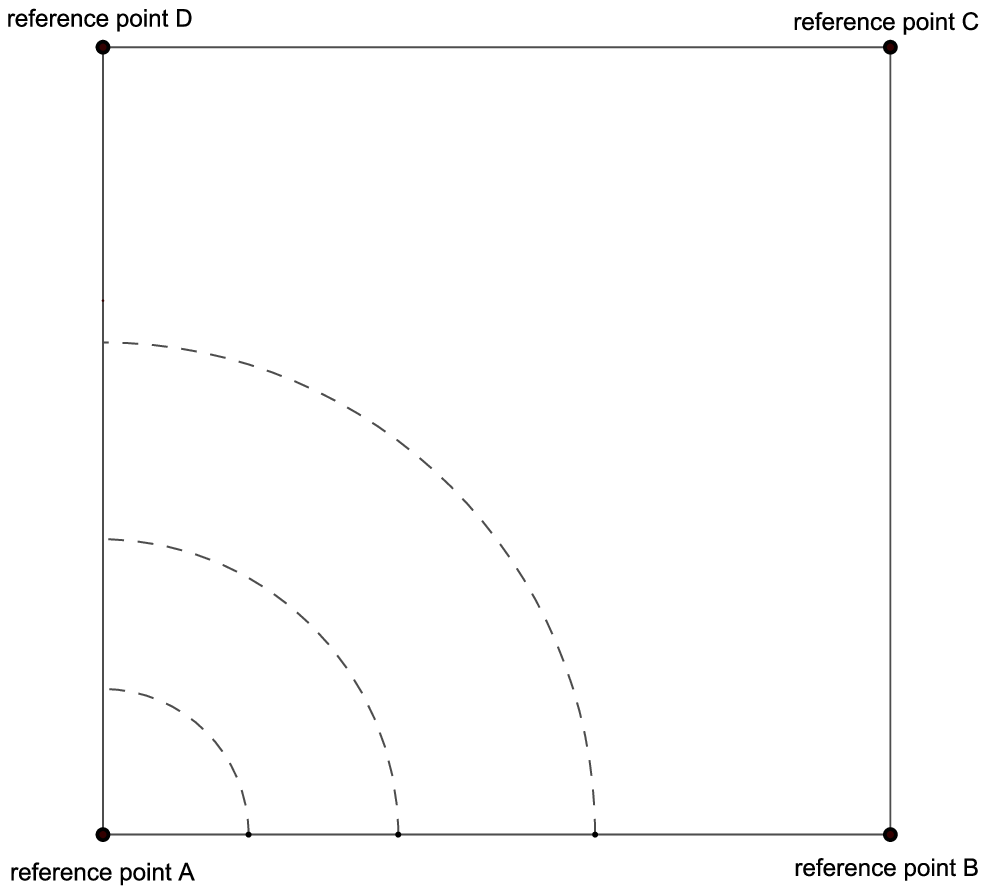}}
\subfigure[]
{\centering\includegraphics[width=2.4in]{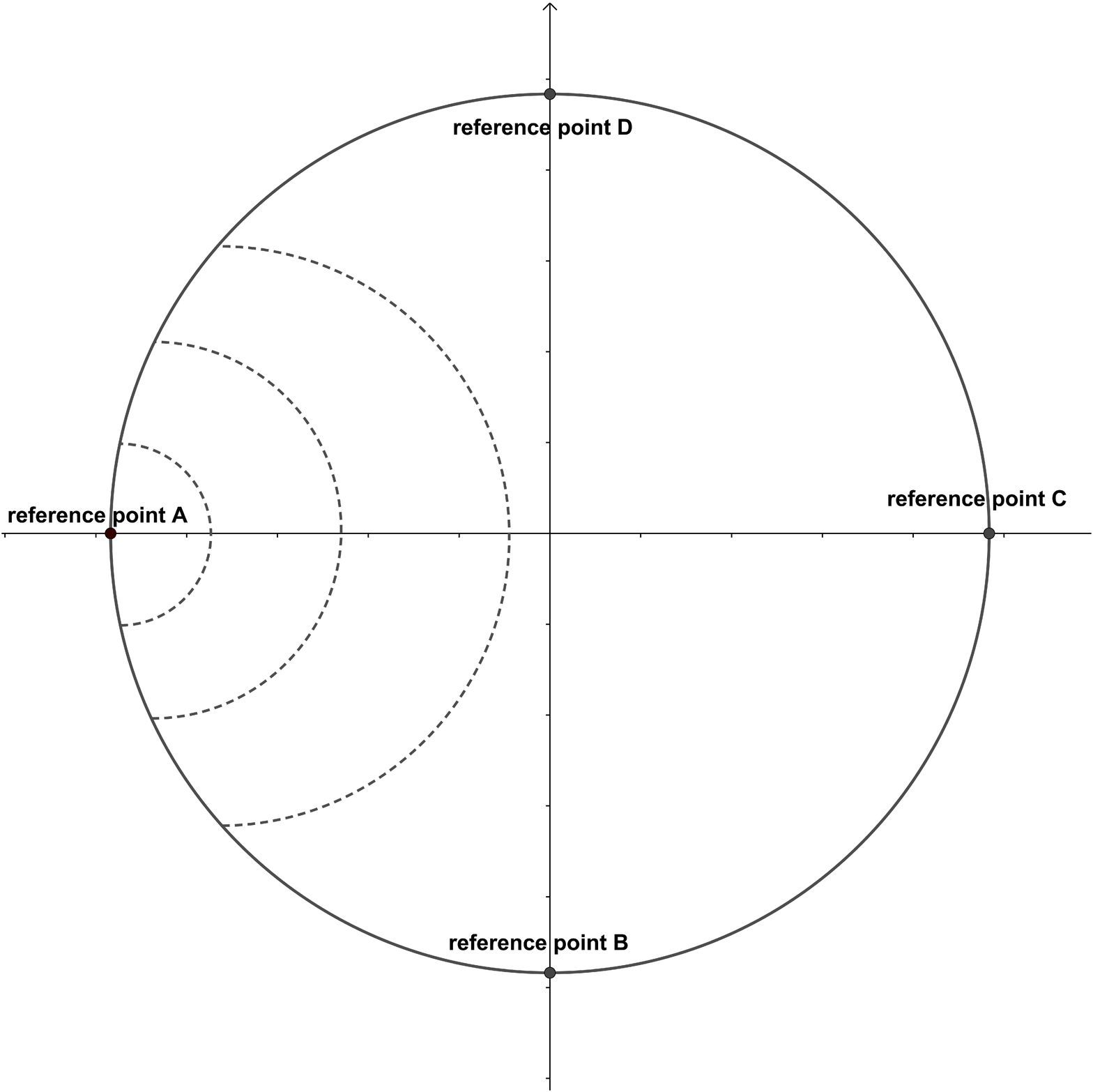}}
\caption{\label{333} Examples of reference points for: (a) a rectangular computational domain; (b) a circular computational domain. A sweeping direction based on the reference point A is shown.}
\end{figure}

Design of alternating sweeping directions on a rectangular mesh is natural. However, it is not straightforward on unstructured meshes. We use the method in \cite{QZZ} for the first order fast sweeping scheme to form the sweeping directions of the high order scheme here, and verify that such alternating sweeping directions are still effective for the high order scheme. The alternating sweeping directions on unstructured meshes are determined by ordering all the cells on the mesh according to their centroids' distances to some reference points on the domain.
As discussed in \cite{QZZ}, the criterion for an optimal choice of reference points and their locations on an unstructured mesh is that all directions of characteristics should be covered with minimal redundancy. In practice, it is better if these reference points are evenly spaced both spatially and angularly in the computational domain.  In the first order fast sweeping method for Eikonal equations on triangular meshes \cite{QZZ}, at least three non-collinear reference points are needed if $l^2-$ distance is used, and the sweepings are performed through all the grid points according to their distances to these reference points in ascent and descent orderings to cover all directions of information propagating along characteristics. For the high order WENO fast sweeping method for hyperbolic conservation laws here, we use four reference points which are evenly (or roughly evenly) distributed on the domain boundary, and numerical experiments show that the resulted alternating sweeping directions are effective for fast convergence of the scheme. For example, we use four corners as reference points if the computational domain is rectangular. See Figure \ref{333} for an illustration of reference points for a rectangular computational domain and a circular computational domain, and sweeping directions based on one of the reference points.

For each reference point $R_{l}$, $l=1,\cdots,4$, we sort all cells on the unstructured mesh according to their centroids' distances to the reference point $R_{l}$ in ascent and descent orders, and put their index numbers into arrays, e.g., the array $S_{l}^{+}(i), i=1,2,\cdots,M$ to store cell index numbers in the resulted ascent order and the array $S_{l}^{-}(i), i=1,2,\cdots,M$ to store cell index numbers in the resulted descent order. With four reference points, we obtain eight sweeping orders. Then in the proposed fixed-point WENO fast sweeping scheme
(\ref{14}), we alternatively take $k_i=S_{l}^{+}(i), i=1,2,\cdots,M$ and $k_i=S_{l}^{-}(i), i=1,2,\cdots,M$,
for $l=1,\cdots,4$, as the sweeping iteration directions. Note that these initial orderings of cells only need to be performed for a fixed mesh once and for all, hence their cost is just a small part of the whole simulation cost. In the numerical tests of this paper, we use a comparison-based sorting method \cite{CLR}.
The complete algorithm is summarized as the following.

{\bf Step 1.} Choose four reference points: $R_{l}, l=1,\cdots,4.$

{\bf Step 2.} For each reference point $R_{l}, l=1,\cdots,4$,  all triangle cells are sorted  according to the $l^{2}-$distances of their centroids to the reference points $R_{l}$ in ascent and descent orders, and store their index numbers in arrays:
$$S_{l}^{+}(i): \text {ascent order}, i=1,2,\cdots,M;$$
$$S_{l}^{-}(i): \text {descent order}, i=1,2,\cdots,M.$$

{\bf Step 3.} Gauss-Seidel iterations using the fast sweeping scheme (\ref{14}) are performed for $n=1, 2, \cdots$, according to the sweeping iteration directions obtained above:
\begin{eqnarray*}
\text{do}& l=1,\cdots, 4 &\\
  & \text{do}&  i=1, \cdots, M \\
& &k_i=S_{l}^{+}(i); \\
& &\text {Apply (\ref{14}) and update}\mbox{ } \overline{\textbf{u}}_{k_i}^{n+1} \mbox{ }\text{on the triangle}\mbox{ } k_i. \\
&\text{enddo} &\\
& &\text {Perform convergence check.} \\
& \text{do}& i=1,\cdots, M \\
& &k_i=S_{l}^{-}(i); \\
& &\text {Apply (\ref{14}) and update}\mbox{ } \overline{\textbf{u}}_{k_i}^{n+1} \mbox{ }\text{on the triangle}\mbox{ } k_i. \\
&\text{enddo} & \\
& &\text {Perform convergence check.} \\
\text{enddo} & &
\end{eqnarray*}

Note that the convergence check is performed via first computing the iteration residue $ResA$ defined in the following section, then judging whether the convergence criterion $ResA \leq \delta$ for $\delta > 0$ given is satisfied. If it is satisfied, we stop the Gauss-Seidel iterations. Since the strategy of alternating direction sweepings utilizes the characteristics property of hyperbolic PDEs, combining with the Gauss-Seidel philosophy, we observe the acceleration of convergence to steady state solutions as shown in the following numerical experiments.
By the Gauss-Seidel philosophy, we use the newest numerical values on the computational stencil of the WENO scheme whenever they are available in the current iteration step. This is why we use the notation
$\overline{\textbf{u}}_{il}^{\ast}, l=1, 2, 3, 11, \cdots, 321$, to represent the numerical values in the scheme
(\ref{14}), and $\overline{\textbf{u}}_{il}^{\ast}$ could be the value $\overline{\textbf{u}}_{il}^{n}$ in the previous iteration step $n$, or the new value $\overline{\textbf{u}}_{il}^{n+1}$ which has been updated and available in the current iteration step, depending on the current sweeping direction of the iteration. To guarantee that the fixed-point iteration (\ref{14}) is a contractive mapping and converges, suitable values of $\triangle{t_{n}}$ which depends on the CFL number, need to be taken. This is similar to choose a suitable CFL number for stability when a high order WENO scheme is used to solve time-dependent hyperbolic PDEs.

\section{Numerical experiments}

In this section, we carry out numerical experiments to test the proposed absolutely convergent fixed-point fast sweeping WENO method on triangular meshes for solving some benchmark steady-state problems of hyperbolic conservation laws. Computational efficiency of the fast sweeping scheme and the other two time marching schemes
is compared. For the convenience of presentation, we call the proposed absolutely convergent fixed-point fast sweeping WENO scheme (\ref{14}) ``FE fast
sweeping scheme'', the forward Euler time marching scheme (\ref{12}) ``FE Jacobi scheme'', and the TVD-RK3 time marching scheme (\ref{13}) ``RK Jacobi scheme''.
Mesh refinement studies are performed to
 compute $L_{1}$ and $L_{\infty}$ numerical errors and accuracy orders of these iterative schemes.  Iteration numbers and CPU times for each iterative scheme to converge are recorded
to compare their computational efficiency.
The convergence of the iterations is measured by the average residue which is defined as
\begin{equation}
ResA=\sum_{i=1}^{M}\frac{|R1_{i}|+|R2_{i}|+|R3_{i}|+|R4_{i}|}{4\times{M}}.
\end{equation}
Here $R\ast_{i}$'s are the local residuals of the conservative variables, namely, $R1_{i}=\frac{\overline{\rho}_{i}^{n+1}-\overline{\rho}_{i}^{n}}{\Delta t_n}, R2_{i}=\frac{(\overline{\rho{u}})_{i}^{n+1}-(\overline{\rho{u}})_{i}^{n}}{\Delta t_n}, R3_{i}=\frac{(\overline{\rho{v}})_{i}^{n+1}-(\overline{\rho{v}})_{i}^{n}}{\Delta t_n}, R4_{i}=\frac{\overline{E}_{i}^{n+1}-\overline{E}_{i}^{n}}{\Delta t_n}$. $M$ is the total number of triangular cells.
The convergence criterion is set to be $ResA \leq \delta$, where the threshold value $\delta$ is taken to be at round off error level $10^{-12}\sim 10^{-11} $.  $\triangle{t_{n}}$ depends on the CFL number as the following
\begin{equation}
\triangle{t_{n}}=\frac{CFL}{\max\limits_{1\leq{i}\leq{M}}\left(\frac{\sum_{ll=1}^{3}\left[|(\overline{u}_{i}^n,\overline{v}_{i}^n)\cdot{\overrightarrow{\mathbf{n}}_{i_{ll}}}|+c_{i}^n\right]|\partial{\triangle_{i_{ll}}}|}{2|\triangle_{i}|}\right)},
\label{15}
\end{equation}
where $c_{i}^n=\sqrt{\frac{\gamma' \overline{p}_{i}^n}{\overline{\rho}_{i}^n}}$ is the sound speed with $\gamma' = 1.4$. In this paper, the number of iterations reported in the numerical simulations counts a complete update of
numerical values in all cells once as one iteration.
Hence these eight alternating directions in the FE fast sweeping scheme are counted as eight iterations.
In the numerical experiments, for various examples we compare the computational efficiency of different schemes by using the largest possible CFL numbers  which lead to iteration convergence with the fastest speed for each method. To identify the largest possible CFL number for a problem, we gradually increase / decrease the values of CFL number from an initial value.

\bigskip
\noindent We first test these schemes by solving two problems with a smooth solution to verify their accuracy, and compare their efficiency to converge to a steady-state solution.

\bigskip
\noindent\textbf{Example 1. An Euler system of equations with source terms}

\begin{figure}%[H]
\centering
\includegraphics[width=3.in]{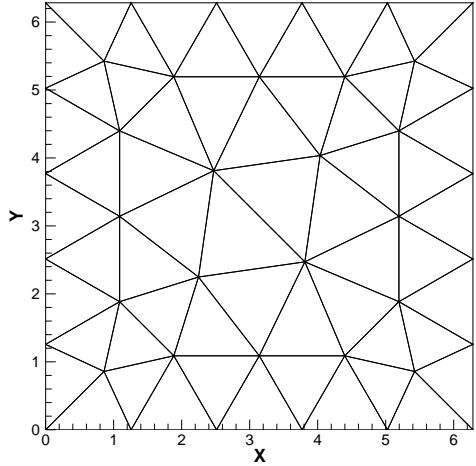}
\caption{The coarsest unstructured triangular mesh for Example 1 and Example 2.}
\label{mesh1}
\end{figure}

\noindent In this example, we solve for  steady-state solution of the following two dimensional Euler system of equations with source terms
\begin{equation}
\frac{\partial}{\partial{t}}\left(\begin{array}{cccc}
 \rho\\ \rho{u}\\ \rho{v} \\ E\\
\end{array}\right)
+\frac{\partial}{\partial{x}}\left(\begin{array}{cccc}
 \rho{u}\\ \rho{u^{2}+p}\\ \rho{uv} \\ u(E+p)\\
\end{array}\right)
+\frac{\partial}{\partial{y}}\left(\begin{array}{cccc}
 \rho{v}\\ \rho{uv}\\ \rho{v^{2}+p} \\ v(E+p)\\
\end{array}\right)
=\left(\begin{array}{cccc}
 0.4\cos(x+y)\\ 0.6\cos(x+y)\\ 0.6\cos(x+y)\\ 1.8\cos(x+y)\\
\end{array}\right).
\end{equation}
The exact steady-state solutions of the problem are $\rho(x,y,\infty)=1 + 0.2\sin(x+y), u(x,y,\infty)=1, v(x,y,\infty)=1, p(x,y,\infty)=1 + 0.2\sin(x+y)$. The computation domain is $(x,y)\in[0,2\pi]\times[0,2\pi]$, and
the exact steady-state solutions are applied
on the domain boundaries. Three different iterative methods, i.e., the FE Jacobi scheme, the RK Jacobi scheme, and  the FE fast sweeping scheme, are used to solve the problem on successively refined unstructured triangular meshes.
The coarsest mesh used here is shown in Figure \ref{mesh1}. In the mesh refinement study to test the accuracy orders of the schemes, the mesh is refined by cutting each triangle in the coarse mesh into four smaller similar ones. As discussed in Section 2, we take four corners of the domain, i.e., $(0,0),(0,2\pi),(2\pi,0)$ and $(2\pi,2\pi)$, as the reference points to form the alternating sweeping directions in the FE fast sweeping scheme.  To start the iterations for the schemes, we take the numerical initial conditions to be the same as the exact steady-state solutions, which does not satisfy the numerical schemes and will be driven by the iterative schemes to the numerical steady states.
In Table \ref{1.1}, we report the numerical accuracy for the density variable, iteration numbers and CPU times when these three different iterative schemes reach the average residue threshold value $10^{-12}$ of the convergence criterion. It is observed that all three schemes achieve basically the same  numerical errors and the fifth order accuracy when the iterations converge. This is as expected since although they are different iterative schemes, they converge to the solution of the same nonlinear algebraic system resulted from the fifth order finite volume WENO local solver. However, these methods exhibit very different computational efficiency by comparing the iteration numbers and CPU times required by them to reach steady state. As shown in Table \ref{1.1}, the FE Jacobi scheme requires a very small CFL number $0.1$ to achieve the convergence. The reason is that a forward Euler time discretization with a high order linear upwind spatial discretization suffers from linear stability issue. When the nonlinearly stable WENO discretization is applied, it alleviates the linear instability problem. As a result of balance, the forward Euler time marching scheme, i.e. the FE Jacobi scheme, converges under a tiny CFL number in this example. That leads to the largest iteration numbers and the most CPU time costs among these three iterative schemes. By using the TVD-RK3 scheme (i.e. the RK Jacobi scheme), both linear and nonlinear stability are maintained. Hence here a larger CFL number $0.6$ can be used to make the iterations converge. Although the RK Jacobi scheme has three stages in one time step, in this example it needs fewer number of iterations and less CPU costs to converge to the steady state, and is more efficient than the FE Jacobi scheme. As shown in Table \ref{1.1}, the proposed fixed-point fast sweeping method (i.e. the FE fast sweeping scheme) is the most efficient one among all three iterative methods. Furthermore, it is important to notice that the FE fast sweeping scheme permits a CFL number $0.6$ which is the same as the TVD-RK3 scheme. This shows that the fixed-point fast sweeping method can improve the linear stability of the forward Euler scheme when it is coupled with a high order spatial discretization. Hence, it suggests that by applying the fixed-point fast sweeping technique to the forward Euler scheme, the forward Euler scheme with a high order WENO spatial discretization becomes practically useful and efficient in solving steady-state hyperbolic conservation laws. In Table \ref{1.1}, it is observed that on refined meshes, the number of iterations for the FE fast sweeping method to converge to the steady state is only about $13\% \sim 15\%$ of that for the forward Euler time marching scheme, and about $25\% \sim 30\%$ of that for the TVD-RK3 scheme. Since the Gauss-Seidel iterations in the FE fast sweeping method (\ref{14}) require that the newest numerical values of $\overline{\textbf{u}}$ are used in the finite volume WENO reconstruction stencils as long as they are available, the WENO reconstructions for both $\textbf{u}^{-}$ and $\textbf{u}^{+}$ in (\ref{4}) need to be performed twice for the same Gaussian quadrature points on the shared boundary of two different target cells. This is the reason that the CPU time cost of the FE fast sweeping method is about $26\% \sim 30\%$ of that for the forward Euler time marching scheme, and about $50\% \sim 60\%$ of that for the TVD-RK3 scheme. Hence on refined meshes in this example, the FE fast sweeping method saves about $70\%$ CPU time cost for the forward Euler time marching scheme, and saves around $50\%$ CPU time cost for the TVD-RK3 scheme.

\begin{table}%[H]
		\centering
\begin{tabular}{|c|c|c|c|c|c|c|}\hline
			\multicolumn{7}{|c|}{FE Jacobi, CFL=0.1 }\\\hline
			$M$&$L_{1}$ error&order&$L_{\infty} $ error&order&iter$\sharp$&CPU time\\\hline
    58    & 1.24E-02 & -  & 3.79E-02     & -  & 5168   & 10.13 \\
    \hline
    232    & 4.95E-04 & 4.64  & 1.63E-03    & 4.54  & 7387   & 49.33 \\
    \hline
    928    & 1.58E-05 & 4.97  & 5.19E-05    & 4.97  & 9306   & 231.73\\
    \hline
    3712   & 5.10E-07 & 4.95   & 1.54E-06  & 5.07   & 13919  & 1350.48 \\
    \hline
    14848   & 1.62E-08 & 4.97  & 5.94E-08   & 4.70  & 24950  & 10003.34 \\
    \hline
      59392   & 5.13E-10 & 4.98  & 2.56E-09   & 4.54  &40437   & 64344.13 \\
    \hline
		\end{tabular}
\begin{tabular}{|c|c|c|c|c|c|c|}\hline
			\multicolumn{7}{|c|}{RK Jacobi,  CFL=0.6 }\\\hline
			$M$&$L_{1}$ error& order&$L_{\infty} $ error&order&iter$\sharp$&CPU time\\\hline
    58    & 1.24E-02 & -  & 3.79E-02     & -  & 2589   & 5.16 \\
    \hline
    232    & 4.95E-04 & 4.64  & 1.63E-02     & 4.54  & 3687   & 25.08 \\
    \hline
    928    & 1.58E-05 & 4.97   & 5.19E-05    & 4.97  & 4638   & 117.20 \\
    \hline
    3712   & 5.10E-07 & 4.95   & 1.54E-06  & 5.07   & 6945  & 696.03 \\
    \hline
    14848   & 1.62E-08 & 4.97  & 5.94E-08   & 4.70  &12465   &5086.49  \\
    \hline
      59392   & 5.13E-10 & 4.98  & 2.55E-09   & 4.54  &20301   &32920.73  \\
    \hline
		\end{tabular}
\begin{tabular}{|c|c|c|c|c|c|c|}\hline
			\multicolumn{7}{|c|}{FE fast sweeping,  CFL=0.6 }\\\hline
			$M$&$L_{1}$ error& order&$L_{\infty} $ error&order&iter$\sharp$&CPU time\\\hline
   58    & 1.24E-02 & -  & 3.79E-02     & -  & 686   & 1.99 \\
    \hline
    232    & 4.95E-04 & 4.64  & 1.63E-03    & 4.54  & 990   & 11.36 \\
    \hline
    928    & 1.58E-05 & 4.97  & 5.19E-05    & 4.97  & 1384   & 59.08\\
    \hline
    3712   & 5.10E-07 & 4.95   & 1.54E-06  & 5.07   & 1870  & 350.07 \\
    \hline
    14848   & 1.62E-08 & 4.97  & 5.94E-08   & 4.70  & 3326  & 2756.83 \\
    \hline
      59392   & 5.14E-10 & 4.98  & 2.56E-09   &  4.54 &  5926 & 21219.23 \\
      \hline
		\end{tabular}
		\caption{Example 1, an Euler system of equations with source terms. Accuracy, iteration numbers and CPU times of three different iterative schemes. CPU time unit: second.}
		\label{1.1}
	\end{table}

\bigskip
\noindent\textbf{Example 2. An Euler system of equations without source terms}

\noindent In this example, we solve for  steady-state solution of the following two dimensional Euler system of  equations. Different from Example 1, the system does not involve any source terms. The equations are
\begin{equation}
\frac{\partial}{\partial{t}}\left(\begin{array}{cccc}
 \rho\\ \rho{u}\\ \rho{v} \\ E\\
\end{array}\right)
+\frac{\partial}{\partial{x}}\left(\begin{array}{cccc}
 \rho{u}\\ \rho{u^{2}+p}\\ \rho{uv} \\ u(E+p)\\
\end{array}\right)
+\frac{\partial}{\partial{y}}\left(\begin{array}{cccc}
 \rho{v}\\ \rho{uv}\\ \rho{v^{2}+p} \\ v(E+p)\\
\end{array}\right)
=\left(\begin{array}{cccc}
 0\\ 0\\ 0\\ 0\\
\end{array}\right).
\end{equation}
We consider the case that the system has the exact steady-state solution $\rho(x,y,\infty)=1 + 0.2\sin(x-y), u(x,y,\infty)=1, v(x,y,\infty)=1, p(x,y,\infty)=1$, and solve the system by the FE Jacobi scheme, the RK Jacobi scheme and  the FE fast sweeping scheme to compare their performance. The computational domain, boundary conditions, triangular meshes used for the simulations, and the reference points to form the alternating sweeping directions in the FE fast sweeping scheme, are all the same as these in Example 1. Also, to start the iterations for the schemes, we take the numerical initial conditions to be the same as the exact steady-state solutions, which does not satisfy the numerical schemes and will be driven by these iterative schemes to their numerical steady states. In Table \ref{2.1}, we report the numerical accuracy for the density variable, iteration numbers and CPU times when these three different iterative schemes reach the average residue threshold value $10^{-12}$ of the convergence criterion. It can be observed that all three schemes achieve basically the same numerical errors and the fifth order accuracy when the iterations converge. Again,
these methods exhibit very different computational efficiency by comparing the iteration numbers and CPU times required by them to reach steady state. As shown in Table \ref{2.1}, the proposed FE fast sweeping scheme is the most efficient one among all three iterative methods. It permits a larger CFL number than the FE Jacobi scheme. On refined meshes, the number of iterations for the FE fast sweeping method to converge to the steady state is only about $13\%$ of that for the FE Jacobi scheme, and about $17\%$ of that for the RK Jacobi scheme. The Gauss-Seidel iterations in the FE fast sweeping method require that the newest numerical values are used in the finite volume WENO reconstruction stencils, so the WENO reconstructions for both $\textbf{u}^{-}$ and $\textbf{u}^{+}$ in (\ref{4}) need to be performed twice for the same Gaussian quadrature points on the shared boundary of two different target cells. Hence, we observe that the CPU time cost of the FE fast sweeping method is about $26\%$ of that for the FE Jacobi scheme, and about $33\%$ of that for the RK Jacobi scheme. On refined meshes in this example, the FE fast sweeping method saves about $74\%$ CPU time cost for the forward Euler time marching scheme, and saves about $67\%$ CPU time cost for the TVD-RK3 scheme.

\begin{table}%[H]
		\centering
\begin{tabular}{|c|c|c|c|c|c|c|}\hline
			\multicolumn{7}{|c|}{FE Jacobi, CFL=0.1 }\\\hline
			$N$&$L_{1}$ error&order&$L_{\infty} $ error&order&iter$\sharp$&CPU time\\\hline
      58    & 2.28E-02 &   & 7.62E-02     &   & 933   & 1.82 \\
    \hline
    232    & 1.08E-03 & 4.40  & 4.88E-03     & 3.97  & 1191   & 8.01 \\
    \hline
    928    & 3.97E-05 & 4.76   & 1.97E-04    & 4.63  & 1577   & 40.14 \\
    \hline
    3712   & 1.52E-06 & 4.71   & 7.85E-06  & 4.65   & 2463  & 252.57 \\
    \hline
    14848   & 5.13E-08 & 4.89  & 2.73E-07   & 4.85  &4532   &1777.16 \\
    \hline
  %    59392   & 1.65E-09 & 4.95  & 8.86E-09   & 4.94  &8689   &13406.40  \\
  %  \hline
		\end{tabular}
\begin{tabular}{|c|c|c|c|c|c|c|}\hline
			\multicolumn{7}{|c|}{RK Jacobi, CFL=0.6 }\\\hline
			$N$&$L_{1}$ error& order&$L_{\infty} $ error&order&iter$\sharp$&CPU time\\\hline
    58    & 2.28E-02 &   & 7.62E-02     &   & 705   & 1.58 \\
    \hline
    232    & 1.08E-03 & 4.40  & 4.88E-03     & 3.97  & 876   & 6.46 \\
    \hline
    928    & 3.97E-05 & 4.76   & 1.97E-04    & 4.63  & 1182   & 31.25 \\
    \hline
    3712   & 1.52E-06 & 4.71   & 7.85E-06  & 4.65   & 1884  & 195.02 \\
    \hline
    14848   & 5.13E-08 & 4.89  & 2.73E-07   & 4.85  &3360   &1372.57 \\
    \hline
  %    59392   & 1.65E-09 & 4.95  & 8.86E-09   & 4.94  &6363   &10279.58  \\
  %  \hline
		\end{tabular}
\begin{tabular}{|c|c|c|c|c|c|c|}\hline
			\multicolumn{7}{|c|}{FE fast sweeping, CFL=0.6 }\\\hline
			$N$&$L_{1}$ error& order&$L_{\infty} $ error&order&iter$\sharp$&CPU time\\\hline
   58    & 2.28E-02 &   & 7.62E-02     &   & 136   & 0.42 \\
    \hline
    232    & 1.08E-03 & 4.40  & 4.88E-03    & 3.97  & 156   & 1.84 \\
    \hline
    928    & 3.97E-05 & 4.76  & 1.97E-04    & 4.63  & 210   & 10.21\\
    \hline
    3712   & 1.52E-06 & 4.71   & 7.85E-06  & 4.65   & 326  & 64.69 \\
    \hline
    14848   & 5.13E-08 & 4.89  & 2.73E-07   & 4.85  & 580  & 460.33 \\
    \hline
   %   59392   & 1.65E-09 & 4.95  & 8.86E-09   &  4.94 &  1742 & 6481.27 \\
   %   \hline
		\end{tabular}
		\caption{Example 2, an Euler system of equations without source terms. Accuracy, iteration numbers and CPU times of three different iterative schemes. CPU time unit: second.}
		\label{2.1}
	\end{table}
~\\
In the following examples, We study and compare these schemes by solving problems with discontinuous solutions.

\bigskip
\noindent\textbf{Example 3. Regular shock reflection}

\noindent In this example, we test the absolute convergence of the proposed fixed-point fast sweeping method on unstructured meshes, by solving the regular shock reflection problem. The problem is a 2D steady-state Euler system of equations with a reflection condition along the bottom boundary, as described in e.g. \cite{WZ, SSCW, SCW}. It is a typical and difficult benchmark problem of using high order schemes to simulate steady flow. As that found in \cite{SSCW, SCW}, even with advanced techniques to improve the steady-state convergence, it is still difficult for the residue of high order WENO schemes to converge to the level of round off errors. The iteration residue of the fast sweeping method in \cite{WZ}, which is based on the WENO scheme with equal-sized sub-stencils, hangs at the level above $10^{-3.5}$.
However, the proposed unstructured fast sweeping method with unequal-sized sub-stencils can drive the iteration residue to the level of round off errors, as shown in the following numerical results.

\begin{figure}%[H]
\centering
\includegraphics[width=4.4in]{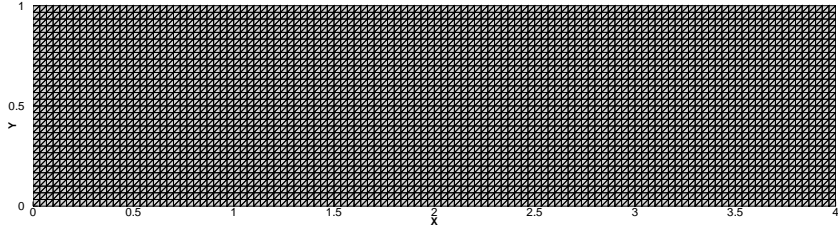}
\caption{\label{7.1} The computational mesh for Example 3, regular shock reflection.}
\end{figure}

\begin{figure}%[H]
\centering
\subfigure[FE Jacobi scheme]{
\begin{minipage}[t]{0.6\linewidth}
\centering
\includegraphics[width=4.0in]{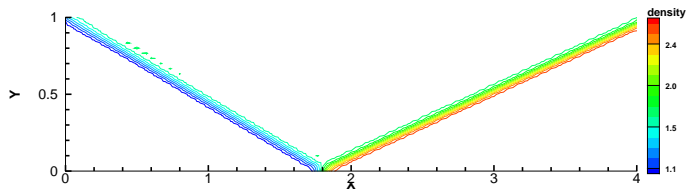}
\end{minipage}}

\centering
\subfigure[RK Jacobi scheme]{
\begin{minipage}[t]{0.6\linewidth}
\centering
\includegraphics[width=4.0in]{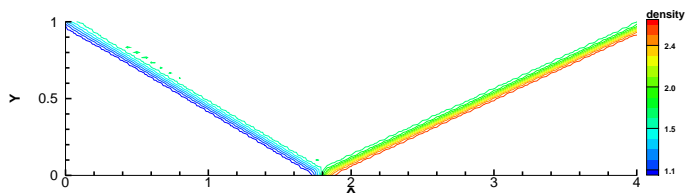}
\end{minipage}}

\centering
\subfigure[FE fast sweeping scheme]{
\begin{minipage}[t]{0.6\linewidth}
\centering
\includegraphics[width=4.0in]{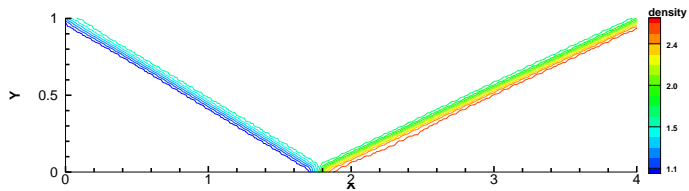}
\end{minipage}}
\caption{\label{7.3}Example 3, regular shock reflection. 30 equally spaced  density contours from 1.1 to 2.6 of the converged steady states of numerical solutions by three different iterative schemes.}
\end{figure}

\begin{figure}%[H]
\centering
\subfigure[FE Jacobi scheme]{
\begin{minipage}[t]{0.3\linewidth}
\centering
\includegraphics[width=2.0in]{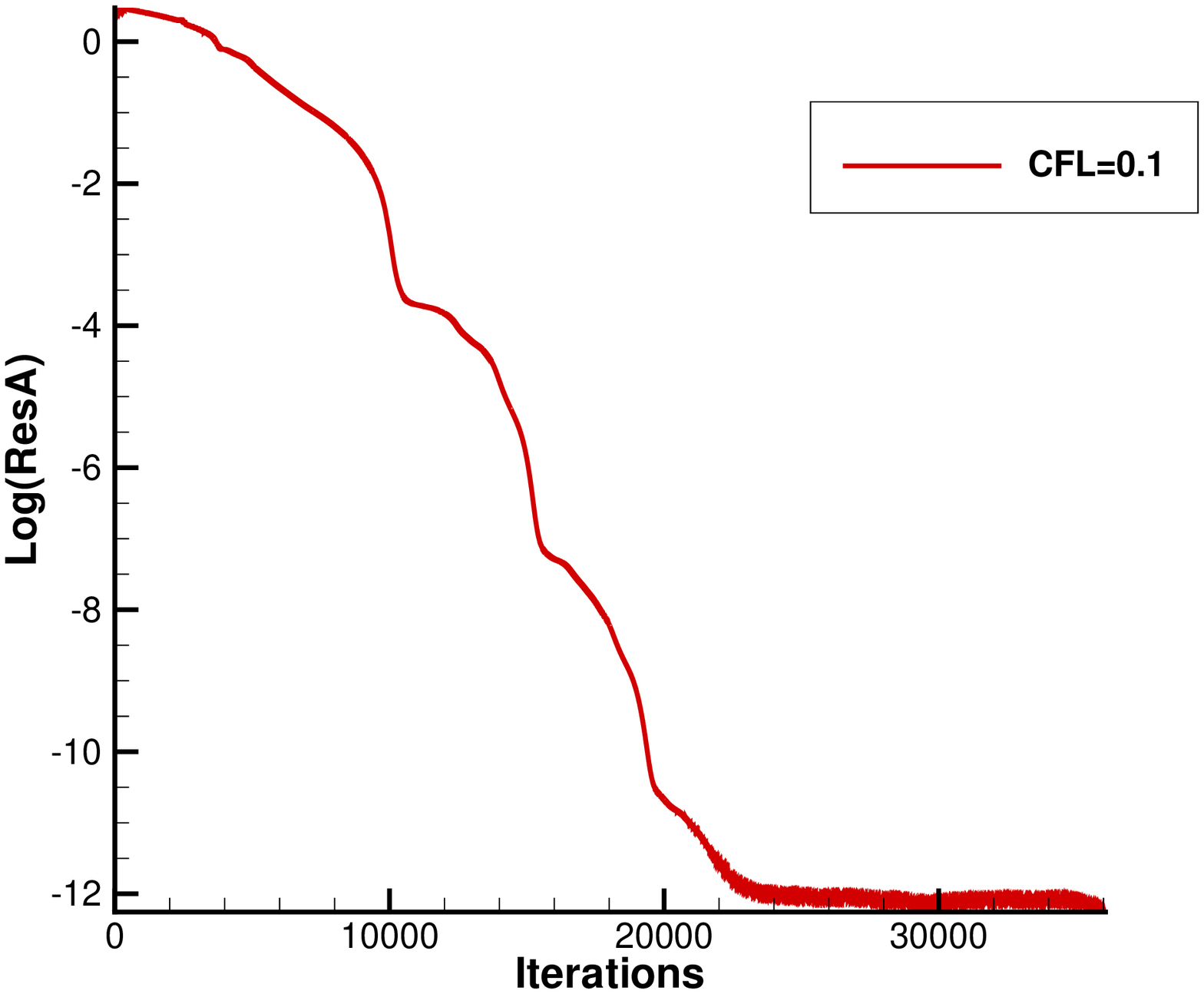}
%\caption{fig1}
\end{minipage}%
}%
\subfigure[RK Jacobi scheme]{
\begin{minipage}[t]{0.3\linewidth}
\centering
\includegraphics[width=2.0in]{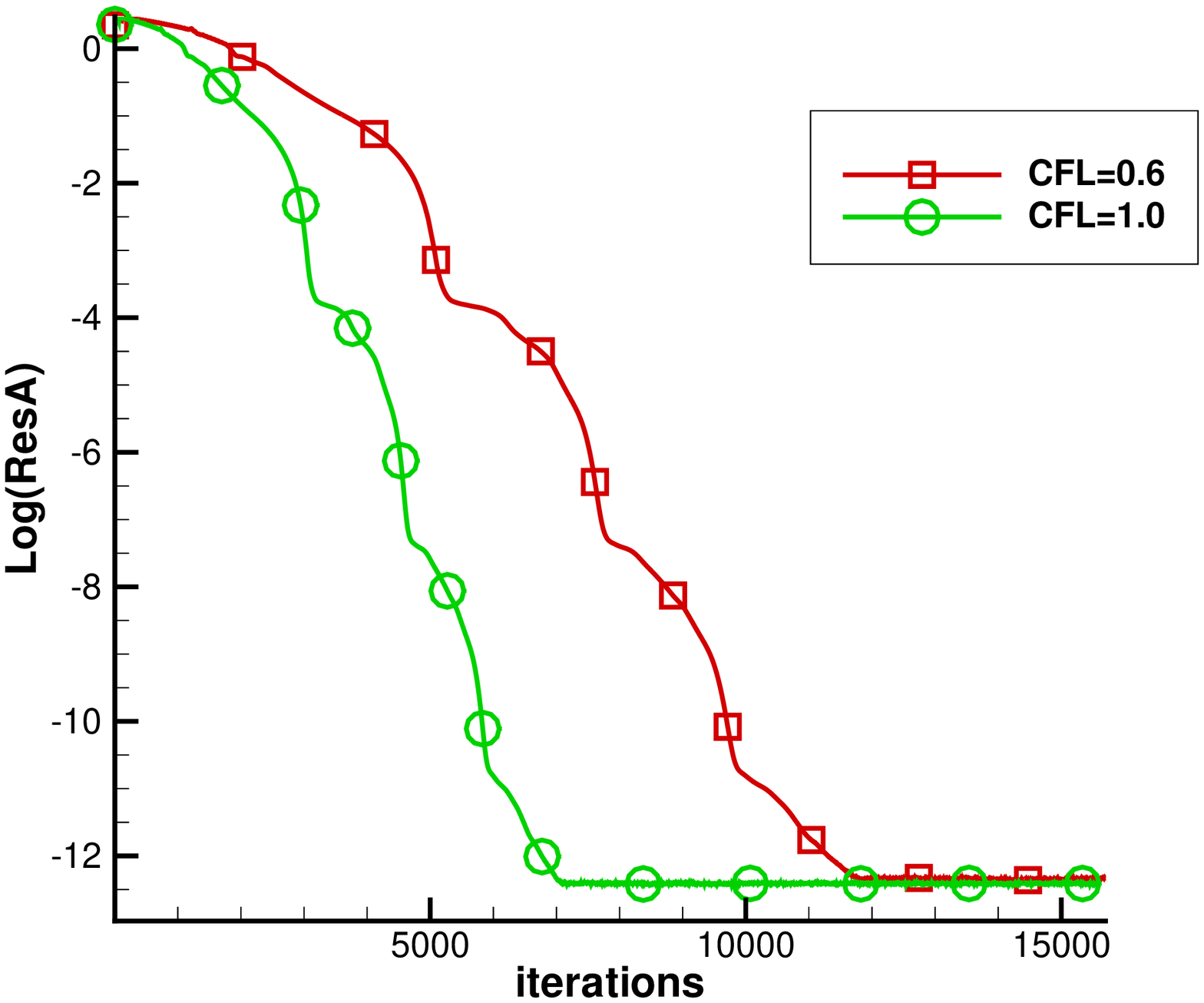}
%\caption{fig1}
\end{minipage}%
}%
\subfigure[FE fast sweeping scheme]{
\begin{minipage}[t]{0.3\linewidth}
\centering
\includegraphics[width=2.0in]{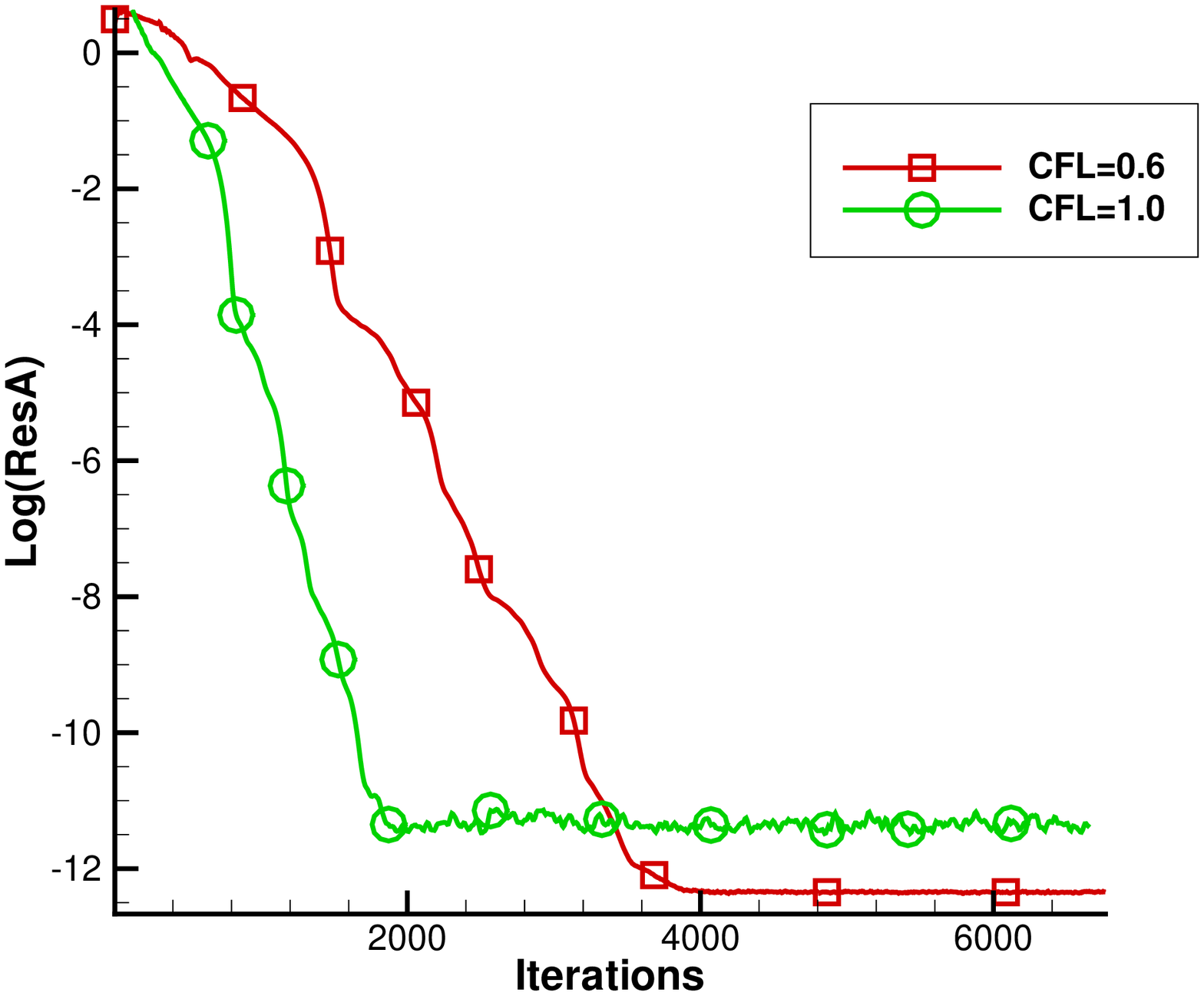}
%\caption{fig2}
\end{minipage}%
}%
\centering
\caption{\label{7.4}Example 3, regular shock reflection. The convergence history of the residue
 as a function of number of iterations for three schemes with different CFL numbers.}
\end{figure}

The computational domain is a rectangle and it has the length 4 and the height 1. The boundary conditions include a reflection condition along the bottom boundary, supersonic outflow along the right boundary, and the Dirichlet conditions on the other two boundaries:
\begin{equation*}
(\rho,u,v,p)^{T}=
\begin{cases}
(1.0,2.9,0,5/7)^{T}\mid_{(0,y)},\\
(1.69997,2.61934,-0.50632,1.52819)^{T}\mid_{(x,1)}.\\
\end{cases}
\end{equation*}
The initial values used to start the iterations in the entire domain are taken to be the same as those at the left boundary.
The computational mesh is shown in Figure \ref{7.1}. Four corners of the domain are taken as the reference points to form the alternating sweeping directions in the FE fast sweeping scheme. In Table \ref{7.2}, number of iterations required
to reach the convergence criterion threshold value $10^{-11}$, and total CPU time when the schemes converge under various CFL numbers are reported for the FE Jacobi scheme, the RK Jacobi scheme and the FE fast sweeping scheme. As the previous examples, the FE Jacobi scheme requires a very small CFL number $0.1$ to achieve the convergence, which results in many
iterations and large CPU time cost. If the RK Jacobi (TVD-RK3) scheme is used, the CFL number can be enlarged to $1.0$ and its iterations
converge much more efficiently than the FE Jacobi scheme. Among these three methods, the FE fast sweeping scheme is still the most efficient one. It can converge under the similar CFL numbers to the RK Jacobi scheme.
With the largest CFL number permitted for each method to reach steady state solution, the proposed FE fast sweeping method on triangular meshes saves about $85\%$ CPU time cost for the forward Euler method, and saves about $50\%$ CPU time cost for the TVD-RK3 method. In Figure \ref{7.3}, density contours of the converged steady state solutions of these three schemes are shown, and the similar numerical results are obtained as expected.
In Figure \ref{7.4}, residue history in terms of iterations for these three schemes with different CFL numbers is shown. It is observed that the residue of iterations settles down to very small values at the level of round off errors for all cases. The absolute convergence for the proposed fast sweeping method on triangular meshes is verified for this difficult benchmark problem.

\begin{table}%[H]
\centering{
\begin{tabular}{|c|c|c|}\hline
\multicolumn{3}{|c|}{FE Jacobi scheme }\\\hline
             CFL number & iteration number  & CPU time \\\hline
	0.1& 20926 	& 6681.16 \\\hline
      0.2   & Not convergent       & - \\
    \hline
\end{tabular}}

\centering{
\begin{tabular}{|c|c|c|}\hline
			\multicolumn{3}{|c|}{RK Jacobi scheme }\\\hline
             CFL number & iteration number  & CPU time \\\hline
	0.6&  10308 	& 3257.98 \\\hline
1.0&  6183 	& 1959.48 \\\hline
           1.1   & Not convergent       & - \\
    \hline
\end{tabular}}

\centering{
\begin{tabular}{|c|c|c|}\hline
			\multicolumn{3}{|c|}{FE fast sweeping scheme }\\\hline
             CFL number & iteration number  & CPU time \\\hline
     	0.6& 3328	& 1851.67 \\\hline
     	1.0& 1808 	& 1020.90 \\\hline
       1.1   & Not convergent       & - \\
    \hline
\end{tabular}}
\caption{\label{7.2}Example 3, regular shock reflection. Number of iterations and total CPU time when convergence is obtained. Convergence criterion threshold value  is $10^{-11}$. CPU time unit: second.}
%\label{tab:Margin_settings}
\end{table}

%~\\

\bigskip
\noindent\textbf{Example 4. The supersonic flow past a circular cylinder}

\begin{figure}%[H]
\centering
\subfigure[]{
\centering
\includegraphics[width=2.4in]{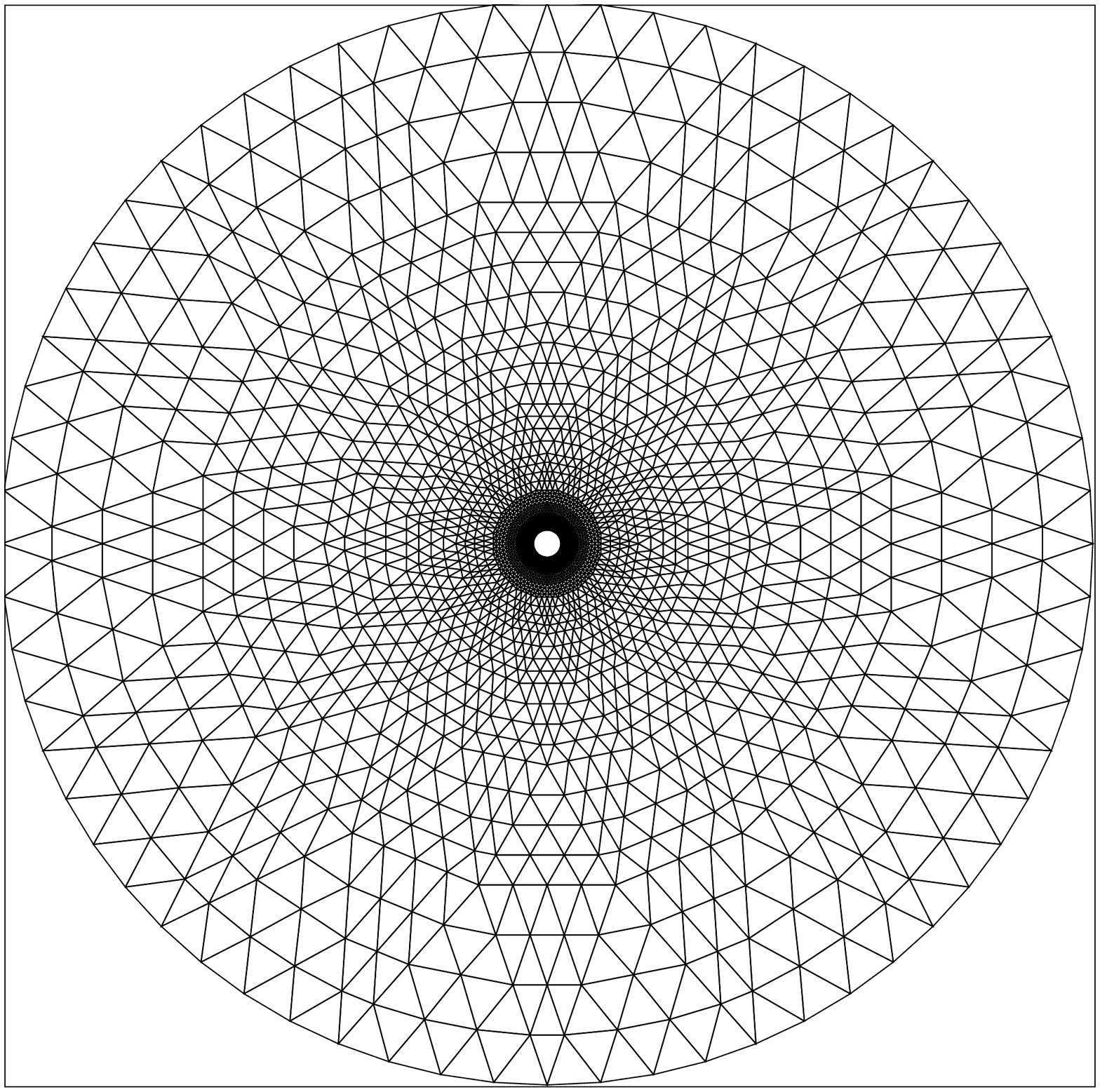}
}
\subfigure[]{
\centering
\includegraphics[width=2.4in]{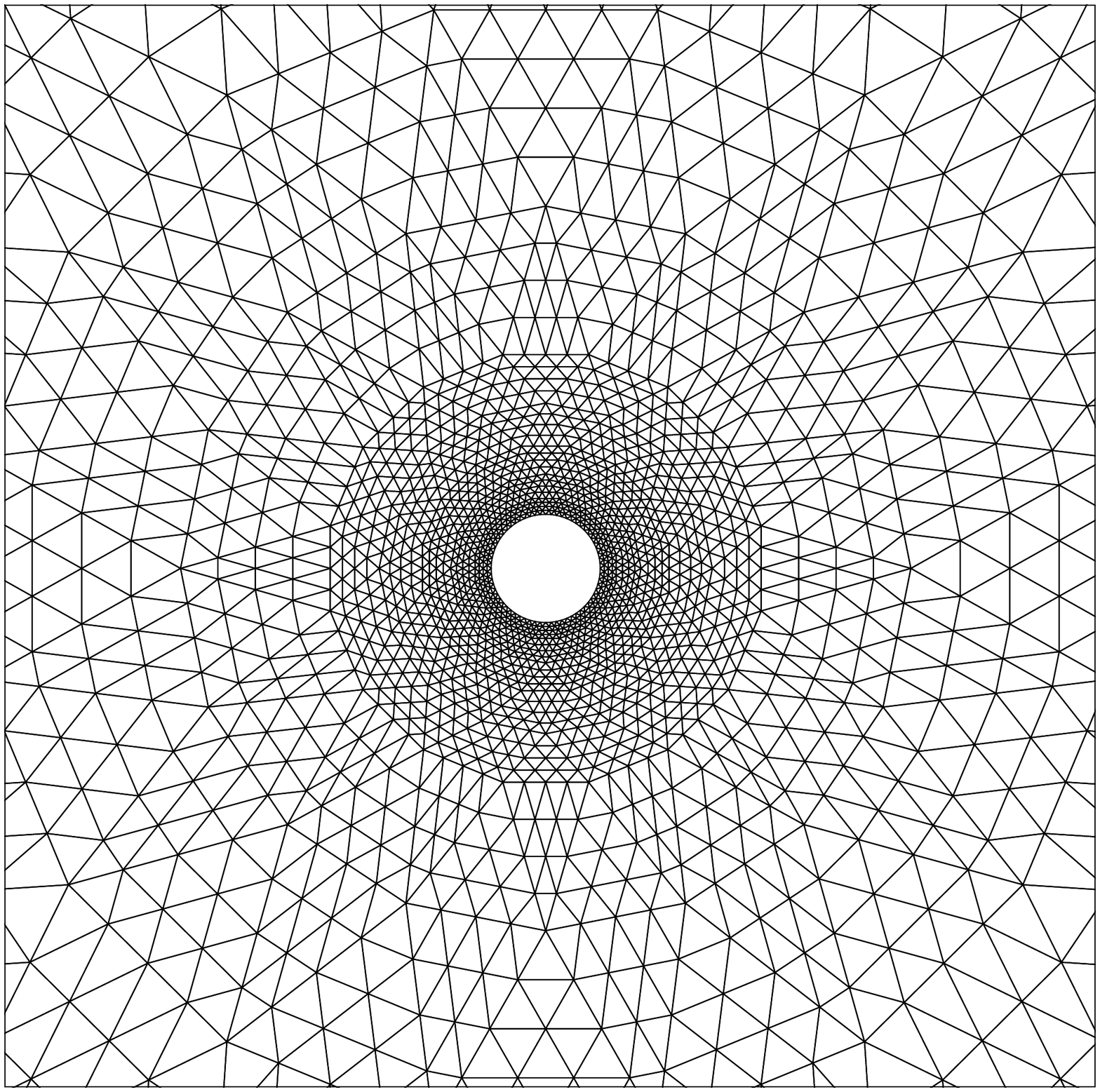}
}
\caption{\label{8.1}The computational mesh for Example 4, the supersonic flow past a circular cylinder. Left: the whole region; right: zoomed near the circular cylinder.}
\end{figure}

\begin{figure}%[H]
\centering
\subfigure[FE Jacobi scheme]{
\begin{minipage}[t]{0.3\linewidth}
\centering
\includegraphics[width=2.0in]{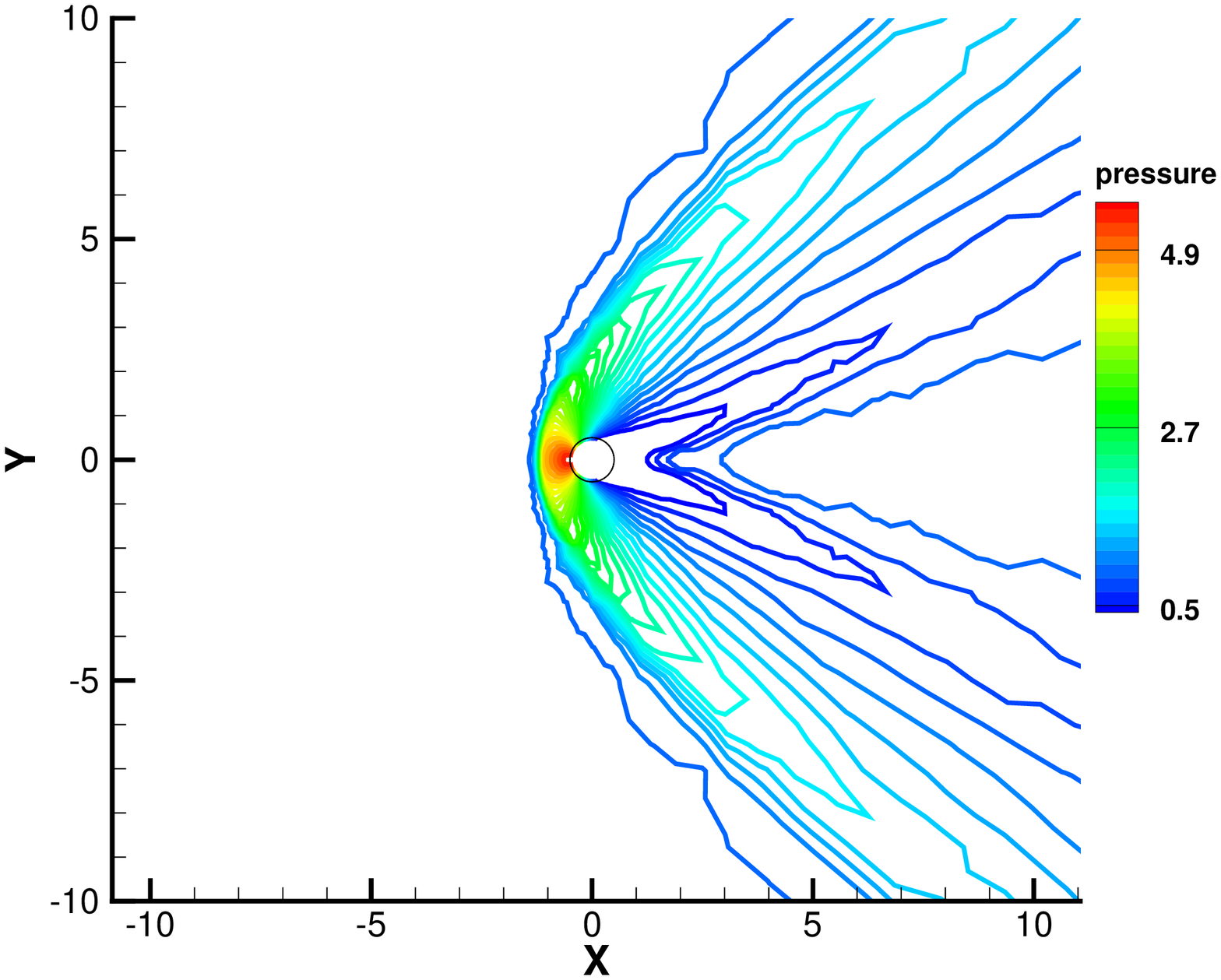}
%\caption{fig1}
\end{minipage}%
}%
\subfigure[RK Jacobi scheme]{
\begin{minipage}[t]{0.3\linewidth}
\centering
\includegraphics[width=2.0in]{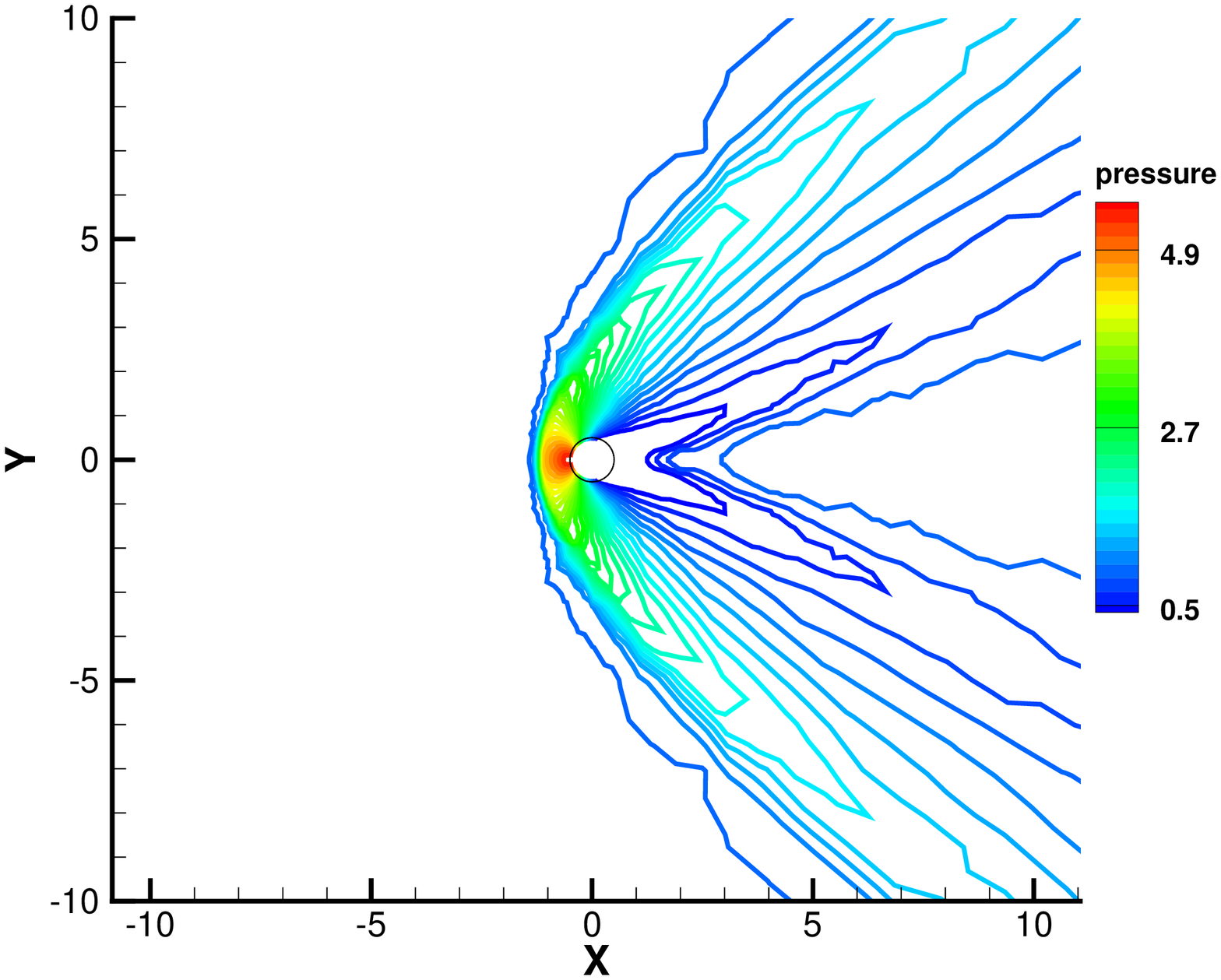}
%\caption{fig1}
\end{minipage}%
}%
\subfigure[FE fast sweeping scheme]{
\begin{minipage}[t]{0.3\linewidth}
\centering
\includegraphics[width=2.0in]{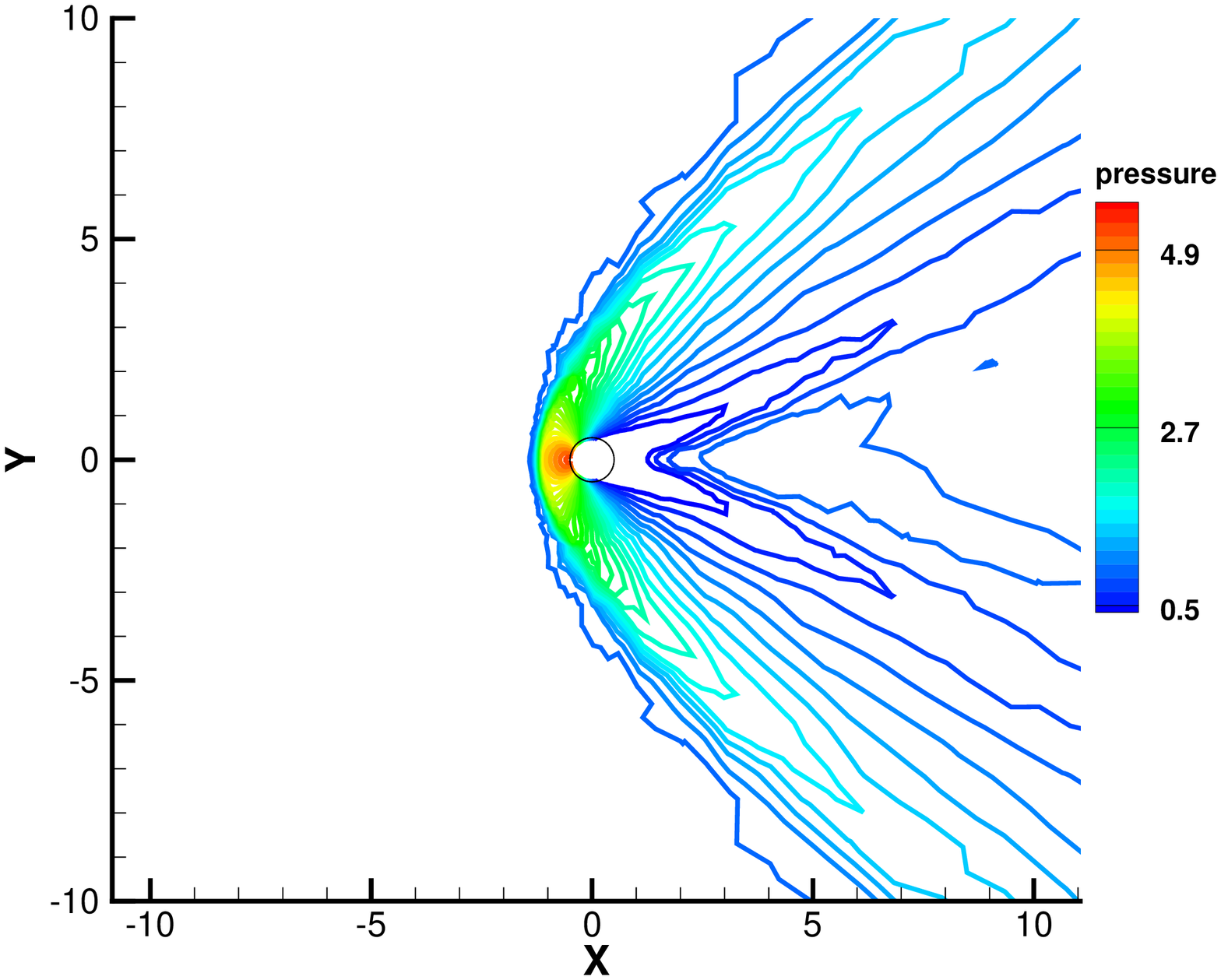}
%\caption{fig2}
\end{minipage}%
}%
\centering
\caption{\label{8.3}Example 4, the supersonic flow past a circular cylinder. 30 equally spaced pressure contours from 0.5 to 5.4 of the converged steady states of numerical solutions by three different iterative schemes.}
\end{figure}

\begin{figure}%[H]
\centering
\subfigure[FE Jacobi scheme]{
\begin{minipage}[t]{0.3\linewidth}
\centering
\includegraphics[width=2.0in]{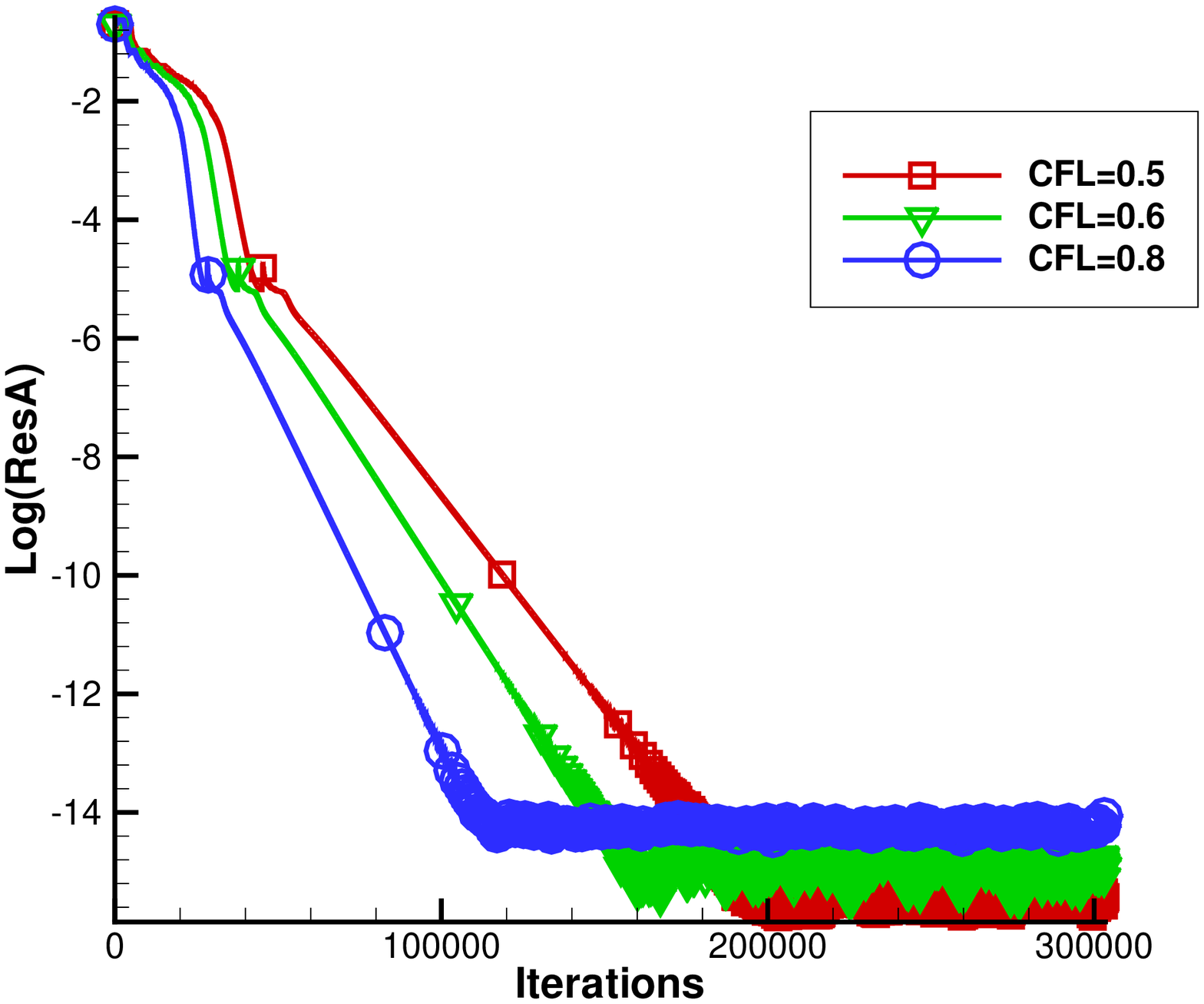}
%\caption{fig1}
\end{minipage}%
}%
\subfigure[RK Jacobi scheme]{
\begin{minipage}[t]{0.3\linewidth}
\centering
\includegraphics[width=2.0in]{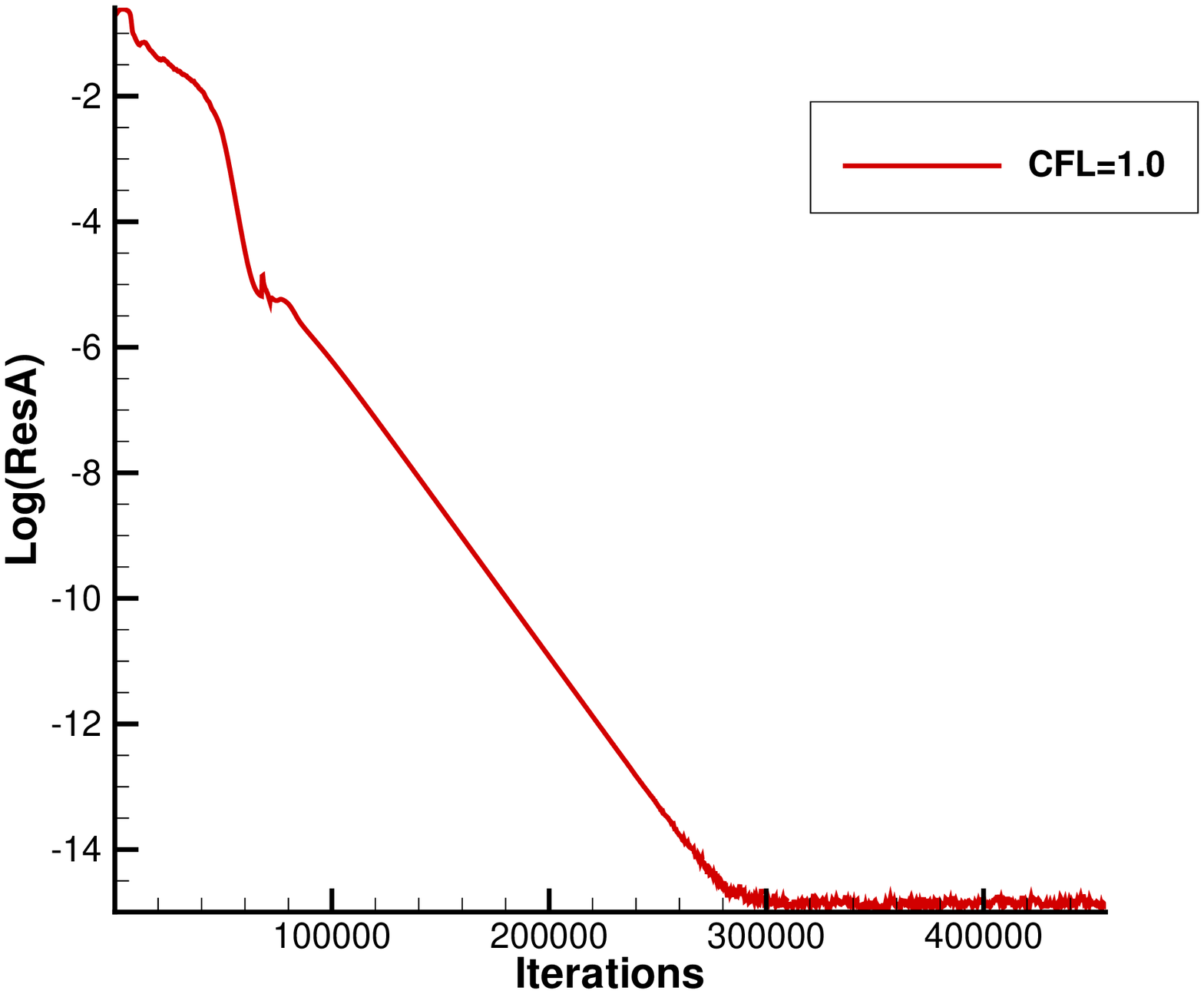}
%\caption{fig1}
\end{minipage}%
}%
\subfigure[FE fast sweeping scheme]{
\begin{minipage}[t]{0.3\linewidth}
\centering
\includegraphics[width=2.0in]{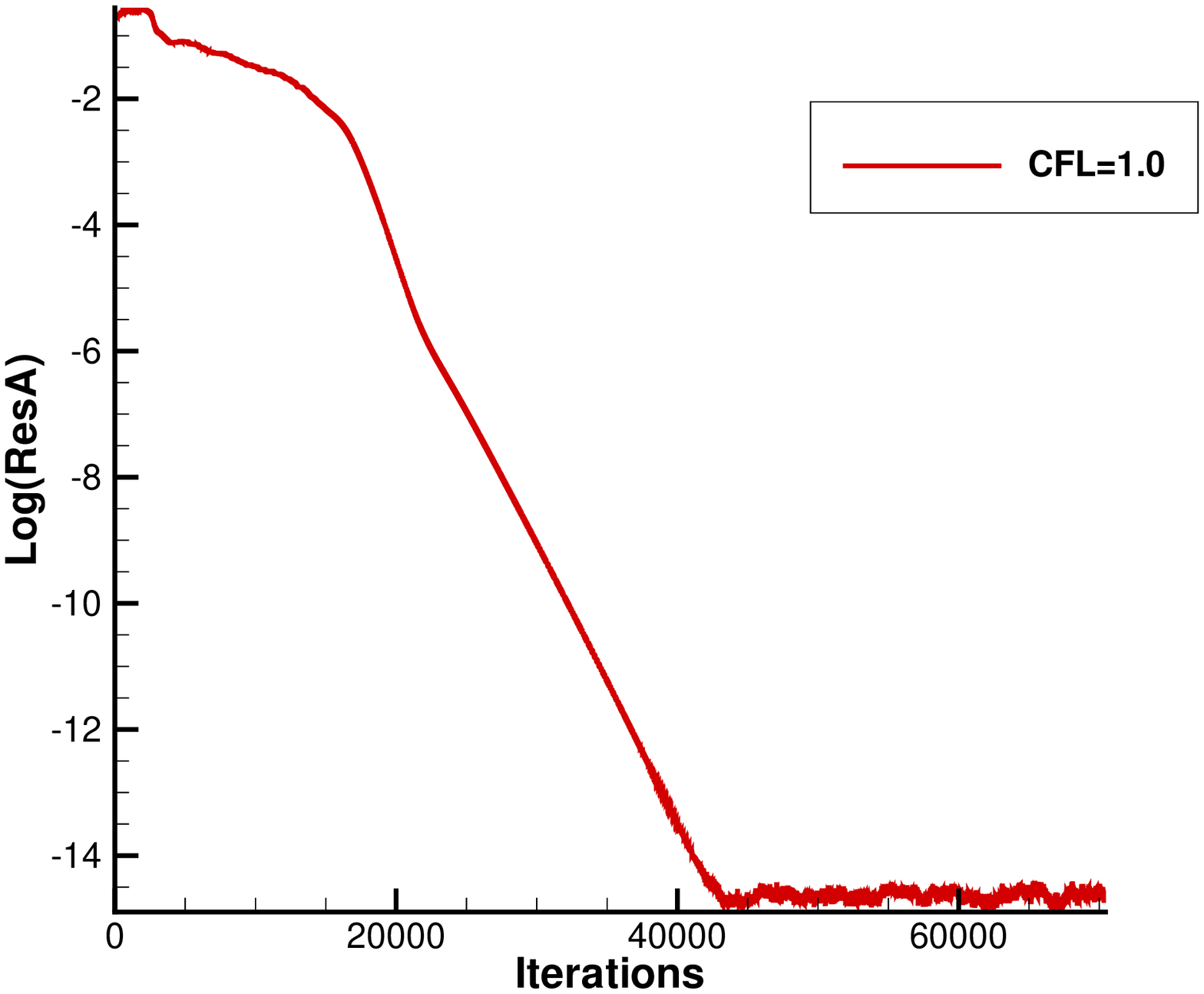}
%\caption{fig2}
\end{minipage}%
}%
\centering
\caption{\label{8.4}Example 4, the supersonic flow past a circular cylinder. The convergence history of the residue
 as a function of number of iterations for three schemes with different CFL numbers. }
\end{figure}

\noindent In this example, we consider an inviscid, compressible flow which initially moves toward a circular cylinder from the left at a Mach number of $Ma=2$ \cite{LBL1,LBL}. The cylinder with the radius 0.5 is located at the origin, and the computational domain is $\{(x,y): 0.5\leq \sqrt{x^{2}+y^{2}}\leq 20\}$. The unstructured triangular mesh used in this example is shown in Figure \ref{8.1}, in which the number of grid points on the boundaries is 64. As discussed in Section 2, for such circular domain here, we take four points which are evenly distributed on the outer boundary, $(-20,0),(0,-20),(20,0)$ and $(0,20)$, as the reference points to form the alternating sweeping directions in the FE fast sweeping scheme. In Table \ref{8.2}, number of iterations required
to reach the convergence criterion threshold value $10^{-11}$,  and total CPU time when the schemes converge under various CFL numbers are reported for the FE Jacobi scheme, the RK Jacobi scheme, and the FE fast sweeping scheme. Different from the previous examples, in this example the FE Jacobi scheme converges under a reasonable CFL number. Although the RK Jacobi scheme still converges under a larger CFL number than the FE Jacobi scheme, it is less efficient here due to its multi-stage structure. Consistent with the previous examples, the FE fast sweeping scheme is still the most efficient one among three methods, and it permits the similar CFL numbers to the RK Jacobi scheme. With the largest CFL number allowed for the methods to reach steady state solution, the FE fast sweeping method on unstructured triangular meshes saves about $70\%$ CPU time cost of that by the RK Jacobi scheme (the TVD-RK3 scheme). In Figure \ref{8.3}, the pressure contours of the converged steady state solutions of these three schemes are presented. We observe similar numerical steady states of these schemes. The residue history of these three schemes with different CFL numbers is reported in Figure \ref{8.4}, which shows that the residue of iterations settles down to very small values at the level of round off errors and verifies the absolute convergence of the proposed high order fast sweeping method on triangular meshes.

\begin{table}%[H]
		\centering
\begin{tabular}{|c|c|c|}\hline
			\multicolumn{3}{|c|}{FE Jacobi scheme }\\\hline
             CFL number & iteration number & CPU time \\\hline
	0.5& 132862    	& 22369.13 \\\hline
	0.6& 110709   	& 20776.28 \\\hline
    0.8& 83015     & 15220.13 \\\hline
    0.9& Not convergent      & - \\
    \hline
		\end{tabular}

\begin{tabular}{|c|c|c|}\hline
			\multicolumn{3}{|c|}{RK Jacobi scheme }\\\hline
             CFL number & iteration number & CPU time \\\hline
	1.0&  201529  	& 36827.34 \\\hline
           1.1   & Not convergent        & - \\
    \hline
		\end{tabular}

\begin{tabular}{|c|c|c|}\hline
			\multicolumn{3}{|c|}{FE fast sweeping scheme }\\\hline
             CFL number & iteration number  & CPU time \\\hline
     	1.0& 34482 & 10957.96 \\\hline
       1.1   & Not convergent        & - \\
    \hline
		\end{tabular}
		\caption{\label{8.2}Example 4, the supersonic flow past a circular cylinder. Number of iterations and total CPU time when convergence is obtained. Convergence criterion threshold value is $10^{-11}$. CPU time unit: second.}
%		\label{tab:Margin_settings}
	\end{table}

~\\
\textbf{Example 5. Supersonic and subsonic flows past an NACA001035 airfoil}

\begin{figure}%[H]
\centering
\subfigure[]{
\centering
\includegraphics[width=2.4in]{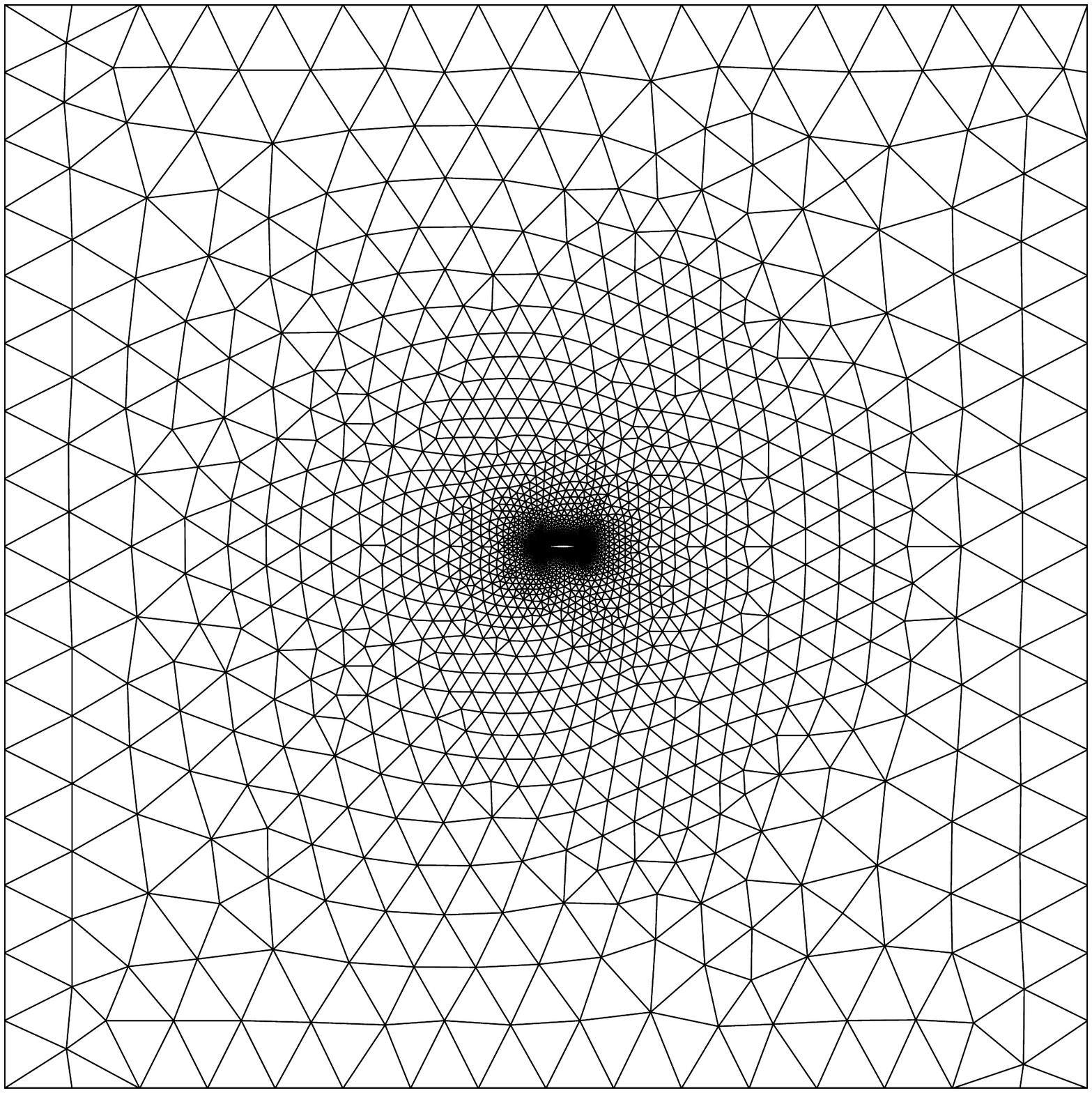}
}
\subfigure[]{
\centering
\includegraphics[width=2.4in]{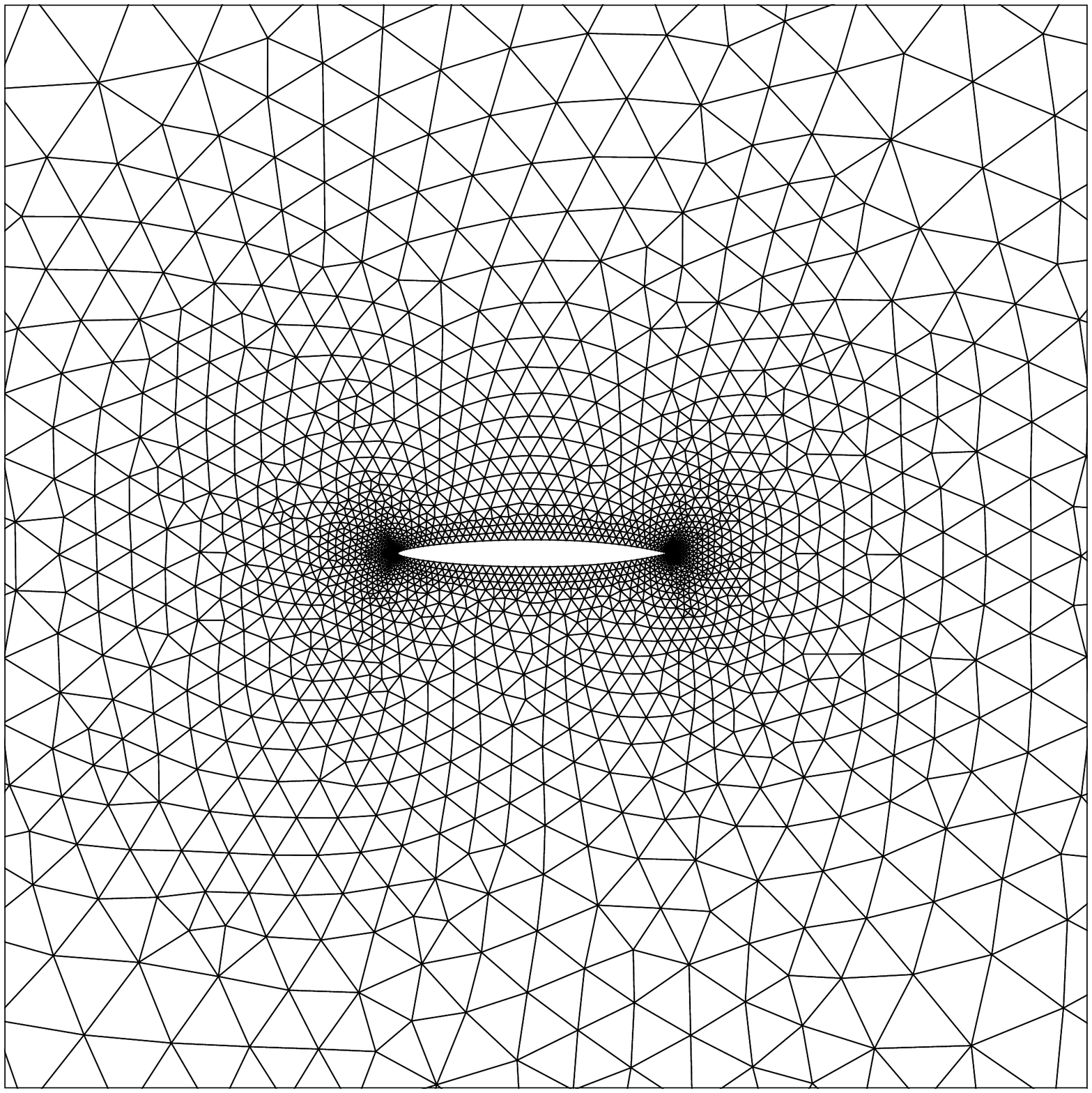}
}
\caption{\label{9.1} The computational mesh for Example 5, supersonic and subsonic flows past an NACA001035 airfoil. Left: the whole domain; right: zoomed region near the airfoil.}
\end{figure}

\begin{figure}%[H]
\centering
\subfigure[FE Jacobi, $Ma=2$]{
\begin{minipage}[t]{0.3\linewidth}
\centering
\includegraphics[width=2.0in]{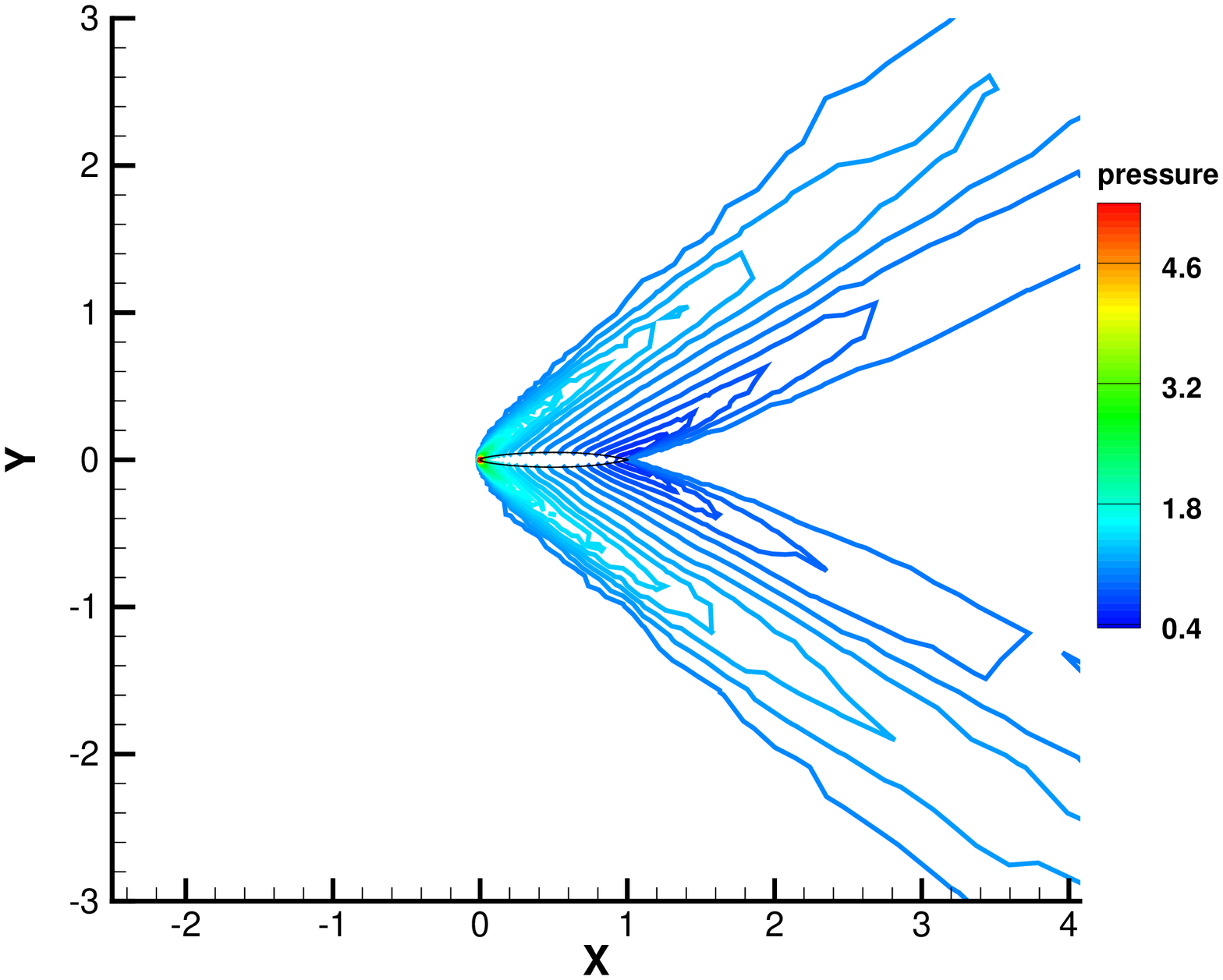}
%\caption{fig1}
\end{minipage}%
}%
\subfigure[RK Jacobi, $Ma=2$]{
\begin{minipage}[t]{0.3\linewidth}
\centering
\includegraphics[width=2.0in]{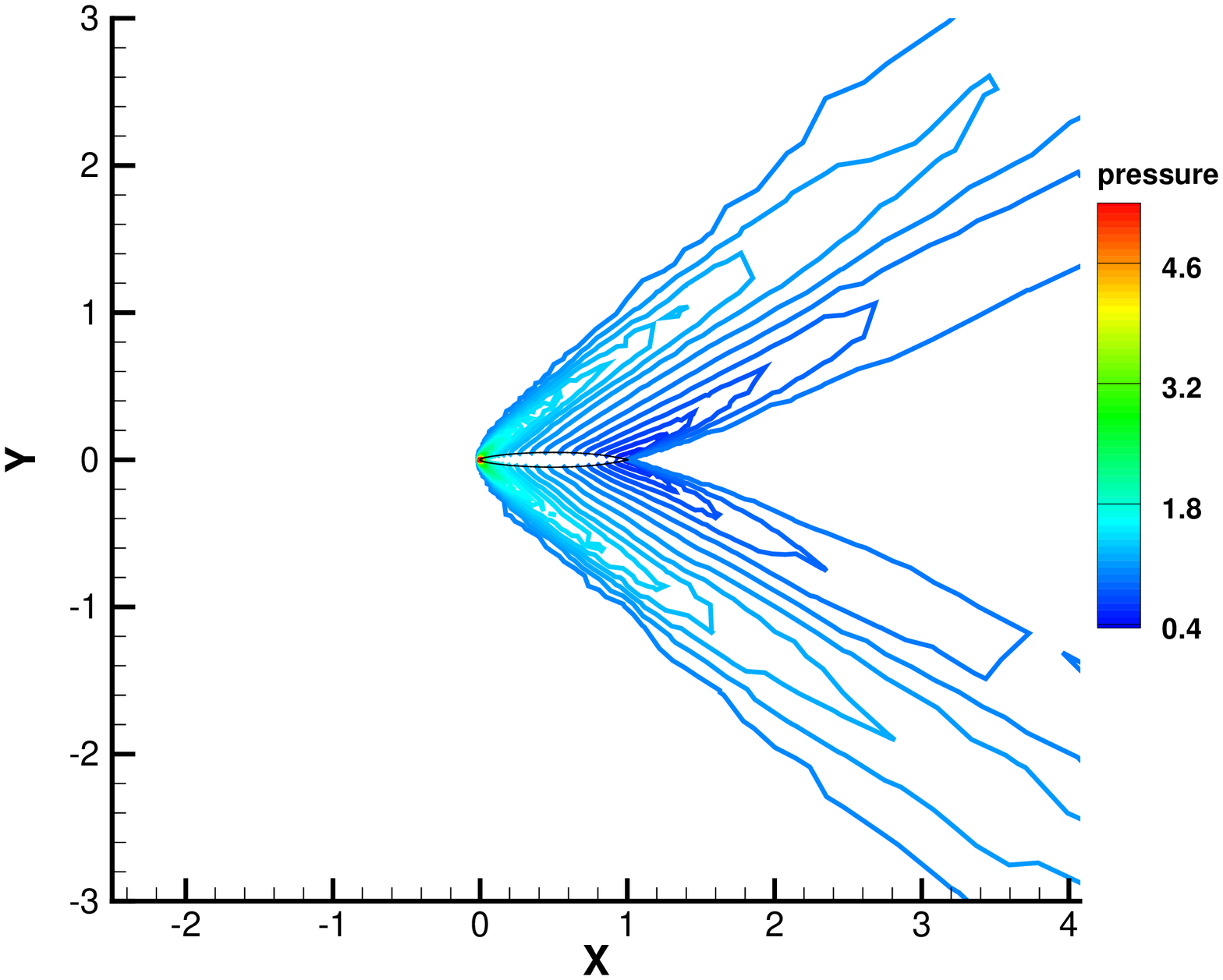}
%\caption{fig1}
\end{minipage}%
}%
\subfigure[FE fast sweeping, $Ma=2$]{
\begin{minipage}[t]{0.3\linewidth}
\centering
\includegraphics[width=2.0in]{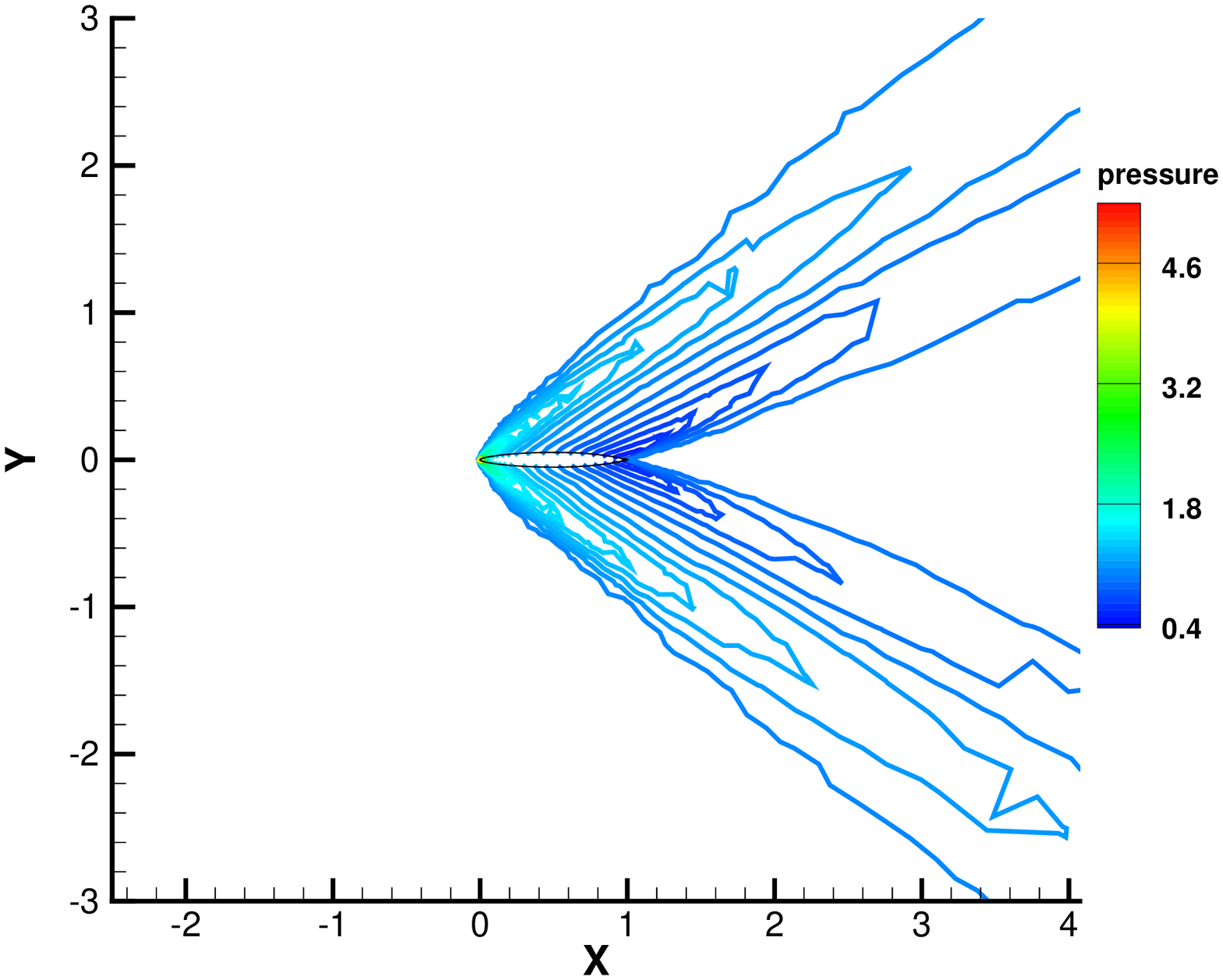}
%\caption{fig2}
\end{minipage}%
}%
\newline
\subfigure[FE Jacobi, $Ma=0.8$]{
\begin{minipage}[t]{0.3\linewidth}
\centering
\includegraphics[width=2.0in]{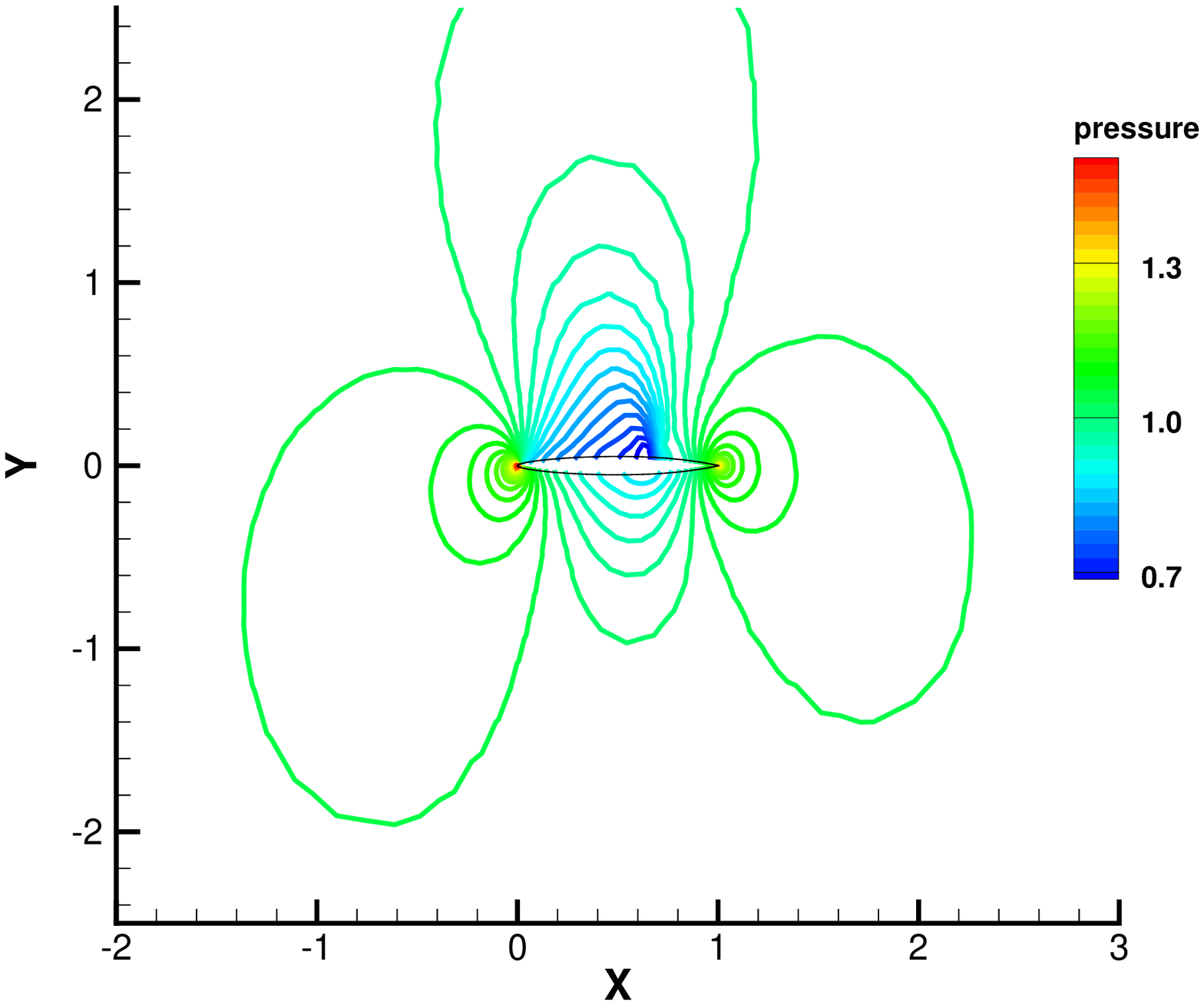}
%\caption{fig1}
\end{minipage}%
}%
\subfigure[RK Jacobi, $Ma=0.8$]{
\begin{minipage}[t]{0.3\linewidth}
\centering
\includegraphics[width=2.0in]{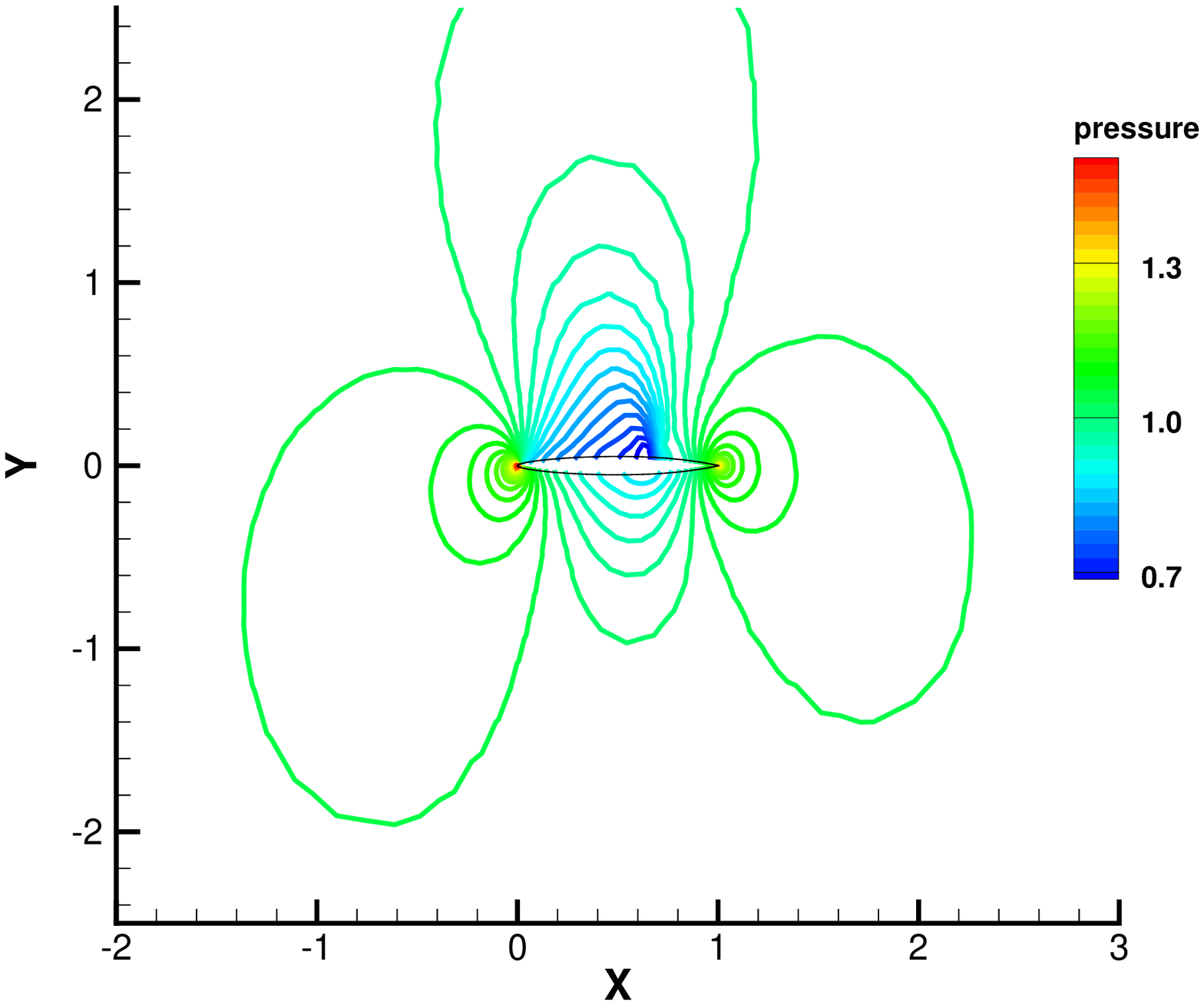}
%\caption{fig1}
\end{minipage}%
}%
\subfigure[FE fast sweeping, $Ma=0.8$]{
\begin{minipage}[t]{0.3\linewidth}
\centering
\includegraphics[width=2.0in]{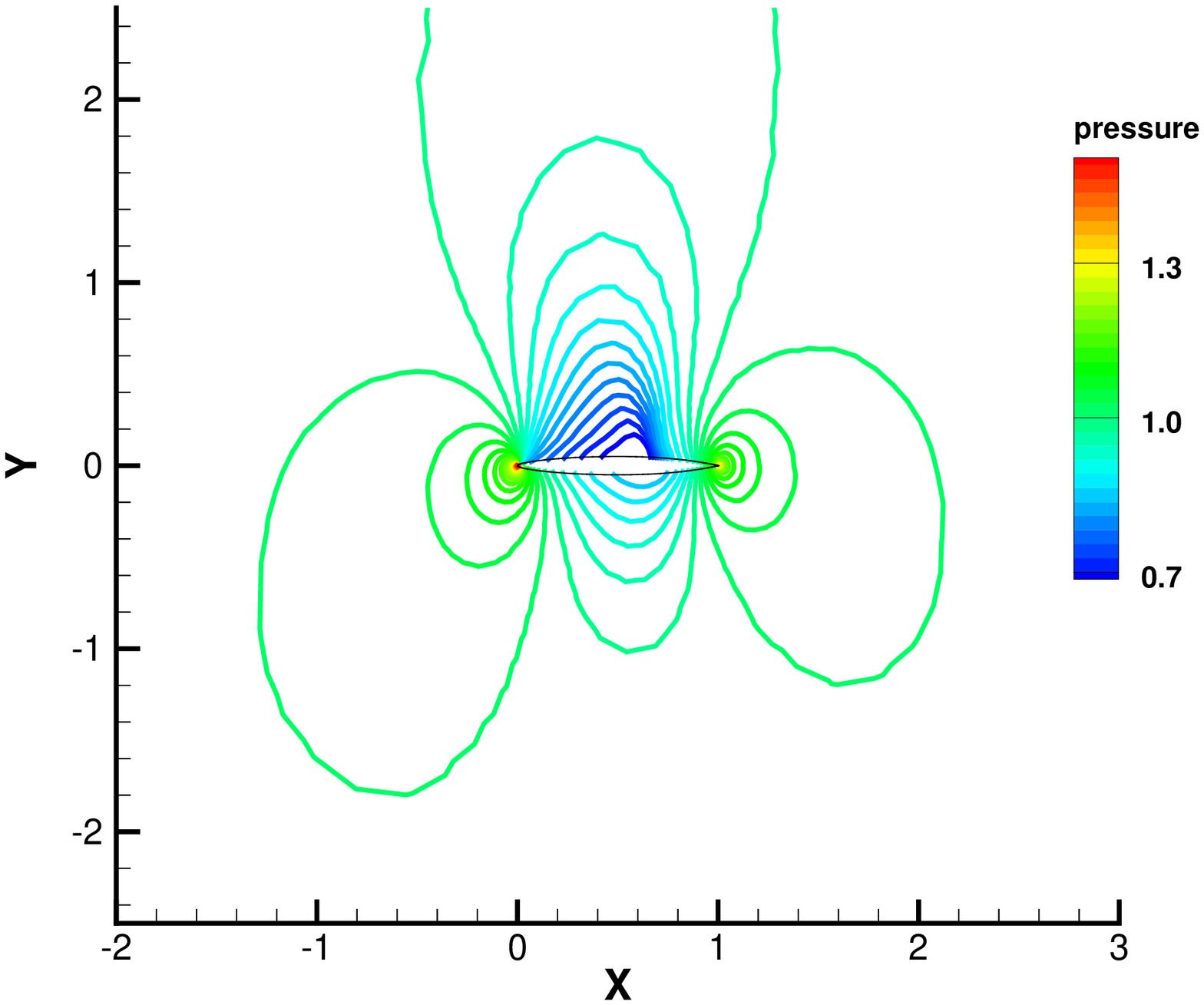}
%\caption{fig2}
\end{minipage}%
}%
\centering
\caption{\label{9.3}Example 5, supersonic and subsonic flows past an NACA001035 airfoil. The converged steady states of numerical solutions by three different iterative schemes. (a) (b) (c): 60 equally spaced pressure contours from 0.36 to 5.27 for the case of $Ma=2$, angle of attack $\alpha=1^{\circ}$; (d) (e) (f): 30 equally spaced pressure contours from 0.7 to 1.45 for the case of $Ma=0.8$, angle of attack $\alpha=1.25^{\circ}$.}
\end{figure}

\begin{figure}%[H]
\centering
\subfigure[FE Jacobi, $Ma=2$]{
\begin{minipage}[t]{0.3\linewidth}
\centering
\includegraphics[width=2.0in]{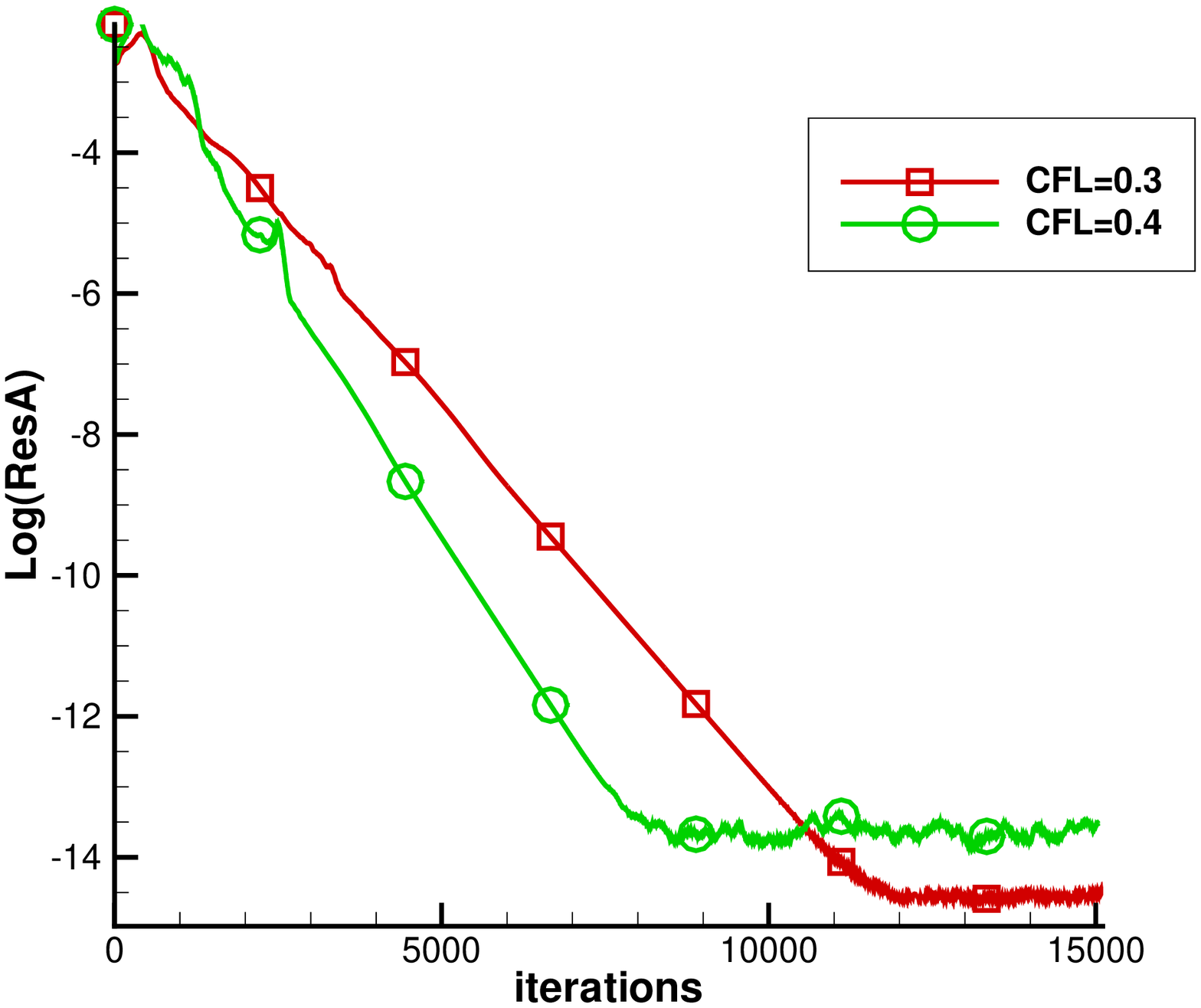}
%\caption{fig1}
\end{minipage}%
}%
\subfigure[RK Jacobi, $Ma=2$]{
\begin{minipage}[t]{0.3\linewidth}
\centering
\includegraphics[width=2.0in]{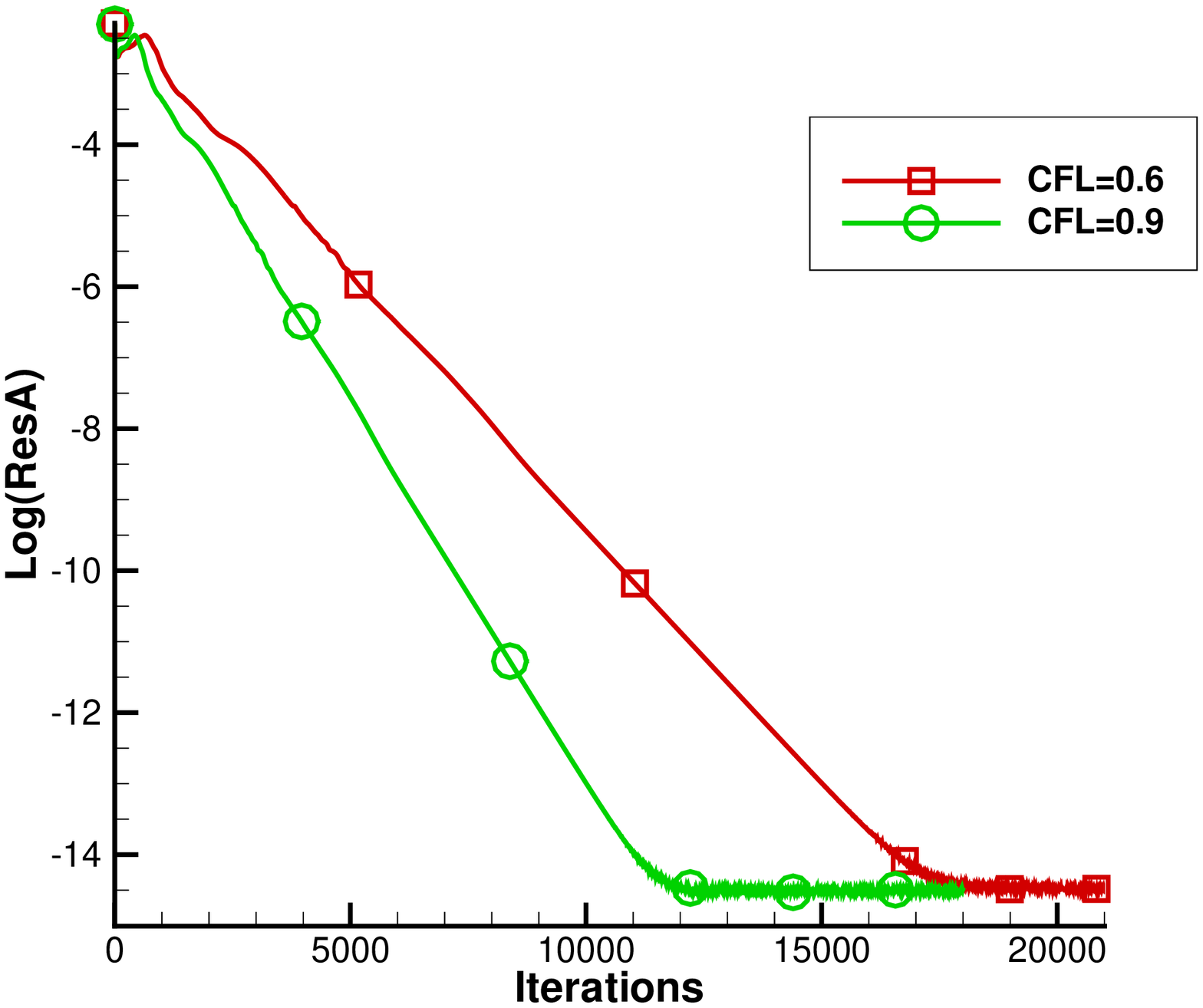}
%\caption{fig1}
\end{minipage}%
}%
\subfigure[FE fast sweeping, $Ma=2$]{
\begin{minipage}[t]{0.3\linewidth}
\centering
\includegraphics[width=2.0in]{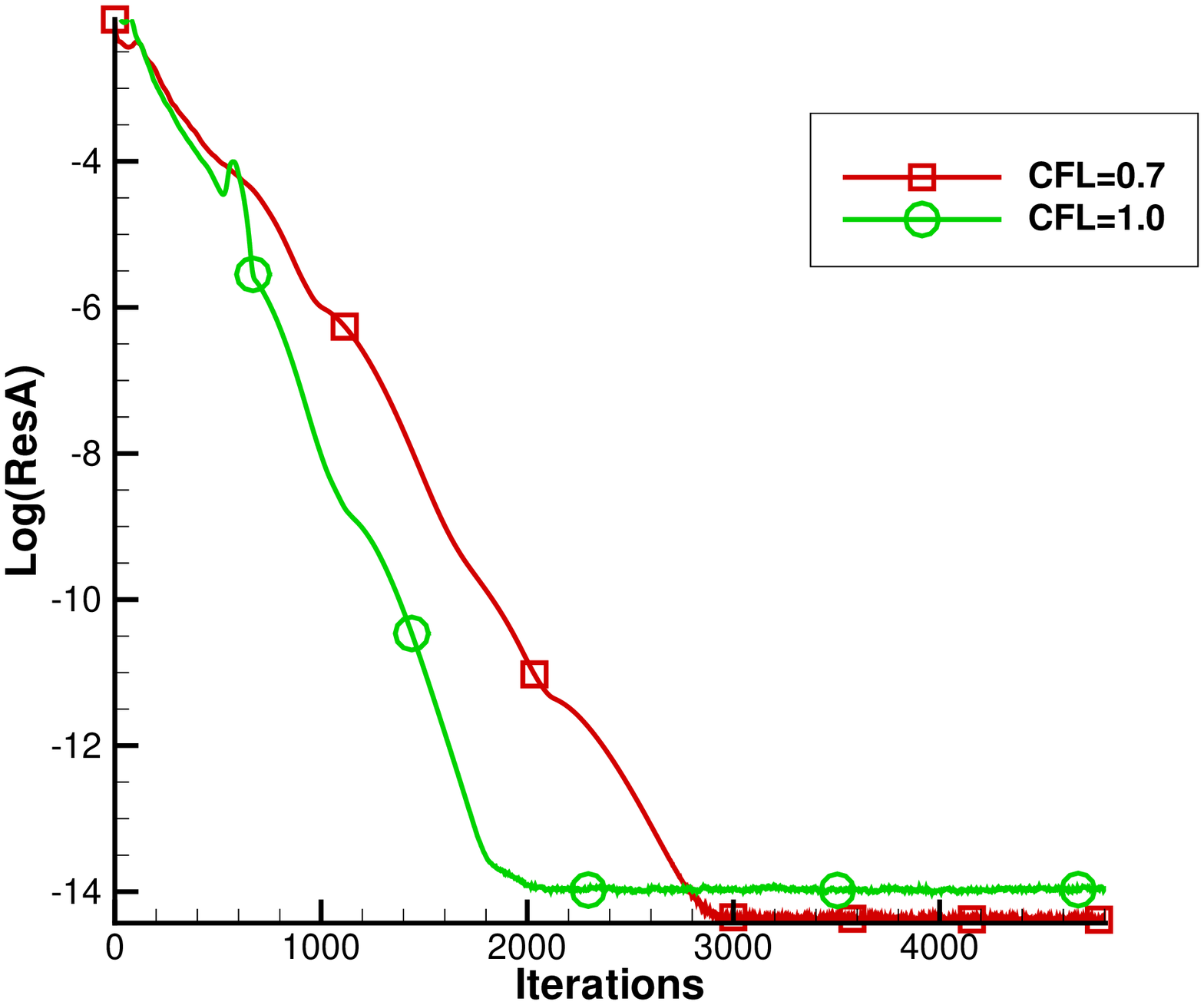}
%\caption{fig2}
\end{minipage}%
}%
\newline
\subfigure[FE Jacobi, $Ma=0.8$]{
\begin{minipage}[t]{0.3\linewidth}
\centering
\includegraphics[width=2.0in]{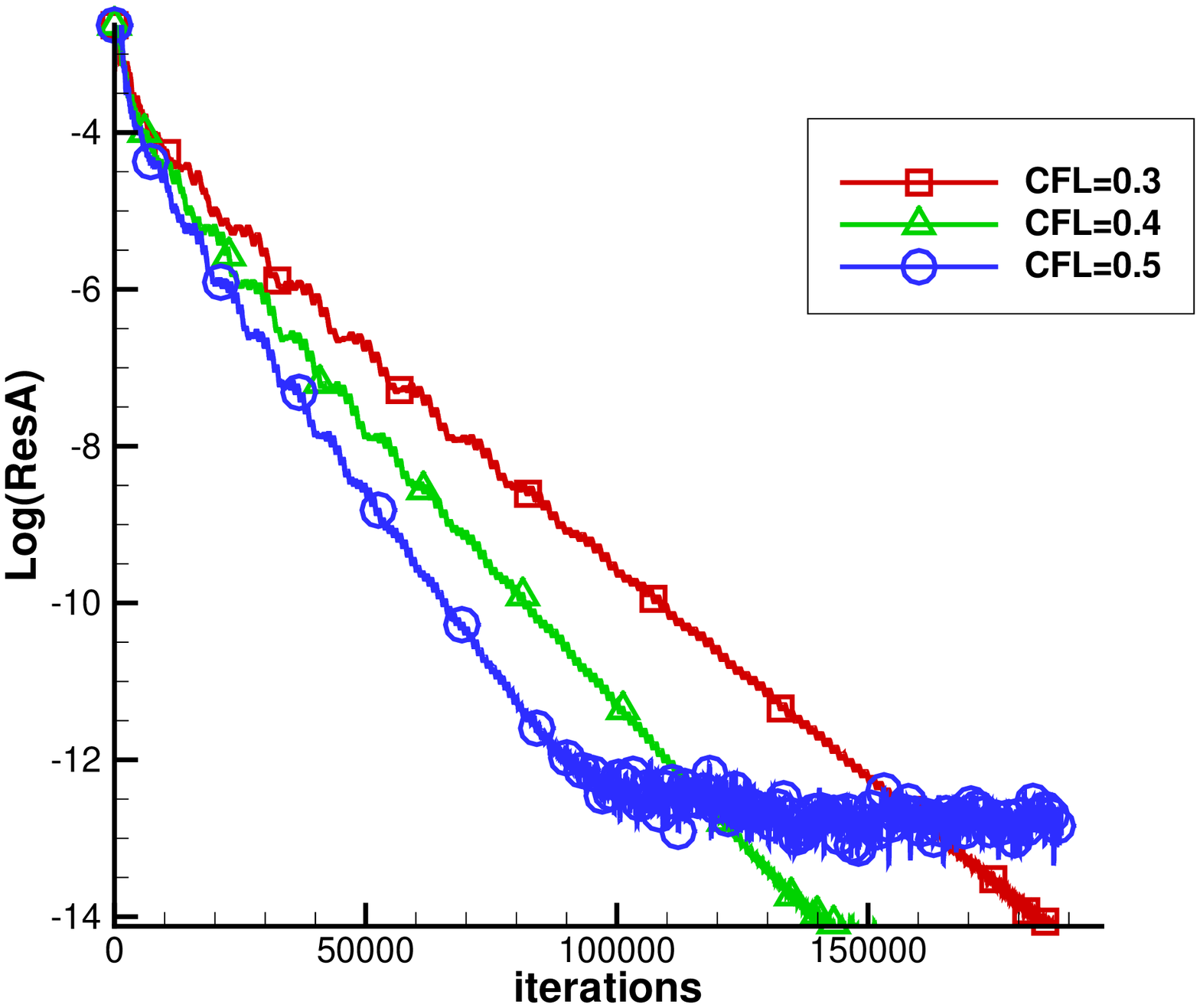}
%\caption{fig1}
\end{minipage}%
}%
\subfigure[RK Jacobi, $Ma=0.8$]{
\begin{minipage}[t]{0.3\linewidth}
\centering
\includegraphics[width=2.0in]{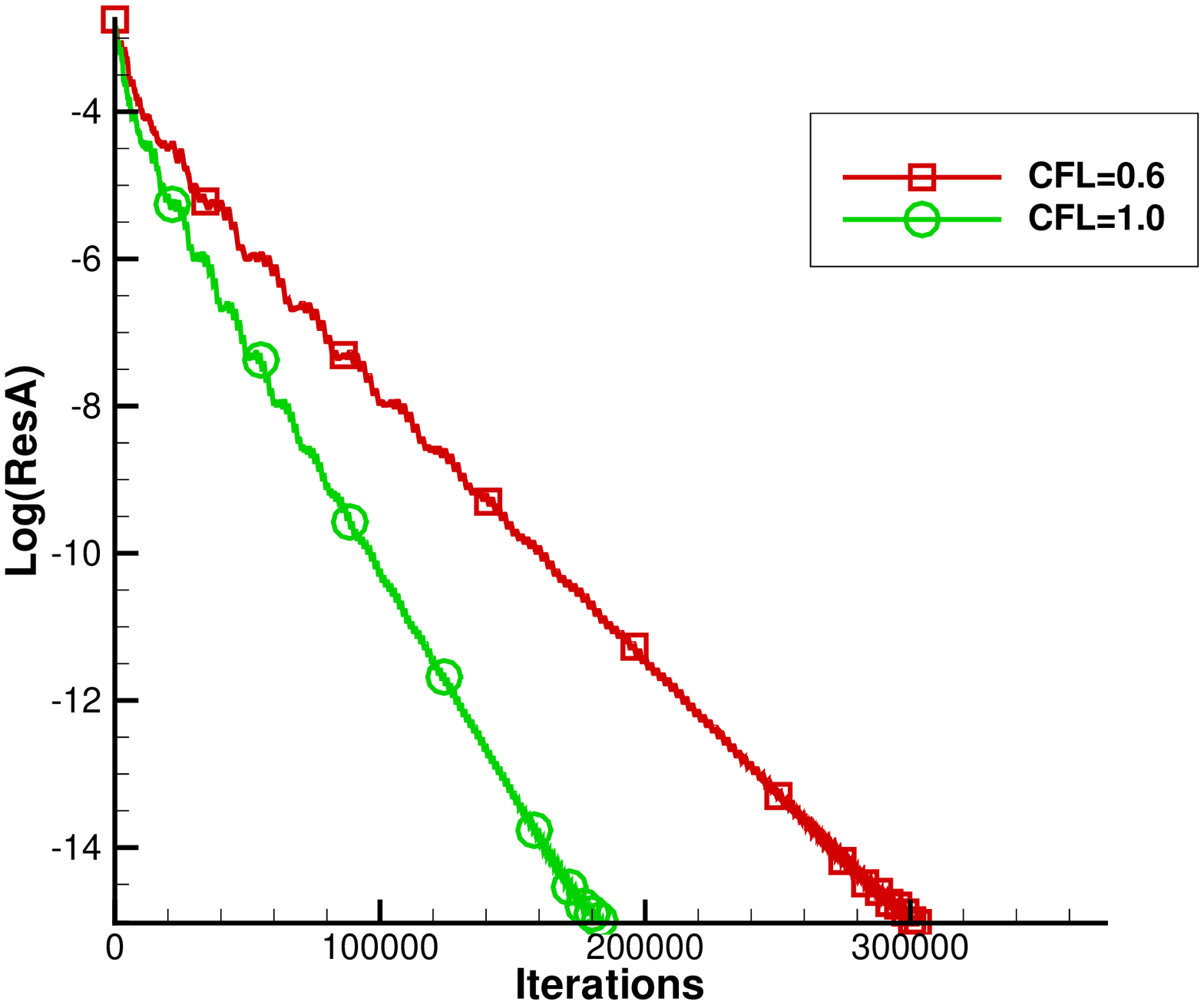}
%\caption{fig1}
\end{minipage}%
}%
\subfigure[FE fast sweeping, $Ma=0.8$]{
\begin{minipage}[t]{0.3\linewidth}
\centering
\includegraphics[width=2.0in]{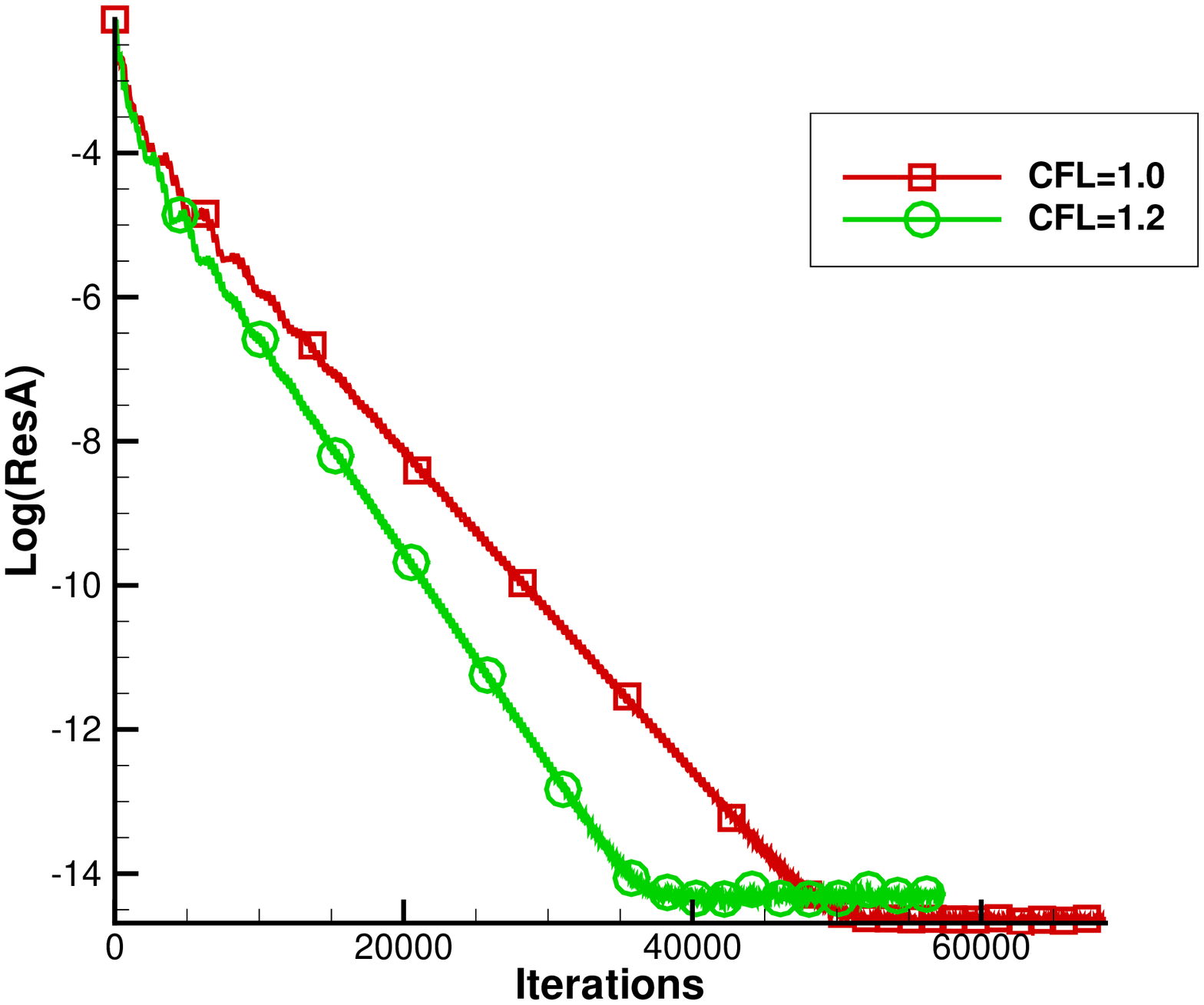}
%\caption{fig2}
\end{minipage}%
}%
\centering
\caption{\label{9.4}Example 5, supersonic and subsonic flows past an NACA001035 airfoil. The convergence history of the residue as a function of number of iterations for three schemes with different CFL numbers. (a), (b), (c): for the case of  $Ma=2$, angle of attack $\alpha=1^{\circ}$; (d) (e) (f): for the case of  $Ma=0.8$, angle of attack $\alpha=1.25^{\circ}$.}
\end{figure}

\noindent In this example, we solve the problem of supersonic and subsonic flows past a single NACA001035 airfoil configuration in \cite{CMMM, ZS3}. Following the setup in \cite{ZS3}, we consider the case of supersonic flow with Mach number $Ma=2$, angle of attack $\alpha=1^{\circ}$; and the case of subsonic flow with Mach number $Ma=0.8$, angle of attack $\alpha=1.25^{\circ}$. The computation domain is $[-16,16]\times[-16,16]$. The unstructured mesh containing 5593 triangles used for this example is shown in Figure \ref{9.1}. Four corners of the computational domain are taken as the reference points to form the alternating sweeping directions in the FE fast sweeping scheme. In Table \ref{9.2}, number of iterations required to reach the convergence criterion threshold value $10^{-11}$,  and total CPU time when the schemes converge under various CFL numbers are reported for the FE Jacobi scheme, the RK Jacobi scheme, and the FE fast sweeping scheme. In this example, the CFL number constraint for the FE Jacobi scheme to converge is not as severe as that  in the examples 1, 2 and 3. Using the largest CFL number permitted and due to its simple one-stage structure, the FE Jacobi scheme takes less number of iterations and less CPU times to converge to steady states than the RK Jacobi scheme, for both the supersonic flow case and the subsonic flow case. Again, the FE fast sweeping scheme is the most efficient iterative method among three methods. Here it permits much larger CFL numbers than the FE Jacobi scheme, and also slightly larger CFL numbers than the RK Jacobi scheme. It also has a simple one-stage structure. With the largest CFL number allowed for the methods to reach steady state solution, the FE fast sweeping method on unstructured triangular meshes saves about $60\% \sim 65\%$ CPU time cost of that by the RK Jacobi scheme (the TVD-RK3 scheme) for both the supersonic flow case and the subsonic flow case. In Figure \ref{9.3}, the pressure contours of the converged steady state solutions of these three schemes are presented for the supersonic flow case and the subsonic flow case. We observe  comparable numerical steady states of these schemes. The residue history of these three schemes with different CFL numbers is reported in Figure \ref{9.4}, which shows that the residue of iterations can settle down to tiny values around $10^{-14}$ to $10^{-15}$, at the level of round off errors. Again, the absolute convergence of the proposed high order fast sweeping method on triangular meshes is verified in this example.

\begin{table}%[H]
		\centering
\begin{tabular}{|c|c|c|}\hline

			\multicolumn{3}{|c|}{\bf{supersonic,} \bm{$Ma=2,\alpha=1^{\circ}$} }\\\hline
			\multicolumn{3}{|c|}{FE Jacobi scheme }\\\hline
             CFL number & iteration number  & CPU time \\\hline
	0.3& 8118 & 2527.66 \\\hline
	0.4& 6077	& 1904.80 \\\hline
      0.5   & Not convergent        & - \\
    \hline
		\end{tabular}

\begin{tabular}{|c|c|c|}\hline
			\multicolumn{3}{|c|}{RK Jacobi scheme }\\\hline
             CFL number & iteration number & CPU time \\\hline
	0.6&  12192 &  3709.35 \\\hline
0.9&  8130 	& 2488.23 \\\hline
           1.0   & Not convergent        & - \\
    \hline
		\end{tabular}

\begin{tabular}{|c|c|c|}\hline
			\multicolumn{3}{|c|}{FE fast sweeping scheme }\\\hline
             CFL number & iteration number  & CPU time \\\hline
     	0.7& 2031 	& 1183.33 \\\hline
     1.0& 1508 	& 867.25 \\\hline
       1.1   & Not convergent      & - \\
    \hline
		\end{tabular}

\begin{tabular}{|c|c|c|}\hline
			\multicolumn{3}{|c|}{\bf{subsonic,} \bm{$Ma=0.8,\alpha=1.25^{\circ}$} }\\\hline
			\multicolumn{3}{|c|}{FE Jacobi scheme }\\\hline
             CFL number & iteration number  & CPU time \\\hline
	0.4& 95877 &  31918.58 \\\hline
	0.5& 76794 &  24608.73  \\\hline
      0.6   & Not convergent       & - \\
    \hline
		\end{tabular}

\begin{tabular}{|c|c|c|}\hline
			\multicolumn{3}{|c|}{RK Jacobi scheme }\\\hline
             CFL number & iteration number  & CPU time \\\hline
	0.6&185661  &   63876.05 \\\hline
1.0&111396 	& 37444.54 \\\hline
           1.1   & Not convergent        & - \\
    \hline
		\end{tabular}

\begin{tabular}{|c|c|c|}\hline
			\multicolumn{3}{|c|}{FE fast sweeping scheme }\\\hline
             CFL number & iteration number & CPU time \\\hline
     	1.0& 32854 &	22258.67 \\\hline
      	1.2& 24876 	& 15613.47 \\\hline
       1.3   & Not convergent        & - \\
    \hline
		\end{tabular}
		\caption{\label{9.2}Example 5, supersonic and subsonic flows past an NACA001035 airfoil.  Number of iterations and total CPU time when convergence is obtained. Convergence criterion threshold value  is $10^{-11}$. CPU time unit: second.}
%		\label{tab:Margin_settings}
	\end{table}

~\\
\noindent\textbf{Example 6. Supersonic flows past an NACA0012 airfoil}

\begin{figure}%[H]
\centering
\subfigure[]{
\centering
\includegraphics[width=2.4in]{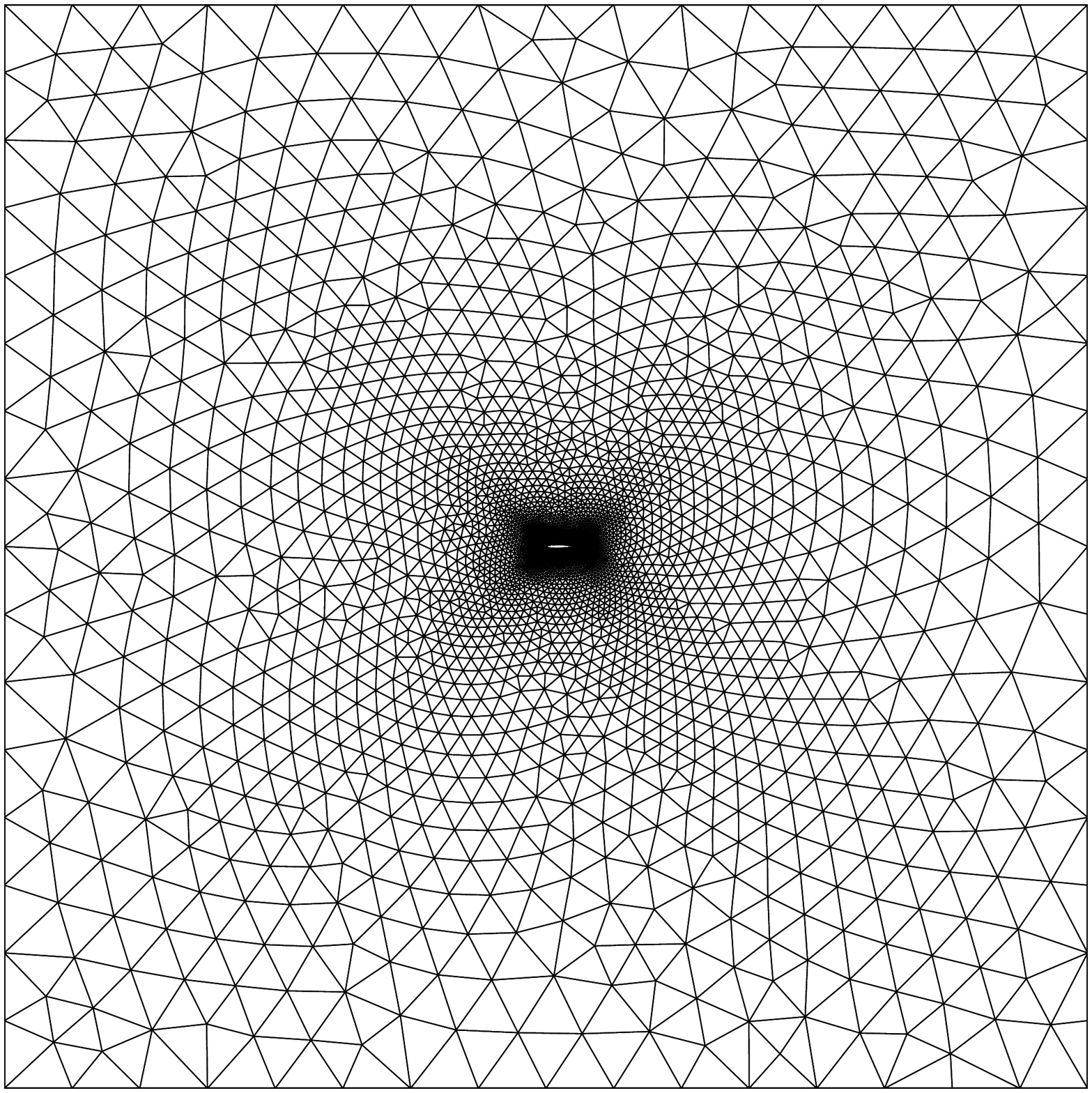}
}
\subfigure[]{
\centering
\includegraphics[width=2.4in]{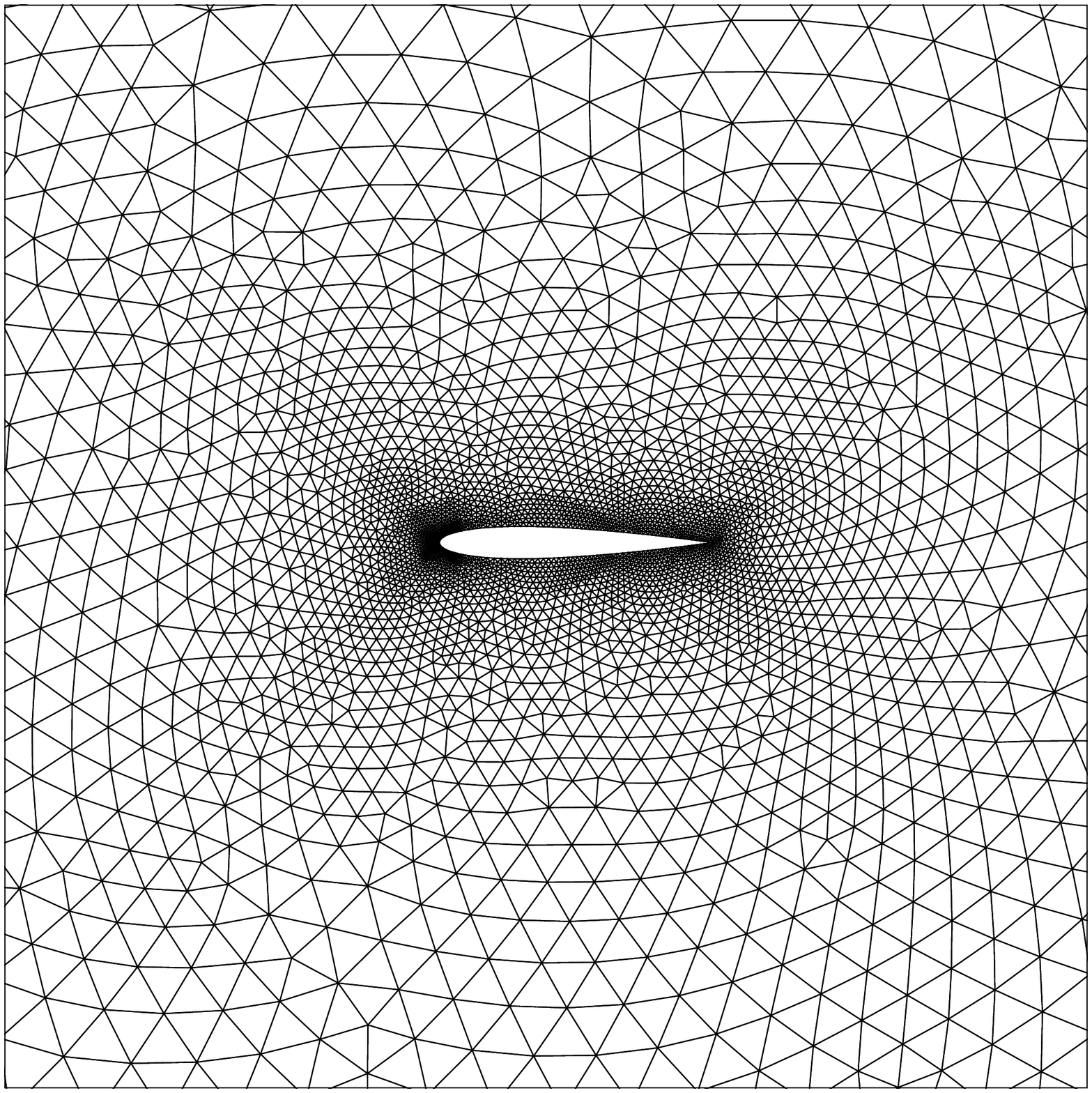}
}
\caption{\label{3.1} The computational mesh for Example 6 and Example 7, supersonic and subsonic flows past an NACA0012 airfoil. Left: the whole domain; right: zoomed region  near the airfoil.}
\end{figure}

\begin{figure}%[H]
\centering
\subfigure[FE Jacobi, $Ma=3$]{
\begin{minipage}[t]{0.3\linewidth}
\centering
\includegraphics[width=2.0in]{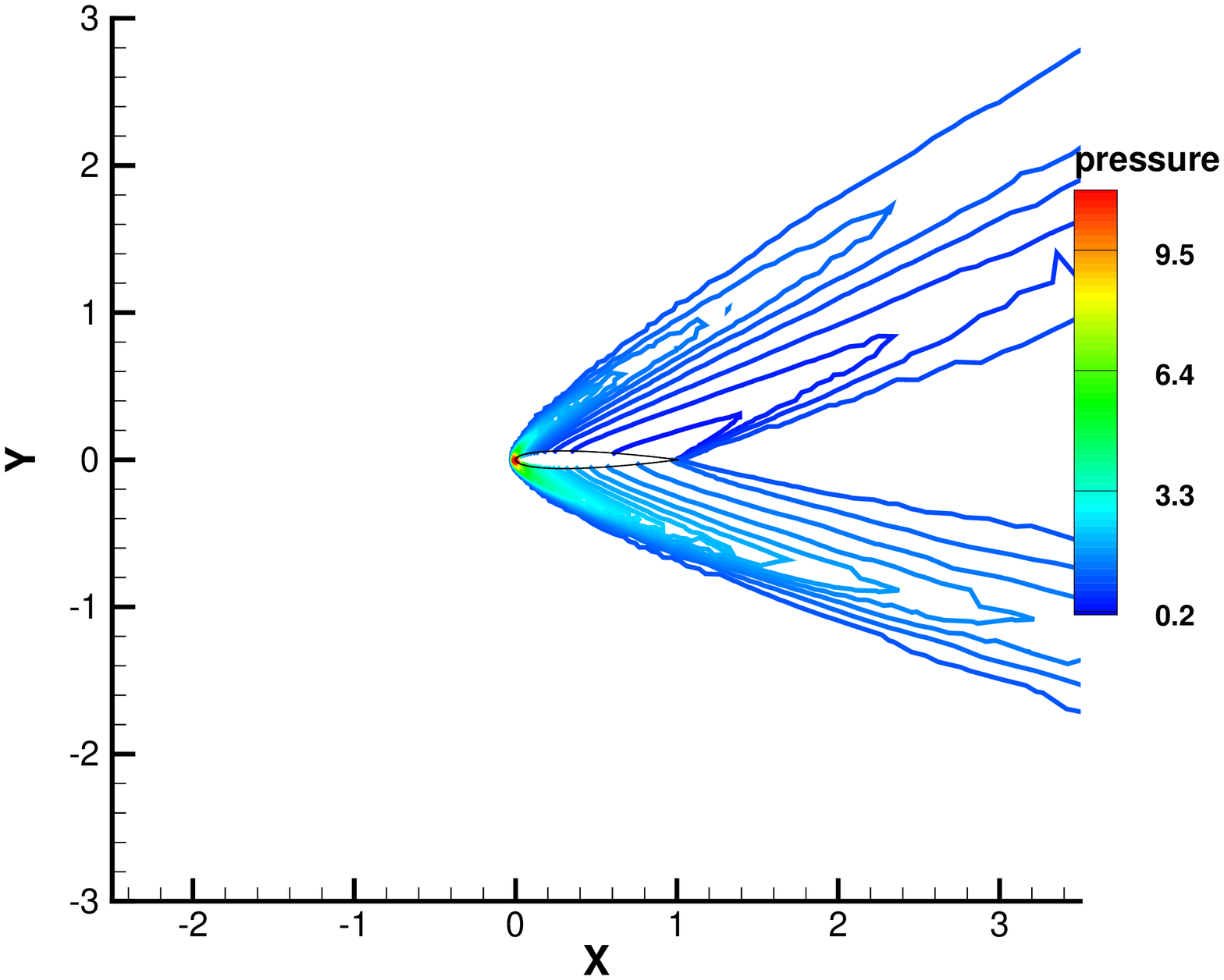}
%\caption{fig1}
\end{minipage}%
}%
\subfigure[RK Jacobi, $Ma=3$]{
\begin{minipage}[t]{0.3\linewidth}
\centering
\includegraphics[width=2.0in]{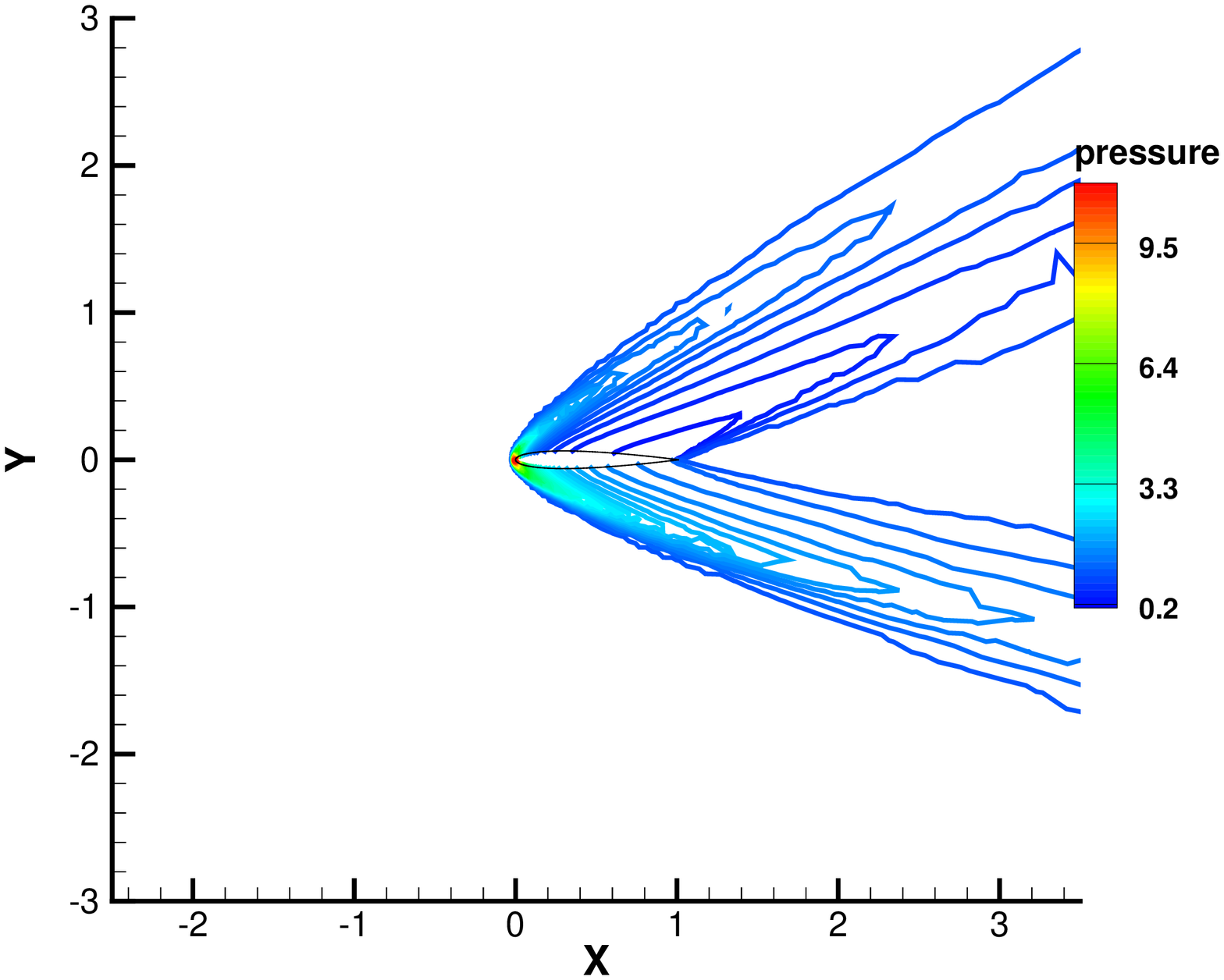}
%\caption{fig1}
\end{minipage}%
}%
\subfigure[FE fast sweeping, $Ma=3$]{
\begin{minipage}[t]{0.3\linewidth}
\centering
\includegraphics[width=2.0in]{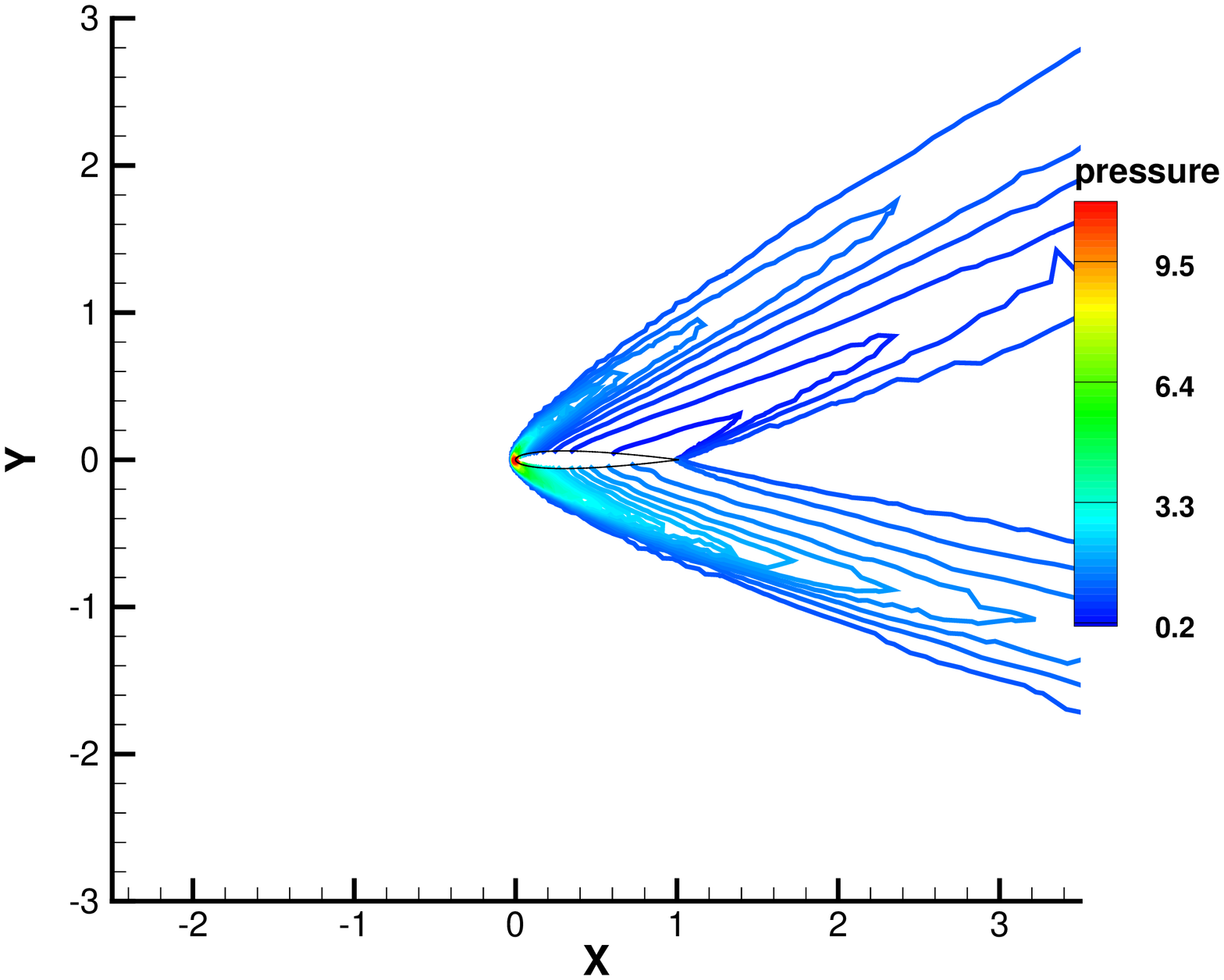}
%\caption{fig2}
\end{minipage}%
}%
\newline
\subfigure[FE Jacobi, $Ma=2$]{
\begin{minipage}[t]{0.3\linewidth}
\centering
\includegraphics[width=2.0in]{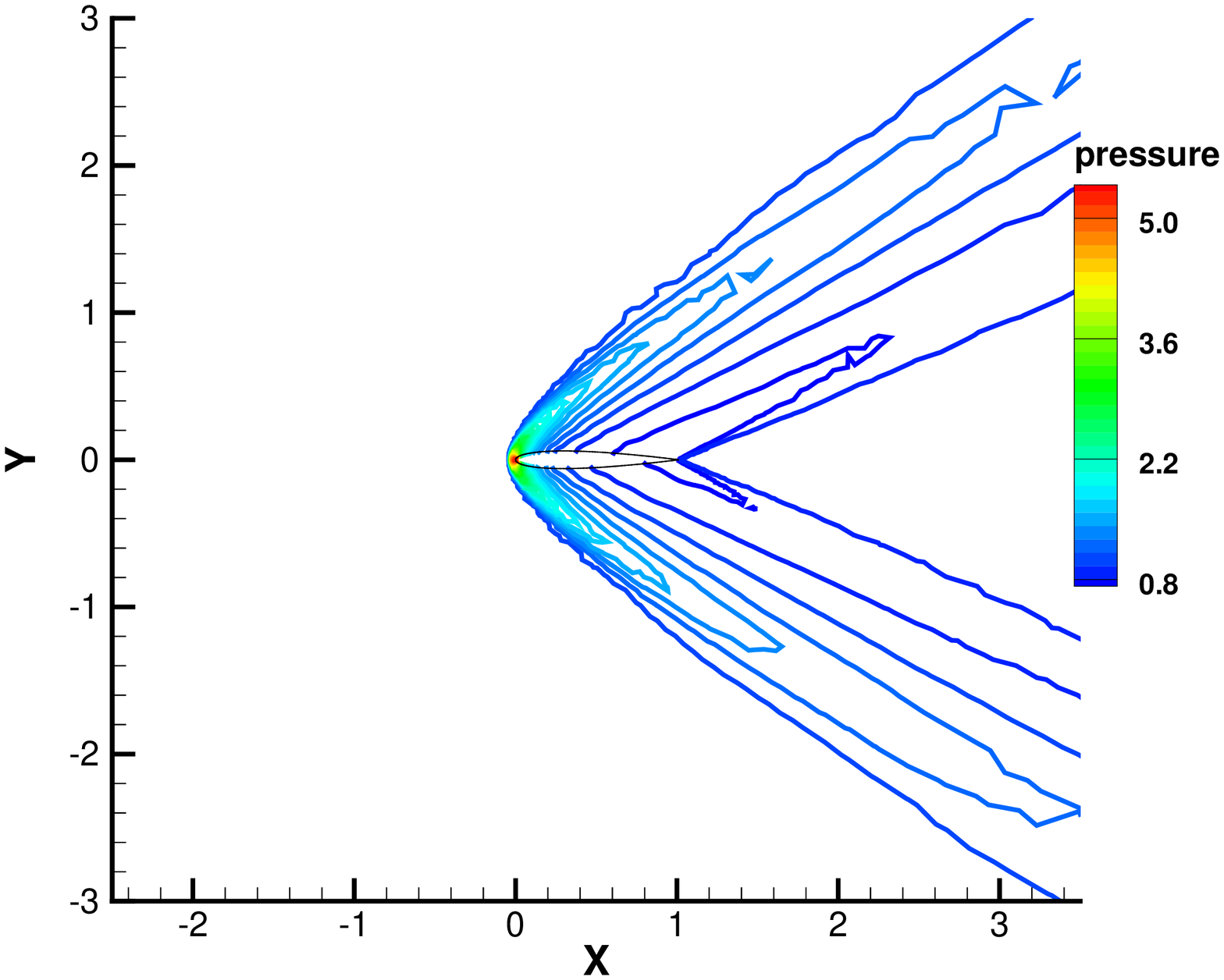}
%\caption{fig1}
\end{minipage}%
}%
\subfigure[RK Jacobi, $Ma=2$]{
\begin{minipage}[t]{0.3\linewidth}
\centering
\includegraphics[width=2.0in]{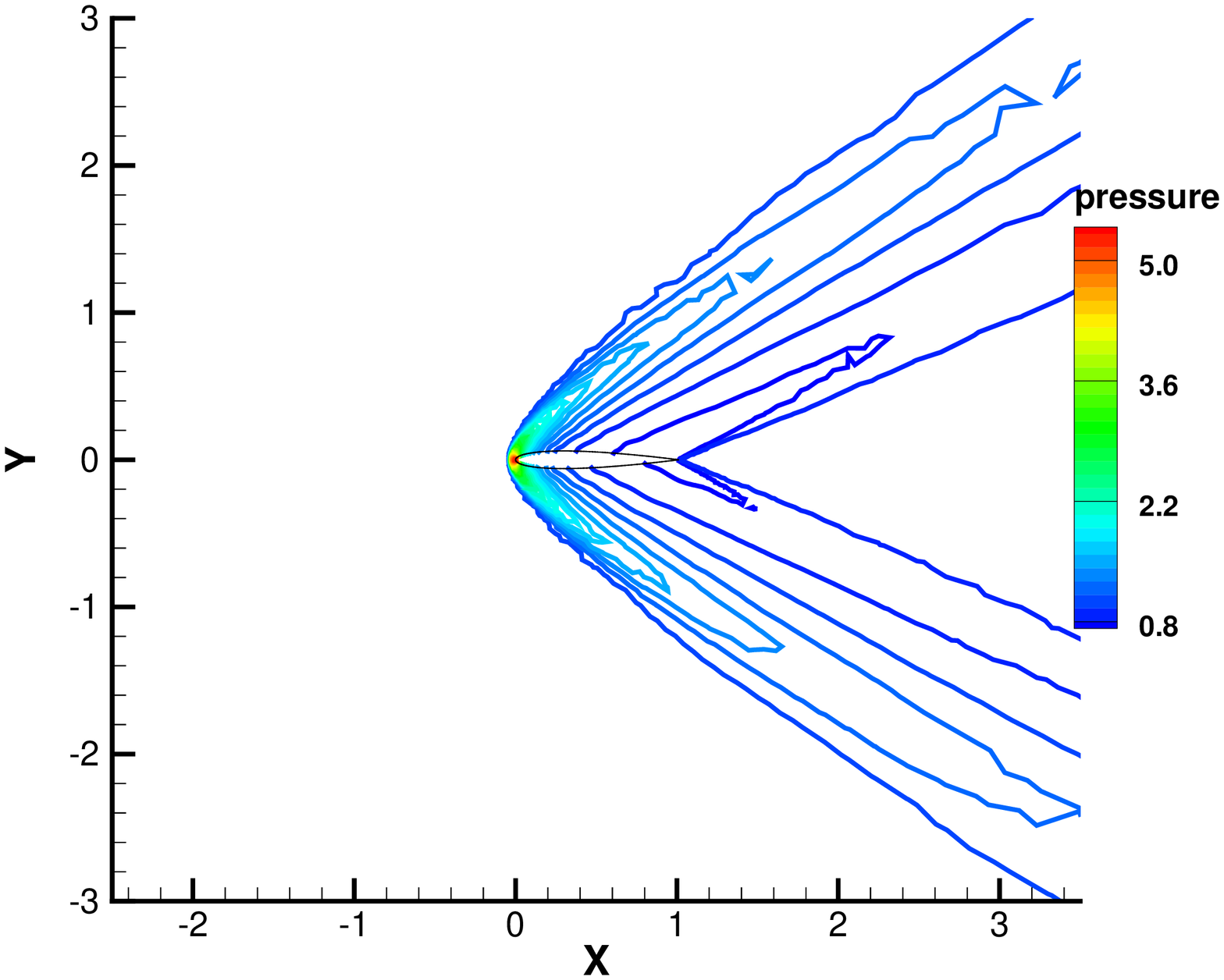}
%\caption{fig1}
\end{minipage}%
}%
\subfigure[FE fast sweeping, $Ma=2$]{
\begin{minipage}[t]{0.3\linewidth}
\centering
\includegraphics[width=2.0in]{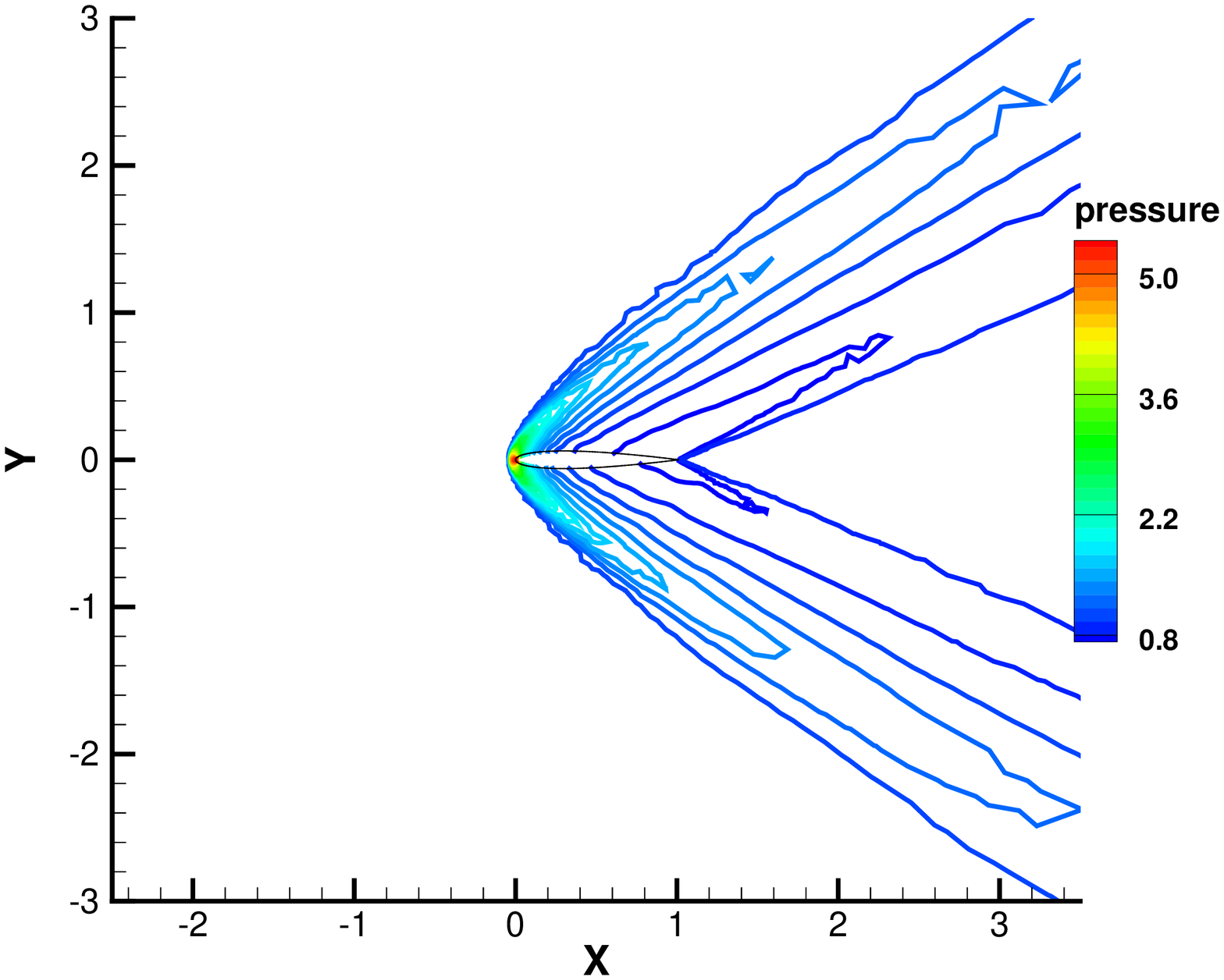}
%\caption{fig2}
\end{minipage}%
}%
\centering
\caption{\label{3.3}Example 6, supersonic flows past an NACA0012 airfoil. The converged steady states of numerical solutions by three different iterative schemes. (a) (b) (c): 30 equally spaced pressure contours from 0.2 to 11 for the case of  $Ma=3$, angle of attack $\alpha=10^{\circ}$; (d) (e) (f): 30 equally spaced pressure contours from 0.8 to 5.2 for the case of $Ma=2$, angle of attack $\alpha=1^{\circ}$.}
\end{figure}

\begin{figure}%[H]
\centering
\subfigure[FE Jacobi, $Ma=3$]{
\begin{minipage}[t]{0.3\linewidth}
\centering
\includegraphics[width=2.0in]{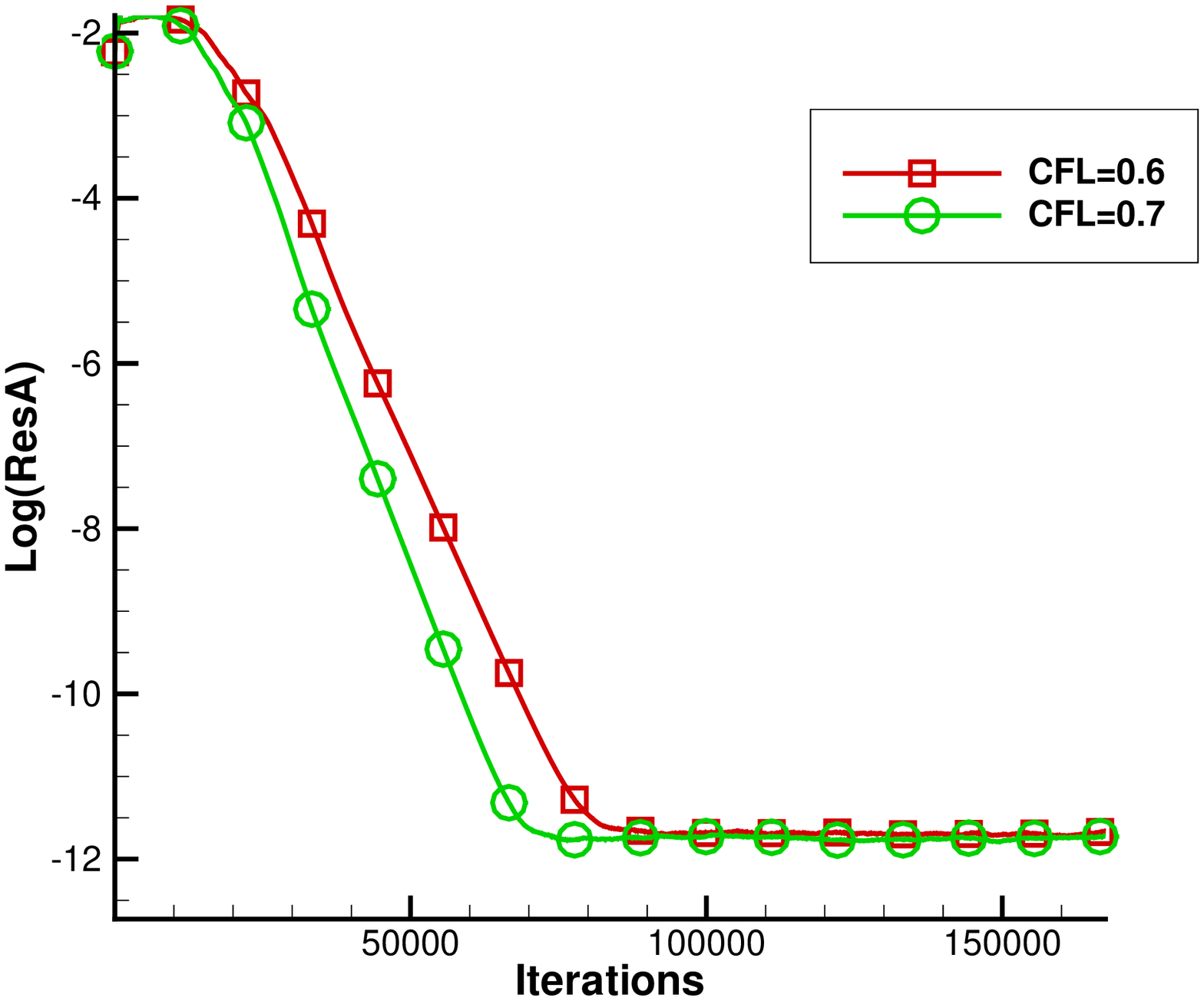}
%\caption{fig1}
\end{minipage}%
}%
\subfigure[RK Jacobi, $Ma=3$]{
\begin{minipage}[t]{0.3\linewidth}
\centering
\includegraphics[width=2.0in]{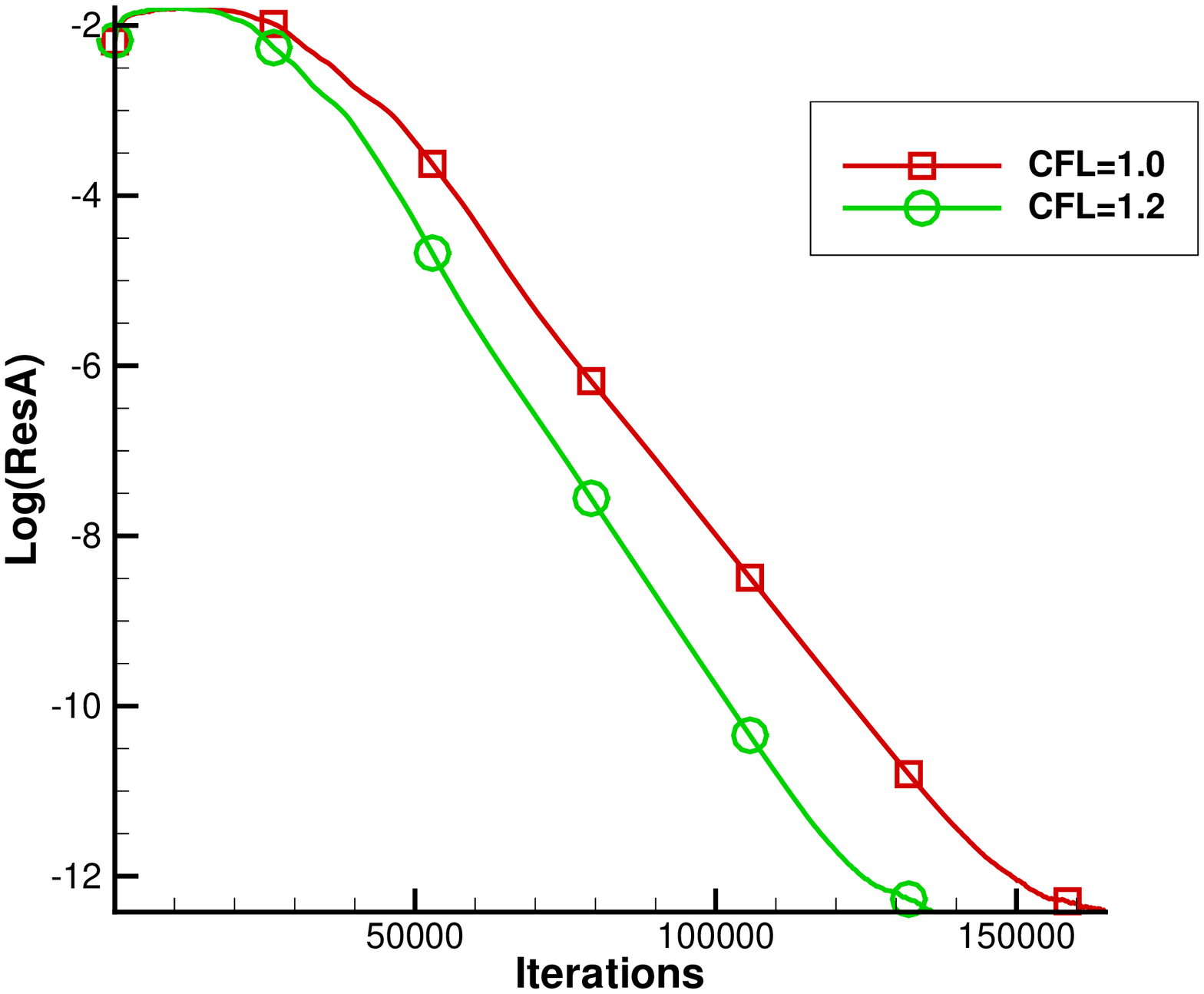}
%\caption{fig1}
\end{minipage}%
}%
\subfigure[FE fast sweeping, $Ma=3$]{
\begin{minipage}[t]{0.3\linewidth}
\centering
\includegraphics[width=2.0in]{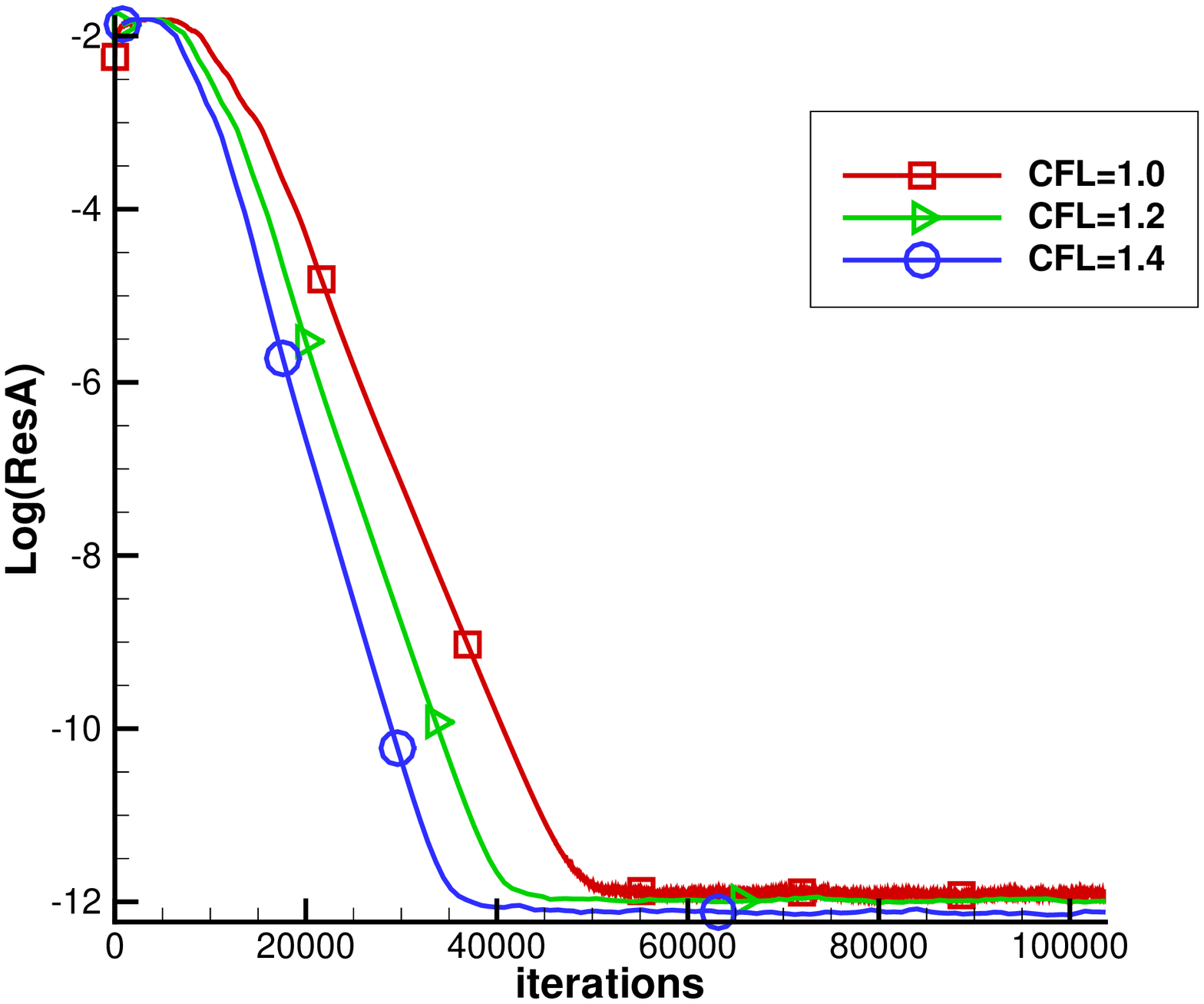}
%\caption{fig2}
\end{minipage}%
}%
\newline
\subfigure[FE Jacobi, $Ma=2$]{
\begin{minipage}[t]{0.3\linewidth}
\centering
\includegraphics[width=2.0in]{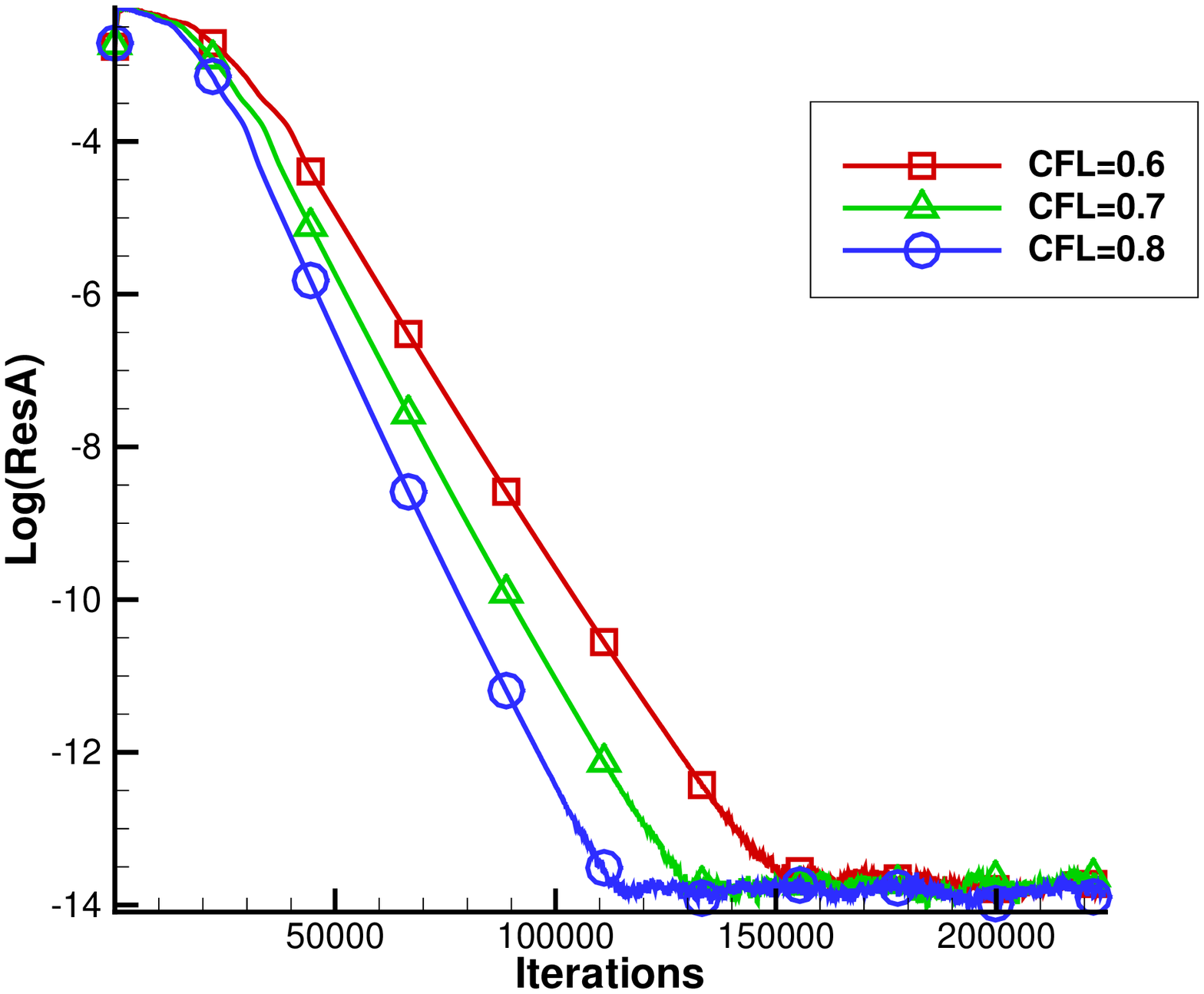}
%\caption{fig1}
\end{minipage}%
}%
\subfigure[RK Jacobi, $Ma=2$]{
\begin{minipage}[t]{0.3\linewidth}
\centering
\includegraphics[width=2.0in]{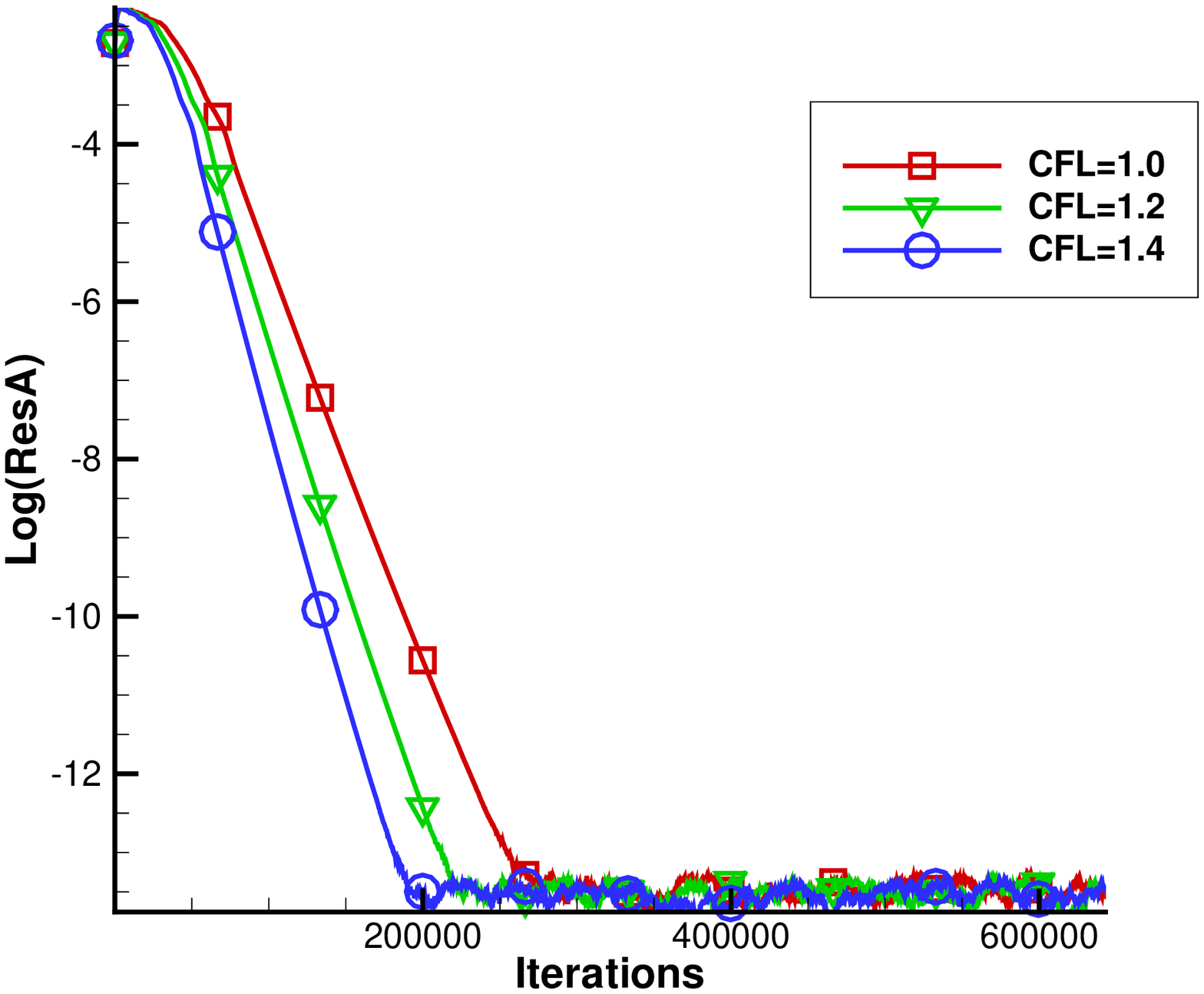}
%\caption{fig1}
\end{minipage}%
}%
\subfigure[FE fast sweeping, $Ma=2$]{
\begin{minipage}[t]{0.3\linewidth}
\centering
\includegraphics[width=2.0in]{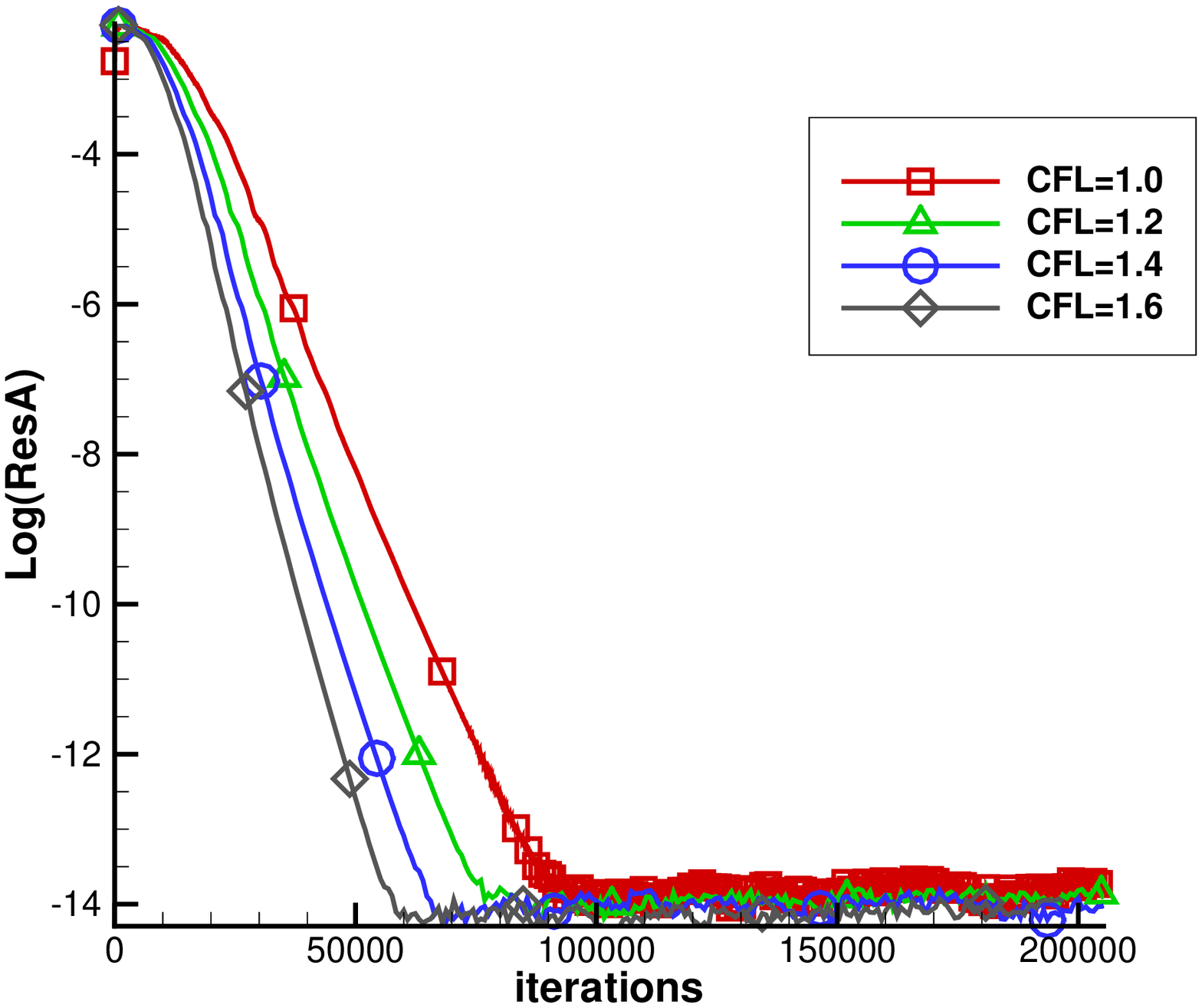}
%\caption{fig2}
\end{minipage}%
}%
\centering
\caption{\label{3.4}Example 6, supersonic flows past an NACA0012 airfoil. The convergence history of the residue as a function of number of iterations for three schemes with different CFL numbers. (a), (b), (c): for the case of $Ma=3$, angle of attack $\alpha=10^{\circ}$; (d) (e) (f): for the case of $Ma=2$, angle of attack $\alpha=1^{\circ}$.}
\end{figure}

\noindent In this example, we consider the problem of inviscid Euler supersonic flows past a single NACA0012 airfoil configuration in \cite{LBL}. Following the setup in \cite{ZS3}, we consider the case of flow with Mach number $Ma=3$, angle of attack $\alpha=10^{\circ}$; and the case of flow with Mach number $Ma=2$, angle of attack $\alpha=1^{\circ}$. The computational domain is $[-15,15]\times[-15,15]$. The unstructured mesh for this example, which consists of 9340 triangles, is shown in Figure \ref{3.1}. As the previous example, four corners of the computational domain are taken as the reference points to form the alternating sweeping directions in the FE fast sweeping scheme. In Table \ref{3.2}, number of iterations required to reach the convergence criterion threshold value $10^{-11}$,  and total CPU time when the schemes converge under various CFL numbers are reported for the FE Jacobi scheme, the RK Jacobi scheme, and the FE fast sweeping scheme. Similar to Example 5, the CFL number constraint for the FE Jacobi scheme to converge is not as severe as that  in the examples 1, 2 and 3. Due to the simple one-stage structure of the FE Jacobi scheme, it takes the FE Jacobi scheme less number of iterations and less CPU times to converge to steady states than the RK Jacobi scheme, for both cases with different Mach numbers. The FE fast sweeping scheme is still the most efficient iterative method among three methods. It allows much larger CFL numbers than the FE Jacobi scheme, and even slightly larger CFL numbers than the RK Jacobi scheme. It also has a simple one-stage structure. With the largest CFL number permitted in each method to reach steady state solution, the FE fast sweeping method on unstructured triangular meshes saves $50\% \sim 55\%$ CPU time cost of that by the RK Jacobi scheme (the TVD-RK3 scheme) for both supersonic flow cases in this example.
In Figure \ref{3.3}, the pressure contours of the converged steady state solutions of these three schemes are presented for both supersonic flow cases. Again, we observe  comparable numerical steady states of these schemes. Figure \ref{3.4}
shows the residue history of these three schemes with different CFL numbers. It is observed that the residue of iterations can settle down to tiny values at the level of round off errors, which verifies the absolute convergence of the developed high order fast sweeping method on triangular meshes here.

\begin{table}%[H]
		\centering
\begin{tabular}{|c|c|c|}\hline
			\multicolumn{3}{|c|}{\bm{$Ma=3,\alpha=10^{\circ}$}  }\\\hline
			\multicolumn{3}{|c|}{FE Jacobi scheme }\\\hline
             CFL number & iteration number  & CPU time \\\hline
	0.6& 75193	& 27432.75 \\\hline
0.7& 64349 	& 23534.16 \\\hline
    0.8   & Not convergent        & - \\
    \hline
		\end{tabular}

\begin{tabular}{|c|c|c|}\hline
			\multicolumn{3}{|c|}{RK Jacobi scheme }\\\hline
             CFL number & iteration number& CPU time \\\hline
	1.0& 134485 &  49427.91 \\\hline
    1.2   &  112165      &42557.77 \\\hline
     1.3   & Not convergent        & - \\
    \hline
		\end{tabular}

\begin{tabular}{|c|c|c|}\hline
			\multicolumn{3}{|c|}{FE fast sweeping scheme }\\\hline
             CFL number & iteration number  & CPU time \\\hline
             1.0   & 44688    &28805.77 \\\hline
        1.2   &37168    &23178.72 \\\hline
         1.4   &31816     &19563.90 \\ \hline
             1.5   &Not convergent     & - \\
    \hline
		\end{tabular}

\begin{tabular}{|c|c|c|}\hline
			\multicolumn{3}{|c|}{\bm{$Ma=2,\alpha=1^{\circ}$} }\\\hline
			\multicolumn{3}{|c|}{FE Jacobi scheme }\\\hline
             CFL number & iteration number &CPU time \\\hline
	0.6& 116201 &  47048.46 \\\hline
0.7& 99598 &  40294.13 \\\hline
    0.8   & 87149   & 36696.27\\\hline
      0.9   & Not convergent        & - \\
    \hline
		\end{tabular}

\begin{tabular}{|c|c|c|}\hline
			\multicolumn{3}{|c|}{RK Jacobi scheme }\\\hline
             CFL number & iteration number  & CPU time \\\hline
	1.0&  209200 &   82761.36 \\\hline
    1.2   &   174319 &  69020.84  \\\hline
     1.4   & 149425 &  61582.79\\\hline
          1.5   & Not convergent       & - \\
    \hline
		\end{tabular}

\begin{tabular}{|c|c|c|}\hline
			\multicolumn{3}{|c|}{FE fast sweeping scheme }\\\hline
             CFL number & iteration number & CPU time \\\hline
     	1.0& 68696 &   47259.59 \\\hline
    1.2   &   57208 & 38592.23  \\\hline
     1.4   & 49000 &   32435.17\\\hline
      1.6   & 42848 &   29783.56\\\hline
          1.7   & Not convergent       & - \\
    \hline
		\end{tabular}
		\caption{\label{3.2}Example 6, supersonic flows past an NACA0012 airfoil. Number of iterations and total CPU time when convergence is obtained. Convergence criterion threshold value   is $10^{-11}$. CPU time unit: second.}
%		\label{tab:Margin_settings}
	\end{table}

~\\
\noindent\textbf{Example 7. Subsonic flows past an NACA0012 airfoil}

\begin{figure}%[H]
\centering
\subfigure[FE Jacobi, $Ma=0.8$]{
\begin{minipage}[t]{0.3\linewidth}
\centering
\includegraphics[width=2.0in]{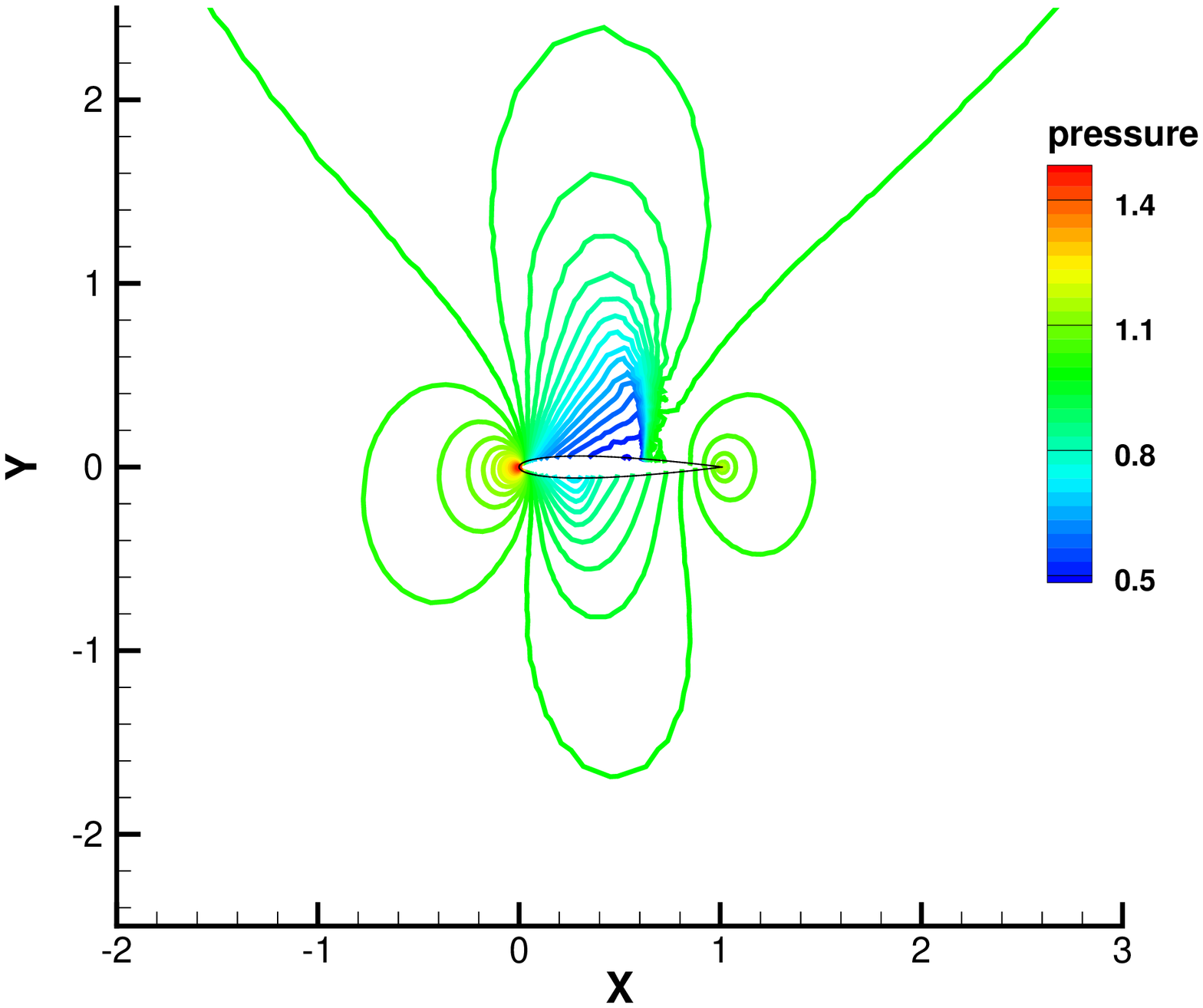}
%\caption{fig1}
\end{minipage}%
}%
\subfigure[RK Jacobi, $Ma=0.8$]{
\begin{minipage}[t]{0.3\linewidth}
\centering
\includegraphics[width=2.0in]{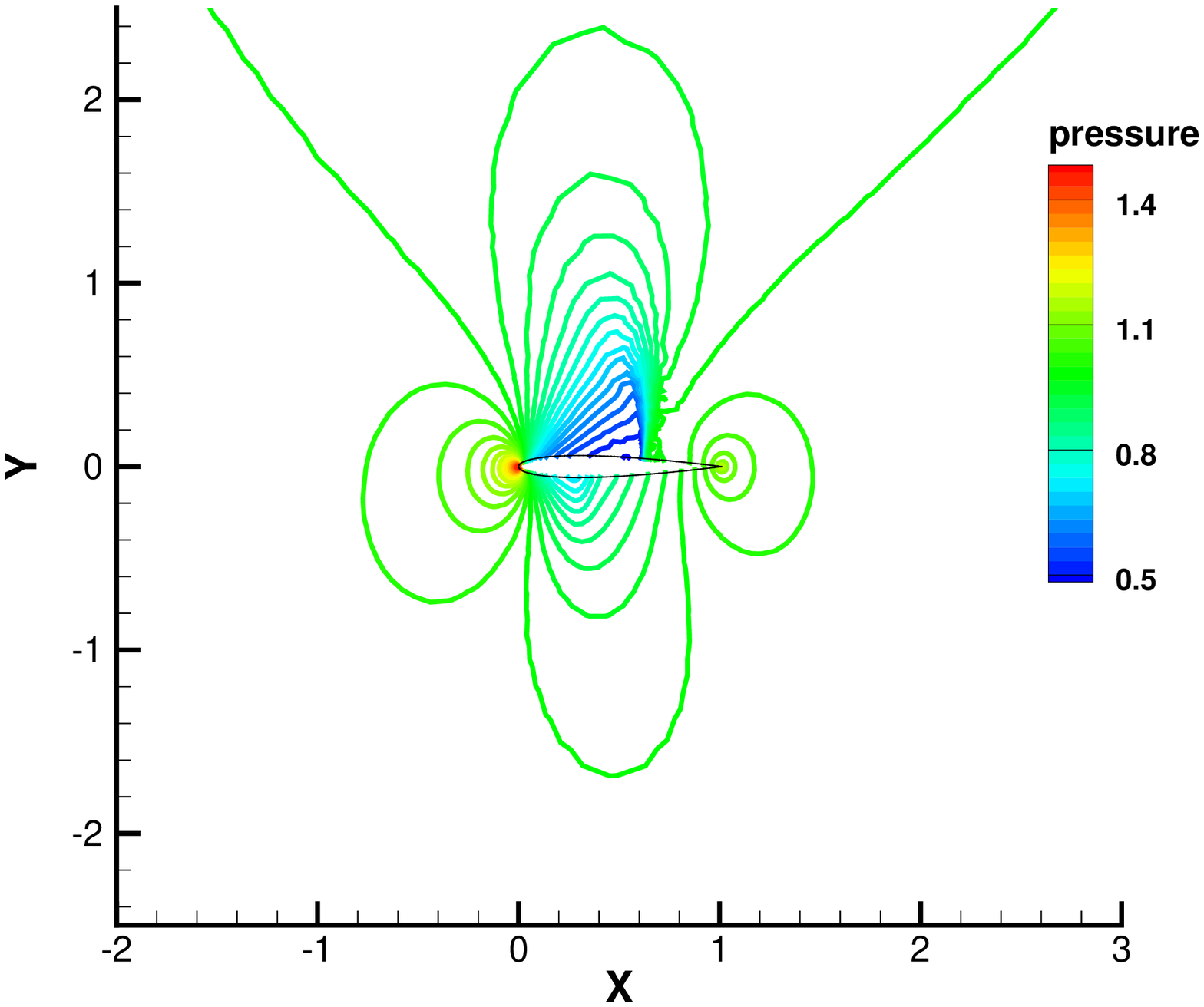}
%\caption{fig1}
\end{minipage}%
}%
\subfigure[FE fast sweeping, $Ma=0.8$]{
\begin{minipage}[t]{0.3\linewidth}
\centering
\includegraphics[width=2.0in]{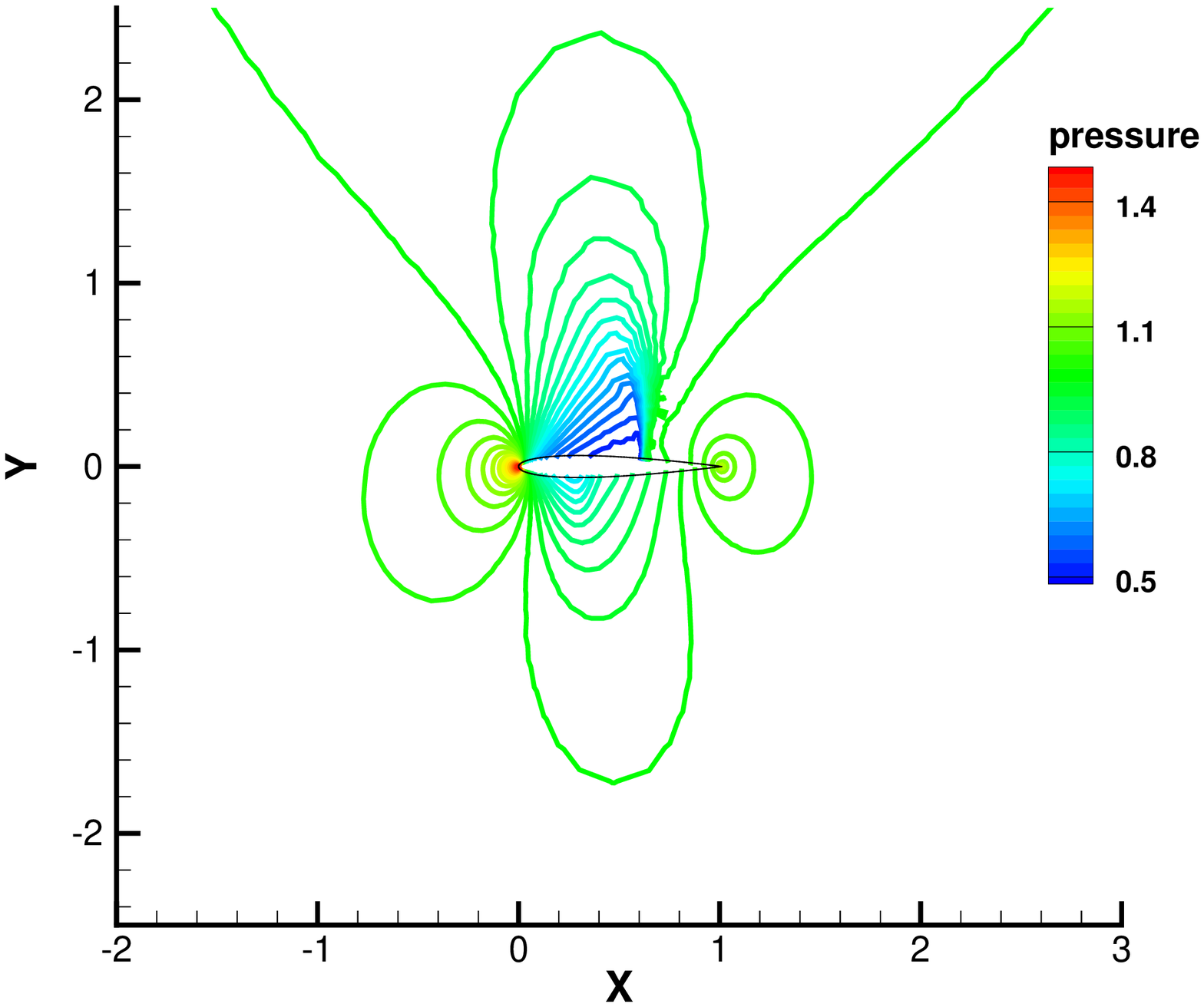}
%\caption{fig2}
\end{minipage}%
}%
\newline
\subfigure[FE Jacobi, $Ma=0.2$]{
\begin{minipage}[t]{0.3\linewidth}
\centering
\includegraphics[width=2.0in]{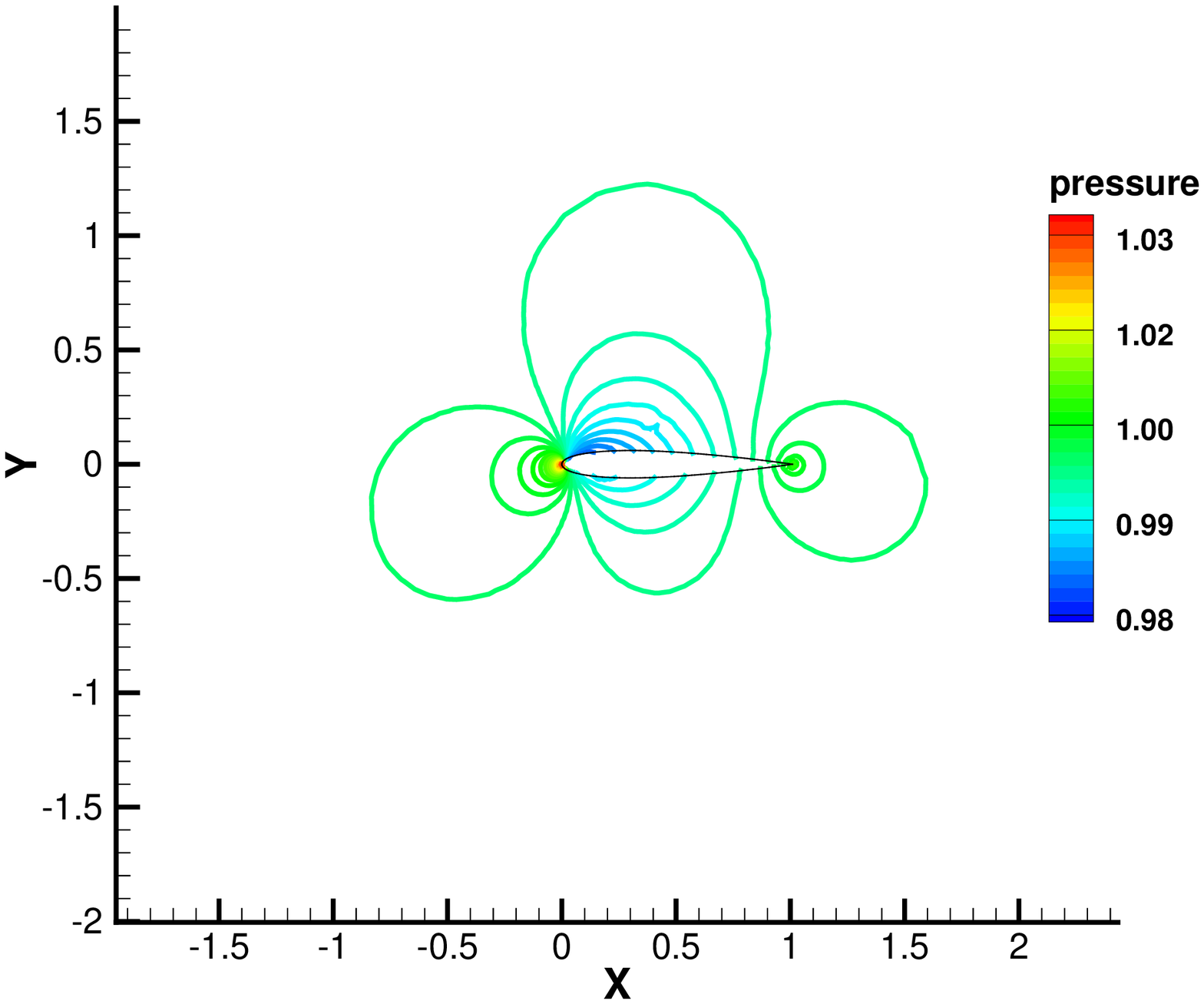}
%\caption{fig1}
\end{minipage}%
}%
\subfigure[RK Jacobi, $Ma=0.2$]{
\begin{minipage}[t]{0.3\linewidth}
\centering
\includegraphics[width=2.0in]{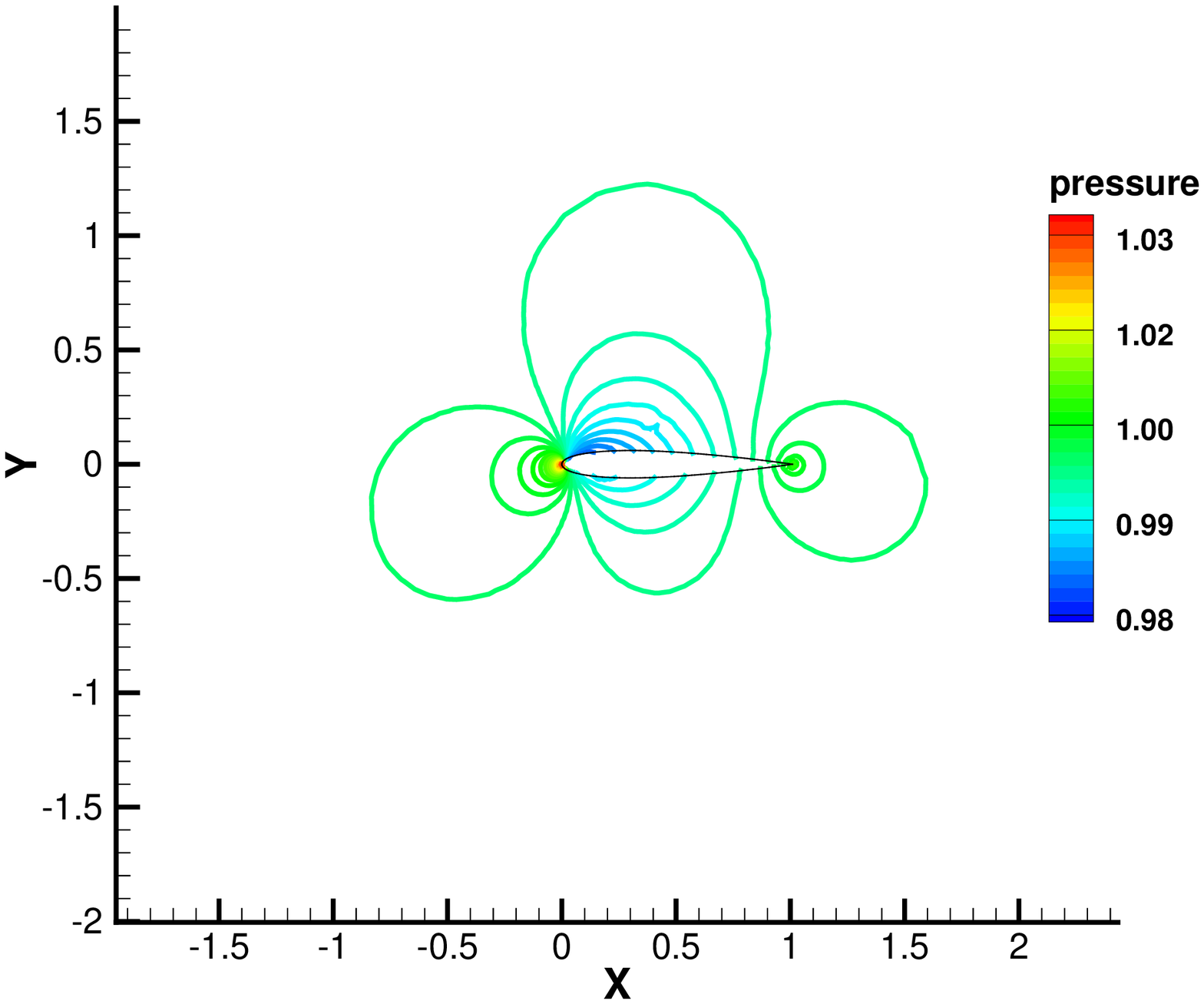}
%\caption{fig1}
\end{minipage}%
}%
\subfigure[FE fast sweeping, $Ma=0.2$]{
\begin{minipage}[t]{0.3\linewidth}
\centering
\includegraphics[width=2.0in]{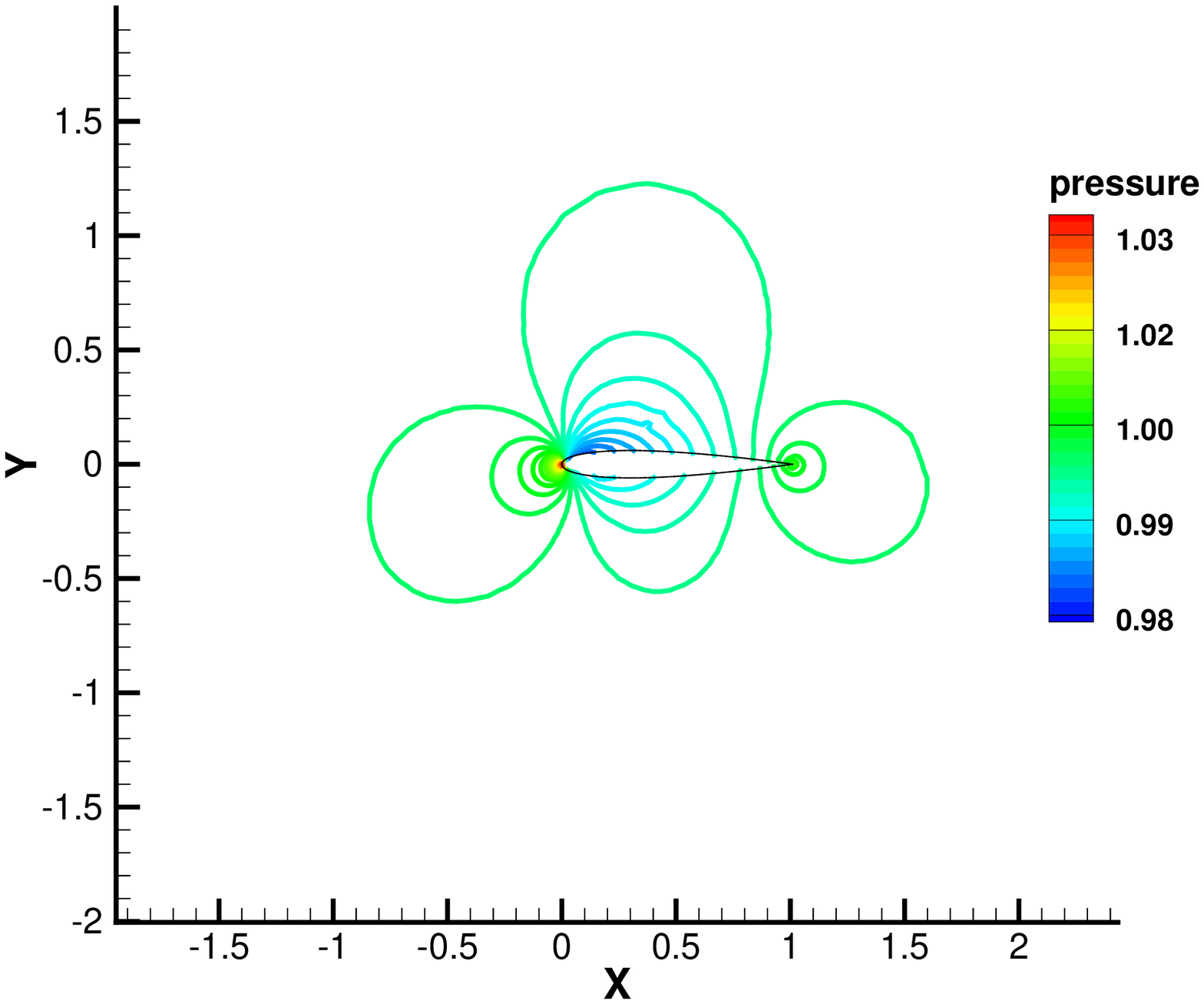}
%\caption{fig2}
\end{minipage}%
}%
\centering
\caption{\label{5.2}Example 7, subsonic flows past an NACA0012 airfoil. The converged steady states of numerical solutions by three different iterative schemes. (a) (b) (c): 30 equally spaced pressure contours from 0.5 to 1.46 for the case of  $Ma=0.8$, angle of attack $\alpha=1.25^{\circ}$; (d) (e) (f): 30 equally spaced pressure contours from 0.98 to 1.03 for the case of  $Ma=0.2$, angle of attack $\alpha=1^{\circ}$.}
\end{figure}

\begin{figure}%[H]
\centering
\subfigure[FE Jacobi, $Ma=0.8$]{
\begin{minipage}[t]{0.3\linewidth}
\centering
\includegraphics[width=2.0in]{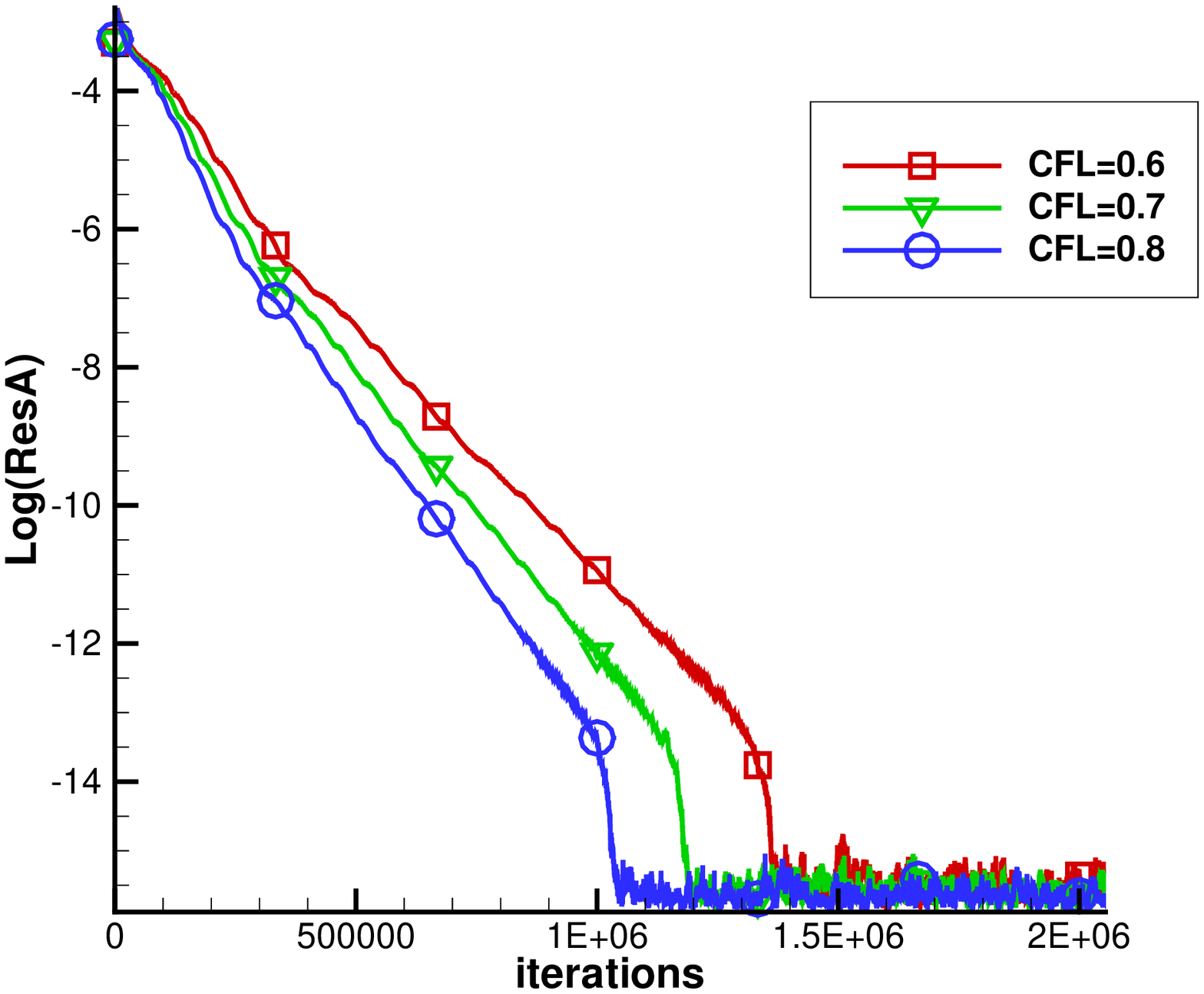}
%\caption{fig1}
\end{minipage}%
}%
\subfigure[RK Jacobi, $Ma=0.8$]{
\begin{minipage}[t]{0.3\linewidth}
\centering
\includegraphics[width=2.0in]{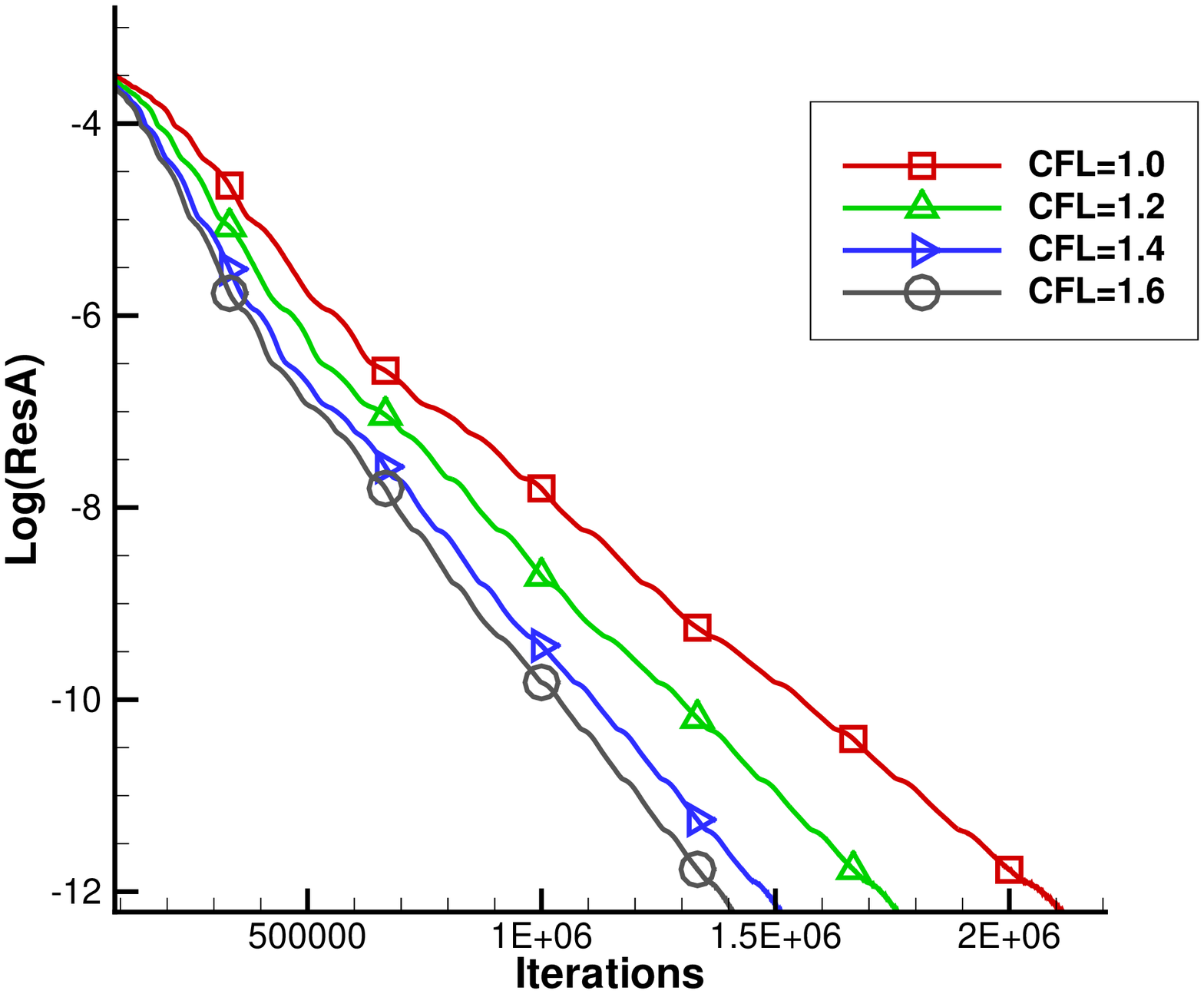}
%\caption{fig1}
\end{minipage}%
}%
\subfigure[FE fast sweeping, $Ma=0.8$]{
\begin{minipage}[t]{0.3\linewidth}
\centering
\includegraphics[width=2.0in]{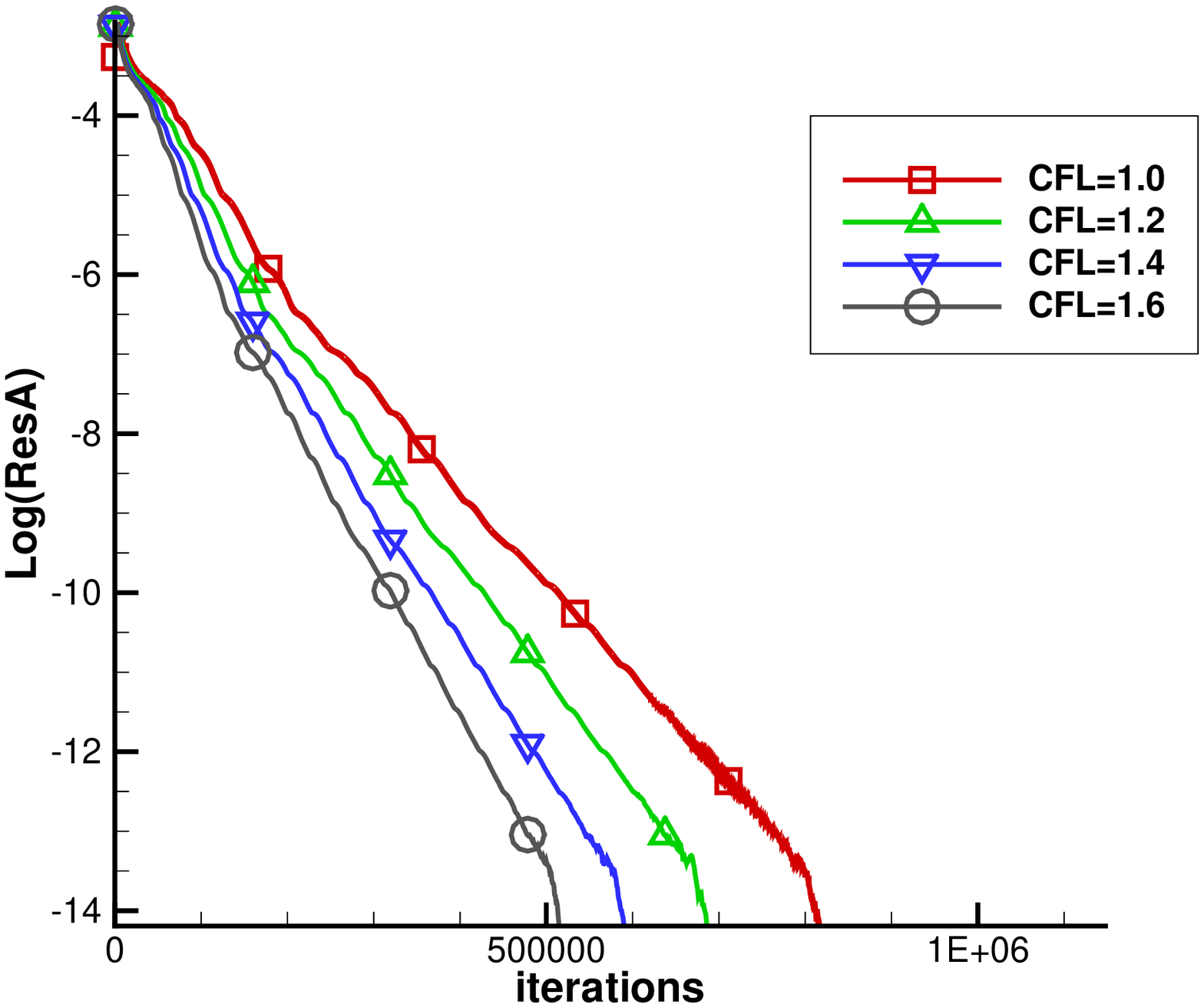}
%\caption{fig2}
\end{minipage}%
}%
\newline
\subfigure[FE Jacobi, $Ma=0.2$]{
\begin{minipage}[t]{0.3\linewidth}
\centering
\includegraphics[width=2.0in]{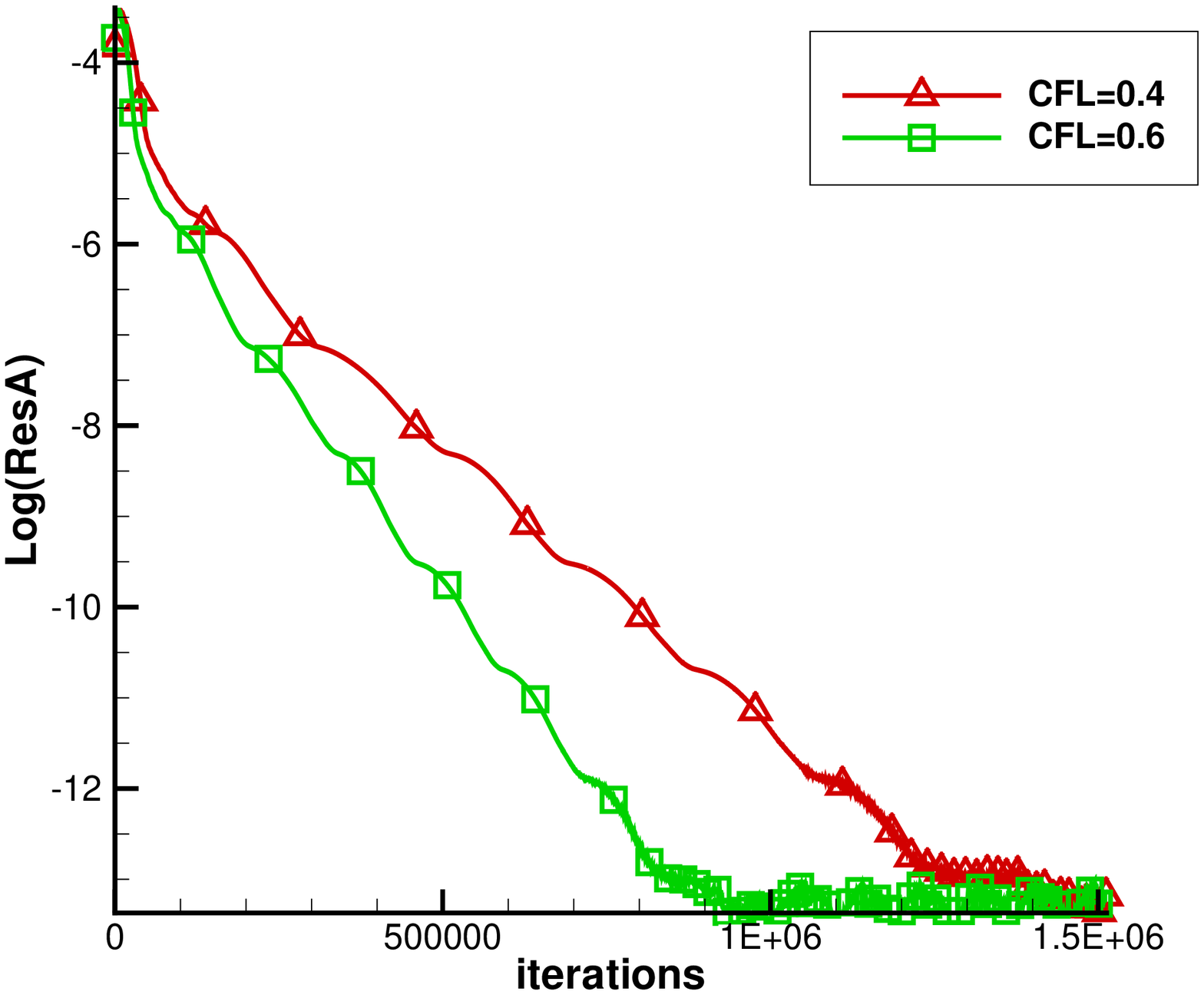}
%\caption{fig1}
\end{minipage}%
}%
\subfigure[RK Jacobi, $Ma=0.2$]{
\begin{minipage}[t]{0.3\linewidth}
\centering
\includegraphics[width=2.0in]{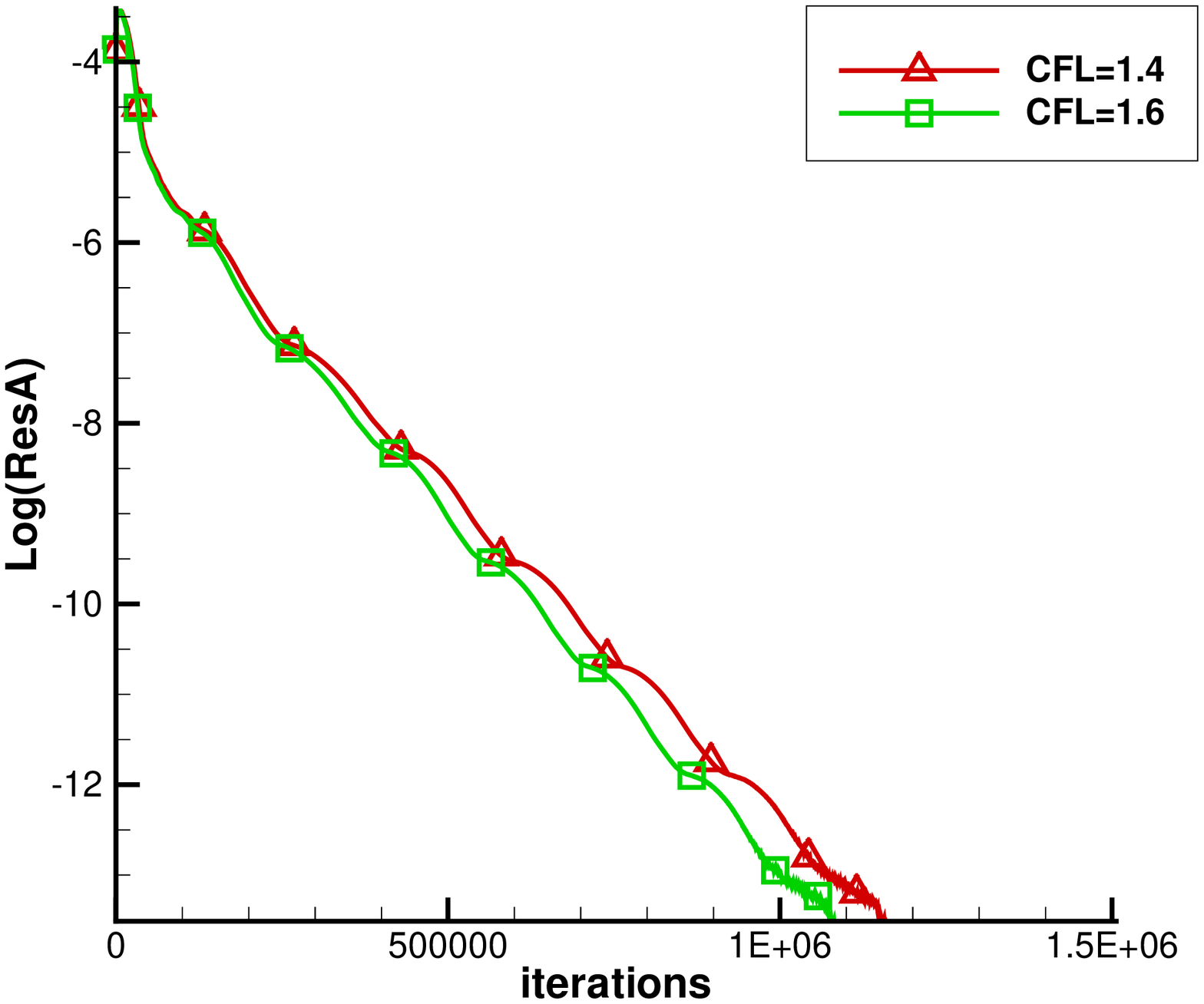}
%\caption{fig1}
\end{minipage}%
}%
\subfigure[FE fast sweeping, $Ma=0.2$]{
\begin{minipage}[t]{0.3\linewidth}
\centering
\includegraphics[width=2.0in]{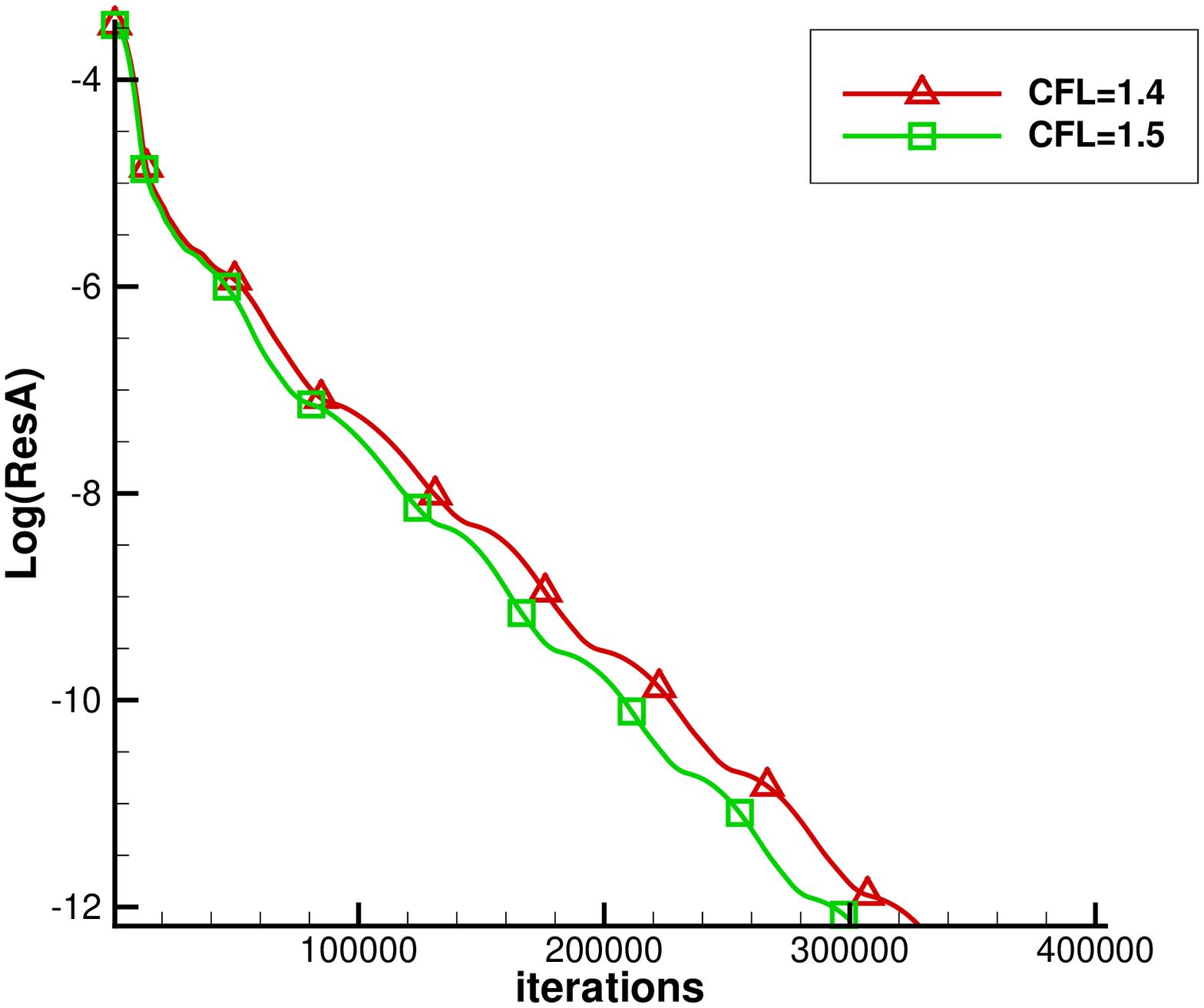}
%\caption{fig2}
\end{minipage}%
}%
\centering
\caption{\label{5.3}Example 7, subsonic flows past an NACA0012 airfoil. The convergence history of the residue as a function of number of iterations for three schemes with different CFL numbers. (a), (b), (c): for the case of $Ma=0.8$, angle of attack $\alpha=1.25^{\circ}$; (d) (e) (f): for the case of $Ma=0.2$, angle of attack $\alpha=1^{\circ}$.}
\end{figure}

\noindent As a continuation of Example 6, the problem of inviscid Euler subsonic flows past a single NACA0012 airfoil configuration in \cite{LBL} is solved in this example. Two cases of subsonic flows are considered, i.e., a flow with Mach number $Ma=0.8$, angle of attack $\alpha=1.25^{\circ}$; and a flow with Mach number $Ma=0.2$, angle of attack $\alpha=1^{\circ}$. The computational domain is $[-15,15]\times[-15,15]$. The unstructured mesh used is the same as that for Example 6, which is shown in Figure \ref{3.1}. Four corners of the computational domain are chosen as the reference points to form the alternating sweeping directions in the FE fast sweeping scheme.
In Table \ref{5.1}, number of iterations required to reach the convergence criterion threshold value $10^{-11}$,  and total CPU time when the schemes converge under various CFL numbers are reported for the FE Jacobi scheme, the RK Jacobi scheme, and the FE fast sweeping scheme. In this example, we also notice that the CFL number constraint for the FE Jacobi scheme to converge is not as severe as that  in the examples 1, 2 and 3, and its simple one-stage structure makes it be more efficient to converge to steady states than the RK Jacobi scheme for both cases with different Mach numbers. Again, the FE fast sweeping scheme is the most efficient iterative method among three methods. It allows much larger CFL numbers than the FE Jacobi scheme, and comparable CFL number sizes with the RK Jacobi scheme. The FE fast sweeping scheme also has a simple one-stage structure. With the largest CFL number permitted in each method to reach steady state solution, the FE fast sweeping method on unstructured triangular meshes saves more than $40\%$ CPU time cost of that by the RK Jacobi scheme (the TVD-RK3 scheme) for both subsonic flow cases in this example.
The contour plots of the pressure variable of the converged steady state solutions for these three schemes are shown in Figure \ref{5.2}. Comparable numerical steady states for these different iterative schemes are
observed. Figure \ref{5.3}
presents the residue history of these three schemes with different CFL numbers. It is observed that the residue of iterations can settle down to tiny values at the level of round off errors, which again verifies the absolute convergence of the developed high order fast sweeping method on triangular meshes.

\begin{table}%[H]
		\centering
\begin{tabular}{|c|c|c|}\hline

			\multicolumn{3}{|c|}{\bm{$Ma=0.8,\alpha=1.25^{\circ}$} }\\\hline
			\multicolumn{3}{|c|}{FE Jacobi scheme }\\\hline
             CFL number & iteration number  & CPU time \\\hline
	0.6& 1007029 	& 425623.98 \\\hline
0.7& 863213 	&  365297.07 \\\hline
    0.8   & 755334     & 316719.83\\\hline
      0.9   & Not convergent        & -\\
    \hline
		\end{tabular}

\begin{tabular}{|c|c|c|}\hline
			\multicolumn{3}{|c|}{RK Jacobi scheme }\\\hline
             CFL number & iteration number & CPU time \\\hline
	1.0&  1812610 	& 742083.62 \\\hline
    1.2   &   1510567    & 623277.44  \\\hline
     1.4   & 1294810      & 533356.95\\\hline
          1.6   & 1208500  & 500951.51\\\hline
           1.7   & Not convergent        & -\\
    \hline
		\end{tabular}

\begin{tabular}{|c|c|c|}\hline
			\multicolumn{3}{|c|}{FE fast sweeping scheme }\\\hline
             CFL number & iteration number  & CPU time \\\hline
     	1.0& 598216 	& 477581.42 \\\hline
    1.2   &   498168    &398820.98  \\\hline
     1.4   & 426688      & 338455.70\\\hline
      1.6   & 373080      & 295754.15\\\hline
       1.7   & Not convergent      & -\\
    \hline
		\end{tabular}

\begin{tabular}{|c|c|c|}\hline
			\multicolumn{3}{|c|}{\bm{$Ma=0.2,\alpha=1^{\circ}$} }\\\hline
			\multicolumn{3}{|c|}{FE Jacobi scheme }\\\hline
             CFL number & iteration number  & CPU time \\\hline
	0.4& 959981 & 392068.22 \\\hline
	0.6& 639947	& 273666.77 \\\hline
      0.7   & Not convergent        & -\\
    \hline
		\end{tabular}

\begin{tabular}{|c|c|c|}\hline
			\multicolumn{3}{|c|}{RK Jacobi scheme }\\\hline
             CFL number & iteration number & CPU time \\\hline
     	1.4& 822721 	& 331972.06 \\\hline
     1.6& 767878 	& 317268.26 \\\hline
           1.7   & Not convergent        & -\\
    \hline
		\end{tabular}

\begin{tabular}{|c|c|c|}\hline
			\multicolumn{3}{|c|}{FE fast sweeping scheme }\\\hline
             CFL number & iteration number  & CPU time \\\hline
     	1.4&   274184 &  203245.83 \\\hline
       1.5&   252120 	& 186856.24\\\hline
       1.6   & Not convergent      & -\\
    \hline
		\end{tabular}
		\caption{\label{5.1}Example 7, subsonic flows past an NACA0012 airfoil. Number of iterations and total CPU time when convergence is obtained. Convergence criterion threshold value   is $10^{-11}$. CPU time unit: second.}
	\end{table}

\section{Concluding remarks}
High order accuracy fast sweeping methods have been well developed on rectangular meshes to efficiently solve steady state solutions of hyperbolic PDEs. However, it was still a open problem how to design high order accuracy fast sweeping methods on unstructured meshes. In this paper, we develop a high order fixed-point fast sweeping method on unstructured triangular meshes for solving steady state solutions of hyperbolic conservation laws.
Multiple reference points on the computational domain are introduced to order all the
cells and form alternating sweeping directions on unstructured meshes. The local solver of the proposed fixed-point fast sweeping method is based on a fifth-order finite volume unstructured WENO scheme with unequal-sized sub-stencils, to achieve the absolute convergence of the iterations. Extensive numerical experiments, including solving difficult problems which are defined on complex domains or challenging for high order schemes to converge to steady states, show that the designed fixed-point fast sweeping scheme on unstructured meshes can significantly enlarge the CFL number of the forward Euler scheme with a high order WENO spatial discretization (e.g., the fifth order WENO scheme here) to the level of the TVD-RK3 scheme.
As a result, up to $70\%$ CPU time can be saved by using the proposed fast sweeping method rather than the TVD-RK3 scheme for iterations to converge to steady states of the WENO scheme. The iteration residues of the new absolutely convergent fast sweeping method on unstructured meshes converge to round off errors for all benchmark problems tested in this paper.

\end{document}